\documentclass[a4paper,12pt,titlepage]{article}
\usepackage{latexsym}
\usepackage{amsbsy}
\usepackage{amssymb}
\usepackage{amscd}
\usepackage{amsmath}
\begin{document}
\newcommand{\ep}{\hspace*{\fill}$\Box$}
\newcommand{\eps}{\varepsilon}
\newcommand{\pr}{{\bf Proof. }}
\newcommand{\ms}{\medskip\\}
\newcommand{\cl}{\mbox{\rm cl}}
\newcommand{\g}{\ensuremath{\mathfrak g} }
\newcommand{\gc}{${\cal G}$-complete }
\newcommand{\sa}{\stackrel{\scriptstyle s}{\approx}}
\newcommand{\prol}{\mbox{\rm pr}^{(n)}}
\newcommand{\prolo}{\mbox{\rm pr}^{(1)}}
\newcommand{\deta}{\frac{d}{d \eta}{\Big\vert}_{_{0}}}
\newcommand{\detas}{\frac{d}{d \eta}{\big\vert}_{_{0}}}
\newcommand{\R}{\mathbb R}
\newcommand{\N}{\mathbb N}
\newcommand{\C}{\mathbb C}
\newcommand{\Z}{\mathbb Z}
\newcommand{\K}{\mathbb K}
\newcommand{\sR}{\mathbb R}
\newcommand{\sN}{\mathbb N}
\newcommand{\gK}{{\cal K}}
\newcommand{\gR}{{\cal R}}
\newcommand{\gC}{{\cal C}}
\newcommand{\Dp}{${\cal D}'$ }                                          
\newcommand{\go}{${\cal G}(\Omega)$ }
\newcommand{\grn}{${\cal G}(\R^n)$ }
\newcommand{\grp}{${\cal G}(\R^p)$ }
\newcommand{\grq}{${\cal G}(\R^q)$ }
\newcommand{\gt}{${\cal G}_\tau$ }
\newcommand{\gto}{${\cal G}_\tau(\Omega)$ }
\newcommand{\gtrn}{${\cal G}_\tau(\R^n)$ }
\newcommand{\gtrp}{${\cal G}_\tau(\R^p)$ }
\newcommand{\gtrq}{${\cal G}_\tau(\R^q)$ }
\newcommand{\gtn}{{\cal G}_\tau(\R^n) }
\newtheorem{thr}{\hspace*{-1.1mm}}[section]
\newcommand{\bt}{\begin{thr} {\bf Theorem. }}
\newcommand{\et}{\end{thr}}
\newcommand{\bp}{\begin{thr} {\bf Proposition. }}
\newcommand{\bc}{\begin{thr} {\bf Corollary. }}
\newcommand{\blem}{\begin{thr} {\bf Lemma. }}
\newcommand{\bex}{\begin{thr} {\bf Example. }\rm}
\newcommand{\bexs}{\begin{thr} {\bf Examples. }\rm}
\newcommand{\bd}{\begin{thr} {\bf Definition. }}
\newcommand{\beast}{\begin{eqnarray*}}
\newcommand{\eeast}{\end{eqnarray*}}
\newcommand{\wsc}[1]{\overline{#1}^{wsc}}
\newcommand{\todo}[1]{$\clubsuit$\ {\tt #1}\ $\clubsuit$}
\newcommand{\rem}[1]{\vadjust{\rlap{\kern\hsize\thinspace\vbox%
                       to0pt{\hbox{${}_\clubsuit${\small\tt #1}}\vss}}}}
\newcommand{\ahat}{\ensuremath{\hat{\mathcal{A}}_0(M)} }
\newcommand{\atil}{\ensuremath{\tilde{\mathcal{A}}_0(M)} }
\newcommand{\aqtil}{\ensuremath{\tilde{\mathcal{A}}_q(M)} } 
\newcommand{\ehat}{\ensuremath{\hat{\mathcal{E}}(M)} } 
\newcommand{\emhat}{\ensuremath{\hat{\mathcal{E}}_m(M)} }
\newcommand{\nhat}{\ensuremath{\hat{\mathcal{N}}(M)} }
\newcommand{\ghat}{\ensuremath{\hat{\mathcal{G}}(M)} } 
\newcommand{\lhat}{\ensuremath{\hat{L}_X} }                    
\newcommand{\comp}{\subset\subset}
\newcommand{\al}{\alpha}
\newcommand{\bet}{\beta} 
\newcommand{\ga}{\gamma}
\newcommand{\Om}{\Omega}\newcommand{\Ga}{\Gamma}\newcommand{\om}{\omega}
\newcommand{\si}{\sigma}\newcommand{\la}{\lambda}
\newcommand{\de}{\delta}
\newcommand{\vphi}{\varphi}\newcommand{\dl}{{\displaystyle \lim_{\eta>0}}\,}
\newcommand{\intl}{\int\limits}\newcommand{\su}{\sum\limits_{i=1}^2}
\newcommand{\D}{{\cal D}}\newcommand{\Vol}{\mbox{Vol\,}}
\newcommand{\Or}{\mbox{Or}}\newcommand{\sign}{\mbox{sign}}
\newcommand{\na}{\nabla}\newcommand{\pa}{\partial}
\newcommand{\ti}{\tilde}\newcommand{\T}{{\cal T}} \newcommand{\G}{{\cal G}}
\newcommand{\DD}{{\cal D}}\newcommand{\X}{{\cal X}}\newcommand{\E}{{\cal E}} 
\newcommand{\CC}{{\cal C}}\newcommand{\vo}{\Vol}
\newcommand{\bat}{\bar t}
\newcommand{\bx}{\bar x}
\newcommand{\by}{\bar y} \newcommand{\bz}{\bar z}\newcommand{\br}{\bar r}
\newcommand{\fr}{\frac{1}}\newcommand{\il}{\int\limits}
\newcommand{\nn}{\nonumber}
\newcommand{\supp}{\mathop{\mathrm{supp}}}

\newcommand{\vp}{\mbox{vp}\frac{1}{x}}\newcommand{\A}{{\cal A}}
\newcommand{\Ll}{L_{\mbox{\small loc}}}\newcommand{\Hl}{H_{\mbox{\small loc}}}
\newcommand{\Lll}{L_{\mbox{\scriptsize loc}}}
\newcommand{\be}{ \begin{equation} }\newcommand{\ee}{\end{equation} }
\newcommand{\beq}{ \begin{equation} }\newcommand{\eeq}{\end{equation} }
\newcommand{\bea}{\begin{eqnarray}}\newcommand{\eea}{\end{eqnarray}}
\newcommand{\beas}{\begin{eqnarray*}}\newcommand{\eeas}{\end{eqnarray*}}
\newcommand{\beqs}{\begin{equation*}}\newcommand{\eeqs}{\end{equation*}}
\newcommand{\lb}{\label}\newcommand{\rf}{\ref}
\newcommand{\GL}{\mbox{GL}}\newcommand{\bfs}{\boldsymbol}
\newcommand{\ben}{\begin{enumerate}}\newcommand{\een}{\end{enumerate}}
\newcommand{\ba}{\begin{array}}\newcommand{\ea}{\end{array}}
\newtheorem{thi}{\hspace*{-1.1mm}}[section]
\newcommand{\bthm}{\begin{thr} {\bf Theorem. }}
\newcommand{\bprop}{\begin{thr} {\bf Proposition. }}
\newcommand{\bcor}{\begin{thr} {\bf Corollary. }}
\newcommand{\bdef}{\begin{thr} {\bf Definition. }}
\newcommand{\brem}{\begin{thr} {\bf Remark. }\rm}
\newcommand{\bth}{\begin{thr}\rm}
\newcommand{\ethi}{\end{thr}}
\newcommand{\ca}{{\cal A}}
\newcommand{\cb}{{\cal B}}
\newcommand{\cc}{{\cal C}}
\newcommand{\cd}{{\cal D}}
\newcommand{\ce}{{\cal E}}
\newcommand{\cg}{{\cal G}}
\newcommand{\ci}{{\cal I}}
\newcommand{\cn}{{\cal N}}
\newcommand{\cs}{{\cal S}}
\newcommand{\ct}{{\cal T}}
\newcommand{\rmd}{\mbox{\rm d}}
\newcommand{\io}{\iota}
\newcommand{\bnot}{\begin{thr} {\bf Notation }}
\newcommand{\lgl}{\langle}
\newcommand{\rgl}{\rangle}
\newcommand{\spp}{\mbox{\rm supp\,}}
\newcommand{\id}{\mathop{\mathrm{id}}}
\newcommand{\pro}{\mathop{\mathrm{pr}}}
\newcommand{\dist}{\mathop{\mathrm{dist}}}
\newcommand{\clb}{\overline{B}_}
\newcommand{\sgn}{\mathop{\mathrm{sgn}}}

\parskip=2mm

\begin{center}
{\bf \Large On the foundations of nonlinear generalized functions I}
\vskip2mm
{\large        E. Farkas, M. Grosser, M. Kunzinger and R. Steinbauer}

        Universit\"at Wien\\ 
        Institut f\"ur Mathematik
\end{center}

{\small
{\sc Abstract.} We construct a diffeomorphism invariant (Colombeau-type) 
differential algebra
ca\-no\-ni\-cally containing the space of distributions in the sense of L.
Schwartz. Employing differential calculus in infinite dimensional
(convenient) vector spaces, previous attempts in this direction are unified
and completed. Several classification results are achieved and applications
to nonlinear differential equations involving singularities are given.


}

\section{Introduction}\lb{intro}
In his celebrated impossibility result (\cite{Schw}), L. Schwartz demonstrated
that the space ${\cal D}'(\Omega)$ of distributions over
some open subset $\Omega$ of $\R^n$ cannot be embedded into an associative commutative 
algebra $({\cal A}(\Omega),+,\circ)$ satisfying
\begin{itemize}
\item[(i)] ${\cal D}'(\Omega)$ is linearly embedded into ${\cal A}(\Omega)$ 
and $f(x)\equiv 1$ is the unity in ${\cal A}(\Omega)$. 
\item[(ii)] There exist derivation operators $\partial_i:{\cal A}(\Omega)
\rightarrow {\cal A}(\Omega)$ ($i=1,\ldots,n$) that are linear and satisfy the 
Leibnitz rule.
\item[(iii)] $\partial_i|_{{\cal D}'(\Omega)}$ is the usual partial 
derivative ($i=1,\ldots,n$).
\item[(iv)] $\circ|_{{\cal C}(\Omega)\times {\cal C}(\Omega)}$
coincides with the pointwise product of functions. 
\end{itemize}
Since this result remains valid upon replacing $\cc(\Om)$ by $\cc^k(\Om)$ for 
any finite $k$, the best possible result would consist in constructing
an embedding of $\D'(\Om)$ as above with (iv) replaced by
\begin{itemize}
\item[(iv')] $\circ|_{{\cal C}^\infty(\Omega)\times {\cal C}^\infty(\Omega)}$
coincides with the pointwise product of functions. 
\end{itemize}
The actual construction of differential algebras satisfying these optimal 
properties is due to J.~F. Colombeau (\cite{c1}, \cite{c2},
\cite{cbull}, \cite{clec}).
The need for algebras of this type arises, for example, from the
necessity of considering non-linear PDEs where either the respective 
coefficients, the data or the prospective solutions are non-smooth. 
Classical linear distribution theory 
clearly
does not permit the treatment of such problems. Colombeau algebras, 
on the other hand, have proven to be a useful tool for analyzing such 
questions (for applications in nonlinear PDEs, cf. e.g., \cite{bmo}, 
\cite{bmo2}, \cite{cmo}, \cite{chmo}, \cite{chmo2}, \cite{kk}, 
\cite{MObook}, for applications to numerics, see e.g., \cite{b},
\cite{bc}, \cite{bcm}, for applications in mathematical physics, 
e.g., \cite{v}, \cite{KS}, \cite{gmot} as well as the literature
cited in these works). For alternative approaches to algebras of 
generalized functions, cf. \cite{r1}, \cite{r2}.

Since Colombeau's monograph \cite{c1}, there have been introduced a considerable
number of variants of Colombeau algebras, many of them adapted to
special purposes. From the beginning, however, the question of the
functor property of the construction was at hand as a crucial one:
If $\mu:\ti\Om\to\Om$ denotes
a diffeomorphism between open subsets $\ti\Om,\Om$ of $\R^s$, is
it possible to extend the operation $\mu^*:f\mapsto f\circ\mu$ on
smooth distributions on $\Om$ to an operation $\hat\mu$ on
the Colombeau algebra such that
$(\mu\circ\nu)\hat{\ }=\hat\nu\circ\hat\mu$ and $(\mbox{\rm id})\hat{\ }=
\mbox{\rm id}$ are satisfied?
To phrase it differently, is it possible to achieve a diffeomorphism
invariant construction of Colombeau algebras? As long as this question
could not be answered in the positive, there remained the serious objection
that there is no way of defining such algebras on manifolds, based
on intrinsic terms, exclusively. 
(This  discussion 
does not take into account the so called  ``special''
or  ``simplified'' variant of Colombeau's algebra whose elements are classes of
nets of smooth functions indexed by $\eps>0$ (cf. \cite{MObook},
p.~109). Although diffeomorphism-invariant, these algebras lack  
a canonical embedding of distributions (\cite{dd}, \cite{roldiss}), so
we do not consider them here.)

The first variants of Colombeau algebras, though serving as a
valuable
tool in the treatment of non-linear problems, indeed did
not have the property of diffeomorphism invariance: Some of the
key ingredients used in defining them (in particular, the ``test
objects'' (see section \rf{scheme1}) being employed as well
as the definition of the subsets 
$\ca_q(\R^s)$ of the set of all test objects)
turned out not to be invariant under the natural action of a
diffeomorphism.

Colombeau and Meril in their paper \cite{CM} made the first decisive
steps to remove this flaw by proposing a construction of
Colombeau algebras which they claimed to be diffeomorphism invariant.
As an essential tool, they had to use calculus on locally convex spaces.
However, they did not give the details of
the application of that calculus; moreover, their definition of the
objects constituting the Colombeau algebra was not unambiguous
and, which amounts to the most serious objection, their notion of test
objects still was not preserved under the action of the 
diffeomorphism. Nevertheless, despite these defects
(which, apparently, went unnoticed by nearly all workers in the
field) their construction was quoted and used many times 
(see, e.g., \cite{kunzdiss}, \cite{tschapitsch}, \cite{we}, \cite{w}, \cite{vw1},
\cite{serbenbuch}, \cite{vw2}, \cite{trio}, \cite{v}, \cite{hb},
\cite{KS}).
It was only in 1998 that J. Jel\'\i nek in \cite{JEL} pointed
out the error in \cite{CM} by giving a (rather
simple) counterexample. In the same paper, he presented another version
of the theory avoiding the shortcomings of \cite{CM} and forming the basis for the 
approach taken here.

The present article is the first in a series of two papers. It is organized as follows:
After fixing notation and terminology in section
\rf{notterm}, a general scheme of construction for diffeomorphism-invariant 
Colombeau-type algebras of generalized functions is introduced
in section \rf{scheme1}. Section \rf{calc} gives a quick
overview of calculus in convenient vector spaces providing the necessary 
results for the development of the theory. Especially with a view to applications
(in particular: partial differential equations) we feel that this approach has several
advantages over the concept of Silva-differentiability employed so far.
Section \rf{basics} introduces a translation formalism that allows to 
freely  switch  between what we call the C- and J- ({\em Colombeau}- and
{\em  Jelinek}-)  formalism  of  the  diffeomorphism invariant theory to be
constructed  in  section  \rf{jelshort}.  In  the  actual construction of this algebra,
smooth functions defined on sets denoted by
$U_\eps(\Om)$ play a central r\^ole. Differentials of such functions
are of utmost importance in the development of the theory. However,
$U_\eps(\Om)$ is not a linear space. Sections \rf{calcuepsom} thus provides the
framework necessary for doing calculus on $U_\eps(\Om)$. 
A complete presentation of 
the resulting diffeomorphism invariant algebra,
based on the general construction scheme of section
\rf{scheme1}, is the focus of section \rf{jelshort}.
The sheaf-theoretic properties of this algebra are discussed in
section \rf{sheaves}.
This is followed by a short 
section on the separation of testing procedures and definition of objects in algebras
of generalized functions. Section \rf{atoz} provides several new
characterizations of the fundamental building blocks ${\mathcal E}_M$ and 
${\mathcal N}$ of the algebra. In particular, 
these characterizations will constitute the key ingredient in
obtaining an intrinsic description of the theory on manifolds (\cite{vi}). 
Finally, we present
some applications to partial differential equations in section
\rf{diffeq}.
                
The second paper of this series
gives a comprehensive analysis of algebras of
Co\-lom\-beau-type generalized functions in the range between the 
diffeomorphism-invariant quotient algebra $\mathcal{G}^d=
\mathcal{E}_M\big/\mathcal{N}$ introduced in section \ref{jelshort}
and  (the smooth version of)
Colombeau's original algebra
$\cg^e$ introduced in \cite{c2} (which, to be sure, is the standard
version among those being independent of the choice of a particular
approximation of the delta distribution).
Three main results
are established: First, a simple criterion describing membership in
$\mathcal{N}$ (applicable to all types of Colombeau algebras) is given
(section 13).
Second, two counterexamples demonstrate that $\mathcal{G}^d$ is not
injectively included in $\mathcal{G}^e$ (section 15); their
construction is based on
a completeness theorem for spaces of smooth functions
in the sense of sections \ref{calc} and \ref{calcuepsom}
(section 14). Finally, it is shown that in
the range  ``between'' $\mathcal{G}^d$ and $\mathcal{G}^e$ only one more
construction leads to a diffeomorphism invariant algebra. In analyzing the
latter, several classification results essential
for obtaining an intrinsic description of $\mathcal{G}^d$ on manifolds
are derived (sections 16, 17).
The concluding section 18
points out that also weaker invariance properties than with respect
to all diffeomorphisms should be envisaged for Colombeau algebras,
in particular regarding applications.
\section{Notation and Terminology}\lb{notterm}
Throughout this paper, $\Om$, $\tilde \Om$ will denote non-empty open subsets
of $\R^s$. For any $A\subseteq \R^s$,
$A^\circ$ denotes its interior. $\mathcal{C}^\infty(\Om)$ is the space
of smooth, complex valued functions on $\Om$. If $f\in \mathcal{C}^\infty(\Om)$ then $Df$
denotes its (total) derivative. Also, we set $\check f(x) = f(-x)$. On any cartesian product, 
$\mathrm{pr_i}$ denotes the projection onto the $i$-th factor. For $r\in \R$, $[r]$ is the largest
integer $\le r$. We set $I = (0,1]$. 
Concerning locally convex spaces our basic reference is 
\cite{schae}. In particular, by a locally convex space we mean a vector space endowed with a
locally convex Hausdorff topology. The space of test functions
(i.e., compactly supported smooth functions) on $\Om$ is denoted by $\D(\Om)$ and is equipped with
its natural (LF)-topology; its dual, the
space of distributions on $\Om$ is termed $\D'(\Om)$. The action of any $u\in \D'(\Om)$ on a test
function $\vphi$ will be written as $\langle u,\vphi\rangle$. $\de$ denotes the Dirac delta
distribution. $K\comp A$ ($A\subseteq \R^s$) means that $K$ is a compact subset of $A^\circ$. 
For $K\comp \Om$, $\D_K(\Om)$ is the space of smooth functions on $\Om$ 
supported in $K$. We set
\beas
&& \A_0(\Om) = \{\vphi \in \D(\Om)| \int\! \vphi(\xi)\,d\xi = 1\}\\
&& \A_q(\Om) = \{\vphi \in \A_0(\Om)| \int\! \xi^\al\vphi(\xi)\,d\xi = 0, 
  \ 1\le |\al| \le q,\ \al \in \N_0^s \} \qquad (q \in \N)
\eeas
$\A_{q0}(\Om)$ is the linear subspace of $\D(\Om)$ parallel to the affine space 
$\A_q(\Om)$ ($q\in \N_0$).
For any maps $f,g,h$ such that $g\circ f$ and $f\circ h$ are defined
we set $g_*(f):=g\circ f$ and $h^*(f):=f\circ h$.
$\pa_i$ resp.\ $\pa^\al$ always stand for $\frac{\pa}{\pa x_i}$ resp.\
$\frac{\pa^{|\al|}}{\pa x_i^\al}$ ($\al\in\N_0^s$).

For any locally convex space $F$ the space $\cc^\infty(\Om,F)$ of smooth functions
from $\Om$ into $F$ will always carry the topology of uniform
convergence in all derivatives on compact subsets of $\Om$.
In particular, a subset ${\cal B}$ of this space will be said to be
bounded if, for any $K\subset\subset\Om$ and any
$\al\in\N_0^s$, the set $\{\pa^\al\phi(x)\mid\phi\in{\cal B},x\in K\}$
is bounded in $F$. Observe that in case that the image of a map
$\phi\in\cc^\infty(\Om,F)$ is contained in some affine subspace
$F_0$ of $F$ then the derivatives of $\phi$ take their values in
the linear subspace parallel to $F_0$. For locally convex spaces $E$, $F$
the space $\cc^\infty(E,F)$ (resp. $\cc^\infty(E)$ for $F = \C$) is introduced in 
section \ref{calc}.

In what follows, $\ca_0(\R^s)$ may be replaced by any closed affine
subspace of $\cd(\R^s)$. By $\cc^{[\infty,\Om]}_b(I\times\Om,\ca_0(\R^s))$
we denote the space of all maps
$\phi:I\times\Om\to\ca_0(\R^s)$ which are smooth with 
respect to the second argument and
bounded in the sense that the corresponding map
$\hat\phi:I\to\cc^\infty(\Om,\ca_0(\R^s))$ has a bounded image
as defined above, i.e., for every $K\subset\subset\Om$
and any $\al\in\N_0^s$, the set 
$\{\pa^\al(\hat\phi(\eps))(x)\mid\eps\in I,x\in K\}$
is bounded in $\ca_0(\R^s)$ resp.\ $\cd(\R^s)$, which, in turn,
is equivalent to saying that
\begin{enumerate}
\item
for every $K$ as above and any
$\al\in\N_0^s$, the supports of all $\pa_x^\al\phi(\eps,x)$
($\eps\in I$, $x\in K$) are contained in some fixed bounded
set (depending only on $K$) and 
\item
$\sup\{|\pa_\xi^\bet(\pa_x^\al(\hat\phi(\eps))(x))(\xi)|
\eps\in I,\ x\in K,\ \xi\in\R^s\}$ (or, expressed in terms of $\phi$
itself)
$\sup\{|\pa_\xi^\bet\pa_x^\al(\phi(\eps,x))(\xi)|
\eps\in I,\ x\in K,\ \xi\in\R^s\}$
is finite.
\end{enumerate}
$\cc^\infty_b( I\times\Om,\ca_0(\R^s))$ is the subspace of  $\cc^{[\infty,\Om]}_b( I
\times\Om,\ca_0(\R^s))$ whose elements are smooth in both arguments. 
Finally, for any $K \comp \Om$ and any $q\ge 1$ 
an element $\phi$ of $\cc_b^\infty( I\times\Om,
\cd(\R^s))$ is said to have asymptotically vanishing moments of order $q$ on $K$ if
\[
\sup_{x\in K} |\int\! \xi^\al\phi(\eps,x)(\xi)\,d\xi| = O(\eps^{q}) \quad  
(1 \le |\al| \le q)\,.
\]
For this notion to make sense it is obviously sufficient for $\phi$
to be defined on $(0,\eps_0]\times K$ for some $\eps_0>0$.

\section{Scheme of construction}\lb{scheme1}
As was already pointed out in section 1, due to the lack of a canonical
embedding  of the space of distributions into ``special''
variants  of  Colombeau algebras we shall not consider these. 
Instead, we focus on  ``full'' algebras 
(in  the  sense  of \cite{kunzdiss}, p.~31), distinguished by the fact that
such  a canonical embedding is always available. Elements of full Colombeau
algebras  are  equivalence  classes  of  functions  $R$ taking as arguments
certain pairs $(\vphi,x)$ consisting of a suitable test function
$\vphi\in\cd(\R^s)$ and a point $x$ of $\Om$.

Every (full) Colombeau algebra is constructed according to the following
blueprint (where {\bf (Di)}, {\bf (Tj)} stand for Definition i and
Theorem j, respectively). {\bf (D5)} and
{\bf (T6)}--{\bf (T8)} are only relevant
if a diffeomorphism invariant type of algebra is to be obtained.

\vskip2mm
{\bf (D1)}\qquad $\E(\Om)$ (the ``basic space'', see the remarks
            below);
            maps $\si:\cc^\infty(\Om)\to\E(\Om)$,\linebreak
\hphantom{{\bf (D1)}\qquad }$\io:\cd'(\Om)\to\E(\Om)$.

{\bf (D2)}\qquad Derivations $D_i$ on $\E(\Om)$ ($i=1,\dots,s$)
                  extending the operators
                  $\frac{\pa}{\pa x_i}$ of partial\linebreak
\hphantom{{\bf (D2)}\qquad }differentiation on $\cd'(\Om)$ resp.\
                on $\CC^\infty(\Om)$, i.e.,\
                  $D_i\circ\io=\io\circ\frac{\pa}{\pa x_i}$ and \linebreak
\hphantom{{\bf (D2)}\qquad }$D_i\circ\sigma=\sigma\circ\frac{\pa}{\pa x_i}$.

{\bf (D3)}\qquad $\E_M(\Om)$ ($\subseteq\E(\Om)$;
the subspace of ``moderate'' functions).

{\bf (D4)}\qquad $\cn(\Om)$ ($\subseteq\E(\Om)$;
the subspace of ``negligible'' functions).

{\bf (T1)\,}\qquad $\io(\cd'(\Om))\subseteq\ce_M(\Om),\quad
  \si(\cc^\infty(\Om))\subseteq\ce_M(\Om),\quad
  (\io-\si)(\cc^\infty(\Om))\subseteq\cn(\Om)$;

\hspace*{-3pt}
\hphantom{{\bf (T1)\,}\qquad }$\io(\cd'(\Om))\cap\cn(\Om)=\{0\}$.

{\bf (T2)\,}\qquad $\E_M(\Om)$ is a subalgebra of $\E(\Om)$.

{\bf (T3)\,}\qquad $\cn(\Om)$ is an ideal in $\E_M(\Om)$.

{\bf (T4)\,}\qquad $\E_M(\Om)$ is invariant under each $D_i$.

{\bf (T5)\,}\qquad $\cn(\Om)$ is invariant under each $D_i$.

{\bf (D5)}\qquad For each diffeomorphism $\mu:\ti\Om\to\Om$,
                  a map $\bar\mu:D_{\ti\Om}\to D_\Om$~\footnote{Concerning $D_{\ti\Om},D_\Om$,
                  see the remark on {\bf (D1)} below.} is defined in
                  a\linebreak
\hphantom{{\bf (D5)}\qquad }functorial way such that its ``transpose''
       $\hat\mu:\ce(\Om)\to\ce(\ti\Om)$, $\hat\mu(R):=R\circ\bar\mu$,
       \linebreak
\hphantom{{\bf (D5)}\qquad }extends the usual effect
$\mu$ has on distributions, i.e., 
$\hat\mu\circ\io=\io\circ\mu^*$ where\linebreak
\hphantom{{\bf (D5)}\qquad }for $u\in\cd'(\Om)$, $\mu^*u$
                  is defined by 
$\lgl\mu^*u,\vphi\rgl:=
                  \lgl u,(\vphi\circ\mu^{-1})\cdot|\det
                  D\mu^{-1}|\rgl$.\linebreak
\hphantom{{\bf (D5)}\qquad }Similarly, we require $\hat\mu\circ\sigma=\sigma\circ\mu^*$
                on $\CC^\infty(\Om)$.

{\bf (T6)}\qquad The class of ``scaled test objects '' (see below)
                  is invariant under the action
\hphantom{{\bf (T6)}\qquad }induced by $\mu$.

{\bf (T7)\,}\qquad $\E_M$ is invariant under $\hat\mu$,
                     i.e., $\hat\mu$ maps
                     $\E_M(\Om)$ into $\E_M(\ti\Om)$.

{\bf (T8)\,}\qquad  $\cn$ is invariant under $\hat\mu$,
                     i.e., $\hat\mu$ maps
                     $\cn(\Om)$ into $\cn(\ti\Om)$.

{\bf (D6)}\qquad $\cg(\Om):=\E_M(\Om)\big/\cn(\Om)$.

For $R\in\E_M(\Om)$, the class $R+\cn(\Om)$ of $R$ in $\G(\Om)$ will be denoted
by $[R]$.

The following comments are intended to motivate and clarify 
the preceding---admittedly very formal---definition
schemes and theorems.

\vskip2mm {\bf ad (D1):}
Here, $\E(\Om)$ denotes some algebra of complex-valued functions
having appropriate smoothness properties on a suitable
domain $D_\Om\subseteq\cd(\R^s)\times\Om$.
$\si$ has to be an injective algebra homomorphism, whereas $\io$
just has to be linear and injective.

{\bf ad (D3), (D4):}
Membership of $R\in\E(\Om)$ in $\cn(\Om)$ respectively $\E_M(\Om)$ 
depends on 
the ``asymptotic'' behaviour of 
$R$ on certain paths in $\cd(\R^s)\times\Om$, where the second component
is constant whereas the
first component, depending on $\eps$ as parameter, tends to
the delta distribution weakly as $\eps\to0$.
Essentially, these paths are obtained by applying the scaling operator
$S_\eps:\vphi\mapsto\frac{1}{\eps^s}\vphi(\frac{.}{\eps})$
(thereby introducing the parameter $\eps$) to so-called test objects.
Typically, a test object is
some fixed element $\vphi\in\cd(\R^s)$ satisfying $\int\vphi=1$
or a suitable bounded family $\phi(\eps,x)\in\cd(\Om)$,
parametrized by $\eps\in I$, $x\in\Om$, where again
$\int\phi(\eps,x)(\xi)\,d\xi\equiv1$.
Roughly speaking, $R$ is defined to be negligible
if the values $R$ attains
on those ``scaled test objects''
tend to zero faster than any positive power of
$\eps$, while it is called moderate if these values are bounded
by some fixed (negative) power of $\eps$.
In both cases, convergence in each derivative, uniformly on
compact subsets of $\Om$, is required.
We will refer to those defining
procedures as {\bf testing for negligibility} resp.\ {\bf moderateness}
(see also section \rf{deftest}).

{\bf ad (T1), (T3):}
$\cn(\Om)$ has to be large enough to contain all $\si(f)-\io(f)$
($f\in\cc^\infty(\Om)$) (this renders $\io\left|_{\cc^\infty(\Om)}\right.$ 
an algebra homomorphism by
passing to a quotient by $\cn(\Om)$), however small enough to intersect
$\cd'(\Om)$ just in $\{0\}$ (this guarantees $\cd'(\Om)$ to be contained
injectively in the quotient by $\cn(\Om)$).
$\E_M(\Om)$, on the other hand,
clearly has to be large enough to contain $\cc^\infty(\Om)$
and $\cd'(\Om)$
(via $\si$ resp.\ $\io$), yet small enough such that $\cn(\Om)$ is an
ideal in it: This will allow us to form the quotient
$\E_M(\Om)\big/\cn(\Om)$.

{\bf ad (D5):}
$\mu^*$ as defined above extends $\mu^*:f\mapsto f\circ\mu$
where the latter is viewed as the action induced by $\mu$
on the smooth distribution $f\in\cc^\infty(\Om)$.
Hence we regard distributions (and, in the sequel, non-linear generalized functions) as
generalizations of {\it functions}
on the respective open set, acting as functionals on (smooth,
compactly supported) densities. This is in agreement
with, for example, \cite{Hoe}, however has to be distinguished
clearly from constructing distributions as distributional
{\it densities}, acting on (smooth, compactly supported) functions,
as it is done, e.g., in \cite{D3}.

{\bf ad (T7), (T8):}
Because of the forms of {\bf (D3)} and {\bf (D4)} as tests to be performed
on the elements $R$ of $\ce(\Om)$, with the appropriate type
of (scaled) test objects being inserted,
{\bf (T7)} as well as
{\bf (T8)} follow immediately from {\bf (T6)},
taking into account {\bf (D5)}.

{\bf ad (D6):}
By this definition, $\cg(\Om)$ is a differential algebra containing $\cd'(\Om)$
via $\io$ followed by the canonical quotient map
({\bf (T1)}--{\bf (T5)}); by abuse of notation, we will denote this
embedding also by $\io$.
Each diffeomorphism $\mu:\ti\Om\to\Om$ induces a map
$\hat\mu:\cg(\Om)\to\cg(\ti\Om)$ extending the usual action of $\mu$ on
distributions such that composition and identities are preserved by
$\mu\mapsto\hat\mu$ ({\bf (T7)}, {\bf (T8)}).

Without the requirement of diffeomorphism invariance (as, for example, in
\cite{c1}), the smoothness property of $R$ mentioned above only needs
to refer to the variable $x$ in the pair $(\vphi,x)$, thus
involving only classical calculus. However, 
as mentioned already in the introduction, to obtain a diffeomorphism
invariant algebra we also have to consider smoothness with respect to
the test function $\vphi$. Therefore, in the following section, we are going
to outline the elements of calculus on locally convex spaces
which are required for the subsequent constructions. The path we will pursue
in this respect is different from the approaches taken
so far and, in our view, has some decisive advantages over these.
\section{Calculus}\lb{calc}
In the first versions of Colombeau algebras (on $\R^n$ or open subsets
thereof), the main ingredient was the algebra of smooth functions 
$\vphi\mapsto R(\vphi)$ on the
((LF)-)space $\cd$ of test functions
(see \cite{c1}). Thus, from the very beginning, there
had to be a theory of differentiation on (certain non-Banach) locally
convex spaces at the basis of the construction of these algebras.

Colombeau's approach in \cite{c1} employs the notion of 
{\em Silva-differentiability} (\cite{Yama}, \cite{c0})
where a map $f:E\supseteq U\to F$ 
from an open subset $U$ of a locally convex space $E$ into another locally
convex space $F$ is called Silva-differentiable in $x\in U$ if there exists a
bounded linear map (called $f'(x)$) $E\to F$ such that 
the restriction of the corresponding remainder function to sufficiently
small homothetic images of bounded subsets may be viewed as a map between
suitable normed spaces and satisfies a condition thereon which is
completely analogous to the classical remainder condition for
Fr\'echet-differentiable maps.
  
In later versions, Colombeau managed to circumvent this necessity
by introducing an additional variable $x\in\R^n$ into $R$  which could carry 
the burden of
smoothness: For the construction of the algebra $\cg(\Om)$ of \cite{c2}
he now used functions $R(\vphi,x)$ which, for each fixed $\vphi$ from
(a certain affine subspace of) $\cd$, are smooth in $x$ (in the usual
elementary sense---hence the title of \cite{c2}) whereas the dependence
on $\vphi$ is completely arbitrary; $\vphi$ just plays the r\^ole of
a parameter in this setting. Apart from simplifying the general setup 
of the theory the introduction of 
$x$ as a separate variable was also crucial for solving differential equations
in ${\mathcal G}(\Omega)$.

However, when Colombeau and Meril
in \cite{CM} began to develop a diffeomorphism
invariant version of the algebra $\cg(\Om)$ of \cite{c2}, they had to
reintroduce the smooth dependence of $R$ on $\vphi$: Under the action
of a diffeomorphism $\mu$, $\vphi$ changes to some $\ti\vphi_x$ {\it
depending on $x$}. For the smoothness of the $\mu$-transform
of $R$
(which, according to {\bf (D5)}, is of the form
$(\hat\mu R)(\vphi,x)=R(\bar\mu(\vphi,x))=R(\ti\vphi_x,\mu x)$)
with respect to $x$, obviously the smooth dependence of $R$
also on its first argument
$\vphi$ is needed (\cite{CM}, p.~263). Concerning calculus on
locally convex spaces, the authors---as
the first of them did already in \cite{c1}---refer to \cite{c0}. 
Omitting any details in this respect, they rather
invite the reader to admit the respective smoothness properties (p.~263).

Jel\'\i nek in \cite{JEL} includes a section on calculus (items 9--16):
In addition to \cite{c0}, he quotes \cite{Yama} as reference
for some results needed. The relevant statements are formulated in
terms of higher Fr\'echet differentials.

Contrary to the above, we prefer to base our presentation on the notion
of smoothness  as it is outlined in \cite{KM}. This approach seems
to us to have a number of striking advantages: 
On the one hand, the basic definition is very simple
and easy to work with, a smooth map between locally convex spaces $E,F$ being
one that takes smooth curves $\R\to E$ to smooth curves $\R\to F$ (by
composition); the notion of a smooth curve into a locally convex space
obviously is without problems. We will denote by $\cc ^{\infty}(E,F)$
the space of smooth maps between $E$ and $F$.
For $\cc^\infty(E,\C)$, we will simply write $\cc^\infty(E)$.
On the other hand, all the 
basic theorems of differential calculus can be reconstructed in this setting
(see, e.g., the version of the mean value theorem given in \ref{mvth})
and more than that (see, e.g., the exponential laws stated in \ref{exp} and 
\ref{linexp} below and the differentiable uniform boundedness 
principle \ref{ubp}). As smooth curves are continuous, the above definition
of smoothness carries over to open subsets of locally convex spaces.
We will make use of this in the sequel and want to note that in any of the
theorems of this section, we may replace the respective locally convex
(domain)
spaces by open subsets thereof whenever their linear structure is not
needed.

This notion of smoothness is a weaker one than Silva-differentiability 
but turns out to be equivalent for a huge class of spaces, e.g. those which 
are complete and Montel so that the two notions coincide in particular on the
regular\footnote{A strict inductive limit
$\lim\limits_{\longrightarrow} E_\al$
is called regular if each bounded subset is contained in some
$E_\al$. Note that every strict inductive limit of an increasing
sequence $E_n$ is regular, as is $\cd(M)$ for any paracompact (not
necessarily separable) smooth manifold $M$.}
strict inductive limit $\cd(\Om)$ of Fr\'echet spaces and each closed 
subspace thereof.

The seeming drawback of this (and any other reasonable such as Colombeau's
above-mentioned) theory of differentiation is the fact that smooth
maps (resp.\ their differentials) need no longer be continuous. The
fundamental r\^ole played by
continuity in the classical context
is taken over by the notion of boundedness: Indeed, the
difference quotients of smooth curves converge in a stronger sense than
the topological one, so that continuity is not a necessary property for
a map to be smooth. In order to be able to test smoothness by composition
with suitable families of linear functionals (see, e.g., \ref{ubp}) 
one needs, in addition, a completeness property which is weaker
than completeness of the locally convex topology. Separated
bornological locally convex spaces which have this property
are called {\em convenient
spaces} and are in some sense the most general class of linear spaces in
which one can perform differentiation and integration. As for each locally
convex space there exists a finer bornological locally convex topology with
the same bornology, i.e., the same system of bounded sets, bornologicity
of the topology is not essential. 
It will be enough for our purpose to confine ourselves to the particular case of
complete locally convex spaces. 

In the sequel, we will endow the space $\cc ^{\infty}(\R,F)$ of smooth curves
into the locally convex space $F$ with the locally convex topology of uniform
convergence on compact intervals in each derivative separately. More
generally, we may consider on the space $\cc ^{\infty}(E,F)$ the initial
locally convex topology induced by the pullbacks along smooth curves 
$\R\to E$. It can be shown that the bounded sets associated with this topology
are the same as the ones associated with the topology of uniform convergence
on compact subsets in each differential (as defined for such maps in 
\ref{chain}) separately. Moreover, as mentioned in the introduction of this
paper, for complete $F$, the latter is again complete; see section
14 of the second part of this series (\cite{ft}) for details.

Testing of smoothness is particularly simple
in the case of a linear map: A linear map is smooth if and only if it is
bounded. $L(E,F)$ will stand for the space of bounded (smooth) linear maps
between $E,F$.

\bt\lb{smoothonD}
A map $f$ from $\cd(\Om)$ into a locally convex space $E$ is
smooth if and only if for each $K\subset\subset\Om$,
the restriction of $f$ to $\cd_K(\Om)$ is smooth.
\et
\pr For the non-trivial part of the proof,
consider a smooth curve $c:\R\to\cd(\Om)$. Its restriction to any
bounded interval $J$ has a relatively compact, hence bounded image.
Therefore, $c$ maps $J$ into some $\cd_K(\Om)$ and the same holds
for each derivative of $c$ since $\cd _K(\Om)$ is a closed subspace of
$\cd$. By assumption, $f\circ c$ is smooth on $J$. Since smoothness is a 
local property, we are done.
\ep
                                                       
The obvious generalization of the
preceding theorem is true for any strict inductive
limit of a sequence of Fr\'echet spaces.
Its trivial part has an important consequence:
$\cd_K(\Om)$ being a Fr\'echet space, the restriction to $\cd_K(\Om)$
of any smooth map $f$ from $\cd(\Om)$ to any metrizable locally convex
space $E$ is continuous:
Both on $\cd_K(\Om)$ and $E$ the so-called
$c^\infty$-topology (see \cite{KM}) coincides with
the metric topology (\cite{KM}, 4.11.(1)); moreover, smooth maps are
continuous with respect to the $c^\infty$-topology
(\cite{KM}, p.~8).
 
One of the particular features of the Fr\"olicher-Kriegl-theory which 
considerably simplify its application is the {\em exponential law} 
(cf. Theorem 3.12 and Corollary 3.13 in \cite{KM}):

\bt \label{exp}
Let $E,F,G$ be locally convex spaces. Then the two spaces 
$\cc ^{\infty}(E\times F,G)$ and $\cc ^{\infty}(E,\cc ^{\infty}(F,G))$
are isomorphic algebraically and bornologically, i.e., they have the same
bounded sets.
\et

Replacing $\cc ^{\infty}$ by $L$ in \ref{exp} yields the exponential law
for linear smooth maps. By iteration one obtains 
(see Proposition 5.2 in \cite{FK}):

\bt \label{linexp}
Let $n,k\in\N$ and $E_i,F$ ($i=1,\dots ,n+k$) locally convex
spaces. Then there is a bornological isomorphism
\[
L(E_1,\dots ,E_{n+k};F)\cong 
L(E_1,\dots ,E_n; L(E_{n+1},\dots ,E_{n+k};F)).
\]                                          
\et

For later use, we present the analoga of items 10--16 in \cite{JEL}
in the setting of \cite{KM}:                                                     

\bt (Theorem 3.18 and Corollary 5.11 in \cite{KM}) \label{chain}
Let $E,F$ be locally convex spaces. Then the differentiation operator
$\mathrm{d}:\cc ^{\infty}(E,F)\to \cc ^{\infty}(E, L(E,F))$ given by
\[
\mathrm{d}f(x)v:=\lim_{t\to 0}
\frac{f(x+tv)-f(x)}{t}
\]
exists and is linear
and bounded (smooth). Hence, for $n\in\N$ one can form the iterated
differentiation operator
\[
\mathrm{d} ^n:\cc ^{\infty}(E,F)\to 
\cc ^{\infty}(E,L(E,\dots ,L(E;F)\dots ))\cong 
\cc ^{\infty}(E,L(E,\dots ,E;F))
\]
which is smooth and linear and has values in
$\cc ^{\infty}(E,L_{sym}(E,\dots ,E;F))$, where
$L_{sym}(E,\dots ,E;F)$ stands for the space of smooth $n$-linear symmetric
maps between $E\times\dots\times E$ and $F$. Also, the chain rule holds:
\[
\mathrm{d}(f\circ g)(x)v=\mathrm{d}f(g(x))\mathrm{d}g(x)v.
\]
\et

It is shown in \cite{KM}, 1.4, that, given a curve which is smooth from 
(an open neighborhood of) $\R\supseteq [a,b]$ to $E$, the difference 
quotient $\frac{c(b)-c(a)}{b-a}$ is an element of 
$\overline{\mathrm{conv}}\{c'(t):t\in [a,b]\}$, where 
$\overline{\mathrm{conv}}$ denotes the closed convex hull. By virtue of the
chain rule given in \ref{chain}, this is equivalent to

\bp  \label{mvth} (Mean Value Theorem)
Let $f:E\supseteq U\to F$ be smooth, where $U$ is an open neighborhood 
of a segment $[x,x+v]\subseteq E$. Then
\[
f(x+v)-f(x)\in\overline{\mathrm{conv}}\{\mathrm{d}f(x+tv)(v):t\in [0,1]\}.
\] 
\et

As a consequence of \ref{exp}, for each smooth map
$f\in \cc ^{\infty}(F,G)$, the maps 
$f _*:\cc ^{\infty}(E,F)\to\cc ^{\infty}(E,G)$ and
$f ^*:\cc  ^{\infty}(G,E)\to\cc ^{\infty}(F,E)$ are smooth. 
In particular, for a smooth map $f\in\cc ^{\infty}(E\times F,G)$
we may define smooth linear ``operators of partial differentials'' 
$\mathrm{d}_1,\mathrm{d}_2$ as 
\[
\mathrm{d}_1:=(\iota _E ^*)_*\circ \mathrm{d}:
\cc ^{\infty}(E\times F,G)\to\cc ^{\infty}(E\times F, L(E,G))
\] and
\[
\mathrm{d}_2:=(\iota _F ^*)_*\circ \mathrm{d}:
\cc ^{\infty}(E\times F,G)\to\cc ^{\infty}(E\times F,L(F,G)),
\]
where $\iota _E,\iota _F$ denote
the natural embeddings of $E$ resp.\ $F$ into $E\times F$. Obviously,
we have 
\[
\mathrm{d}_1f(x)(v)=\mathrm{d}f(x)(\io _E(v))=
\lim_{t\to 0}\frac{f(x+t\io_E(v))-f(x)}{t},
\]
which yields an alternative definition of $\mathrm{d}_1$, which makes sense
also for maps $f:E\times F\to G$ which are not a priori known to be 
smooth on $E\times F$.

\bp \label{part}
A map on $E\times F$ is smooth if and only if both partial 
differentials $\mathrm{d}_1,\mathrm{d}_2$ exist and are smooth as maps
on $E\times F$. In this case the differential $\mathrm{d}$ equals the
sum $(\pro_1^*)_*\circ\rmd_1+(\pro_2^*)_*\circ\rmd_2$
of the partial
differentials; the iterated mixed second derivatives coincide via the
isomorphism $L(E,L(F,G))\cong L(F,L(E,G))$ which is a consequence of 
\ref{exp}.
\et

\pr
Necessity follows by what has been remarked above together with the symmetry
of iterated derivatives stated in \ref{chain}. For sufficiency, consider the map
$\tilde{\mathrm{d}}f\in \cc  ^{\infty}(E\times F, L(E\times F,G))$ defined by
$\tilde{\mathrm{d}}f(x)(v_1,v_2):= \mathrm{d}_1 f(x)(v_1)+\mathrm{d}_2 f(x)(v_2)$. Then
obviously for fixed $x$ the map $(t,v)\mapsto\tilde{\mathrm{d}}f(x+tv)(v)$ is smooth
from $[0,1]\times E\times F\to G$ and hence can be viewed as an element of 
$\cc ^{\infty}([0,1],\cc ^{\infty}(E\times F,G))$. By \cite{KM}, 2.7, a smooth curve
is Riemann integrable, the Riemann integral leads again into 
$\cc ^{\infty}(E\times F,G)$ and commutes with the application of smooth 
linear maps. It follows that the map 
\[
v\mapsto f(x) + \int _0 ^1\tilde{\mathrm{d}}f(x+tv)(v)dt
\]  
is smooth on $E\times F$ and it suffices to verify that the expression on the right
hand side equals $f(x+v)$ in order to obtain smoothness of $f$ on $E\times F$. For this,
note that for each fixed segment $[x,x+v]$, we can recover the claimed identity from the
finite dimensional one by 
composing the restriction of $f$ to the segment with bounded linear 
functionals.                      
\ep

The {\em differentiable uniform boundedness principle} 
(see 4.4.7 in \cite{FK})
constitutes an extremely useful tool for testing smoothness of linear maps
into spaces of smooth functions:

\bt \label{ubp} Let $E,F,G$ be locally convex spaces,
$E,G$ complete.
A linear map $E\to \cc  ^{\infty}(F,G)$ is smooth if and and only if its
composition with the evaluation $\mathrm{ev}_x$ for each $x\in F$ is
smooth.
\et

If we endow the space
$\cc^\infty(X,\R)$ (in the present paper, $X$ will be one of the spaces
$\cd(\Om),\cd(\Om)\times\Om,\ca_0(\Om)\times \Om$ or 
$\ca_0(\R^s)\times\Om$) with the topology of uniform convergence on
compact subsets in each derivative, i.e., in each iterated differential
separately, then by considering the corresponding seminorms one sees that
taking the differential constitutes a continuous linear operation.
To be precise, the space $\ca_0(\R ^s)$ is not a linear space itself but
the affine image of the closed linear subspace
$E:=\ca_{00}(\R^s)\subseteq \cd(\R^s)$ and may be identified 
with the latter. A map on $\ca _0(\R ^s)$ is then said to be smooth if it is
the pullback of a smooth map on $E$ under the affine isomorphism. We say
that the smooth structure on $\ca _0(\R ^s)$ is induced by its isomorphism
with $E$.
This is a simple example of the notion of a {\em smooth space} as 
introduced in \cite{FK}. Locally convex spaces may be viewed as smooth
spaces with a compatible linear structure.
\bp\lb{smoothmaps}        
The following maps (to be defined in
section \rf{basics}) are smooth: The linear maps
$S_\eps:\cd (\R ^s)\to\cd (\R ^s)$, 
$T_x:\cd (\R ^s)\to\cd (\R ^s)$ and
$(\vphi,x)\mapsto\bar\mu^X(\vphi,x)$,
$R\mapsto\hat\mu^X R$
$(X\in\{C,J\})$, as well as the non-linear maps
$S$, $T$, $x\mapsto T_x$, $x\mapsto T_x\vphi$, $(\vphi,x)\mapsto T_x\vphi$.
\et

\pr
Smoothness of $S_{\eps},T_x,\bar\mu^X,\hat\mu^X$
follows by our remarks preceding \ref{smoothonD}
and following \ref{mvth}, 
respectively, as each of these maps is essentially a pullback of a smooth
map by definition. As the map $S:(\eps,\vphi)\mapsto S_{\eps}\vphi)$ is
linear in $\vphi$, it follows by the exponential law \ref{exp} and the 
uniform boundedness principle \ref{ubp} that $S$ is smooth iff it is
separately smooth, i.e., if and only if the maps $S_{\eps}$ 
and $(\eps\mapsto S_{\eps}\vphi)$ are smooth. While smoothness of the 
former is already established, the latter is a curve which is obviously
smooth off $0$ and we are done. In a similar fashion, we obtain smoothness
of $T$ and all the maps associated with it. 
\ep 

\section{C- and J-formalism}\lb{basics}
Colombeau in \cite{c2} and in \cite{CM} (together with Meril) on the one
hand and Jel\'\i nek in \cite{JEL}
on the other hand used different, yet equivalent formalisms
to describe their respective constructions of Colombeau algebras: For
embedding the space $\cd'(\R^s)$ of distributions on $\R^s$ into
the space $\ce_M(\R^s)$ of representatives of generalized functions,
they chose different (linear injective) maps which we denote by $\io^C$
(\cite{c2}, \cite{CM}) and $\io^J$ (\cite{JEL}, compare also~\cite{c1}),
respectively. On a distribution given by a smooth function $f$
on $\R^s$, $\io^C$ and $\io^J$ are defined by
\be\lb{iocf}
(\io^Cf)(\vphi,x):=\int f(y)\vphi(y-x)\,dy
\end{equation}
resp.
\be\lb{iojf}
(\io^Jf)(\vphi,x):=\int f(y)\vphi(y)\,dy.
\end{equation}
Here, $\vphi$ denotes a test function from the
subspace $\ca_0(\R^s)$ of $\cd(\R^s)$
while $x\in\R^s$. There are good reasons for either of these choices
of the embedding; we are going to discuss their respective merits below. 
In this section we show that both formalisms are actually equivalent and
establish a translation formalism allowing to change from one setting to the 
other at any stage of the presentation.

\bd 
For $\eps\in I$ and $x\in\R^s$ define the following operators:
\bea T_x:\cd(\R^s)\ni\vphi&\mapsto& T_x\vphi:=\vphi(.-x)\in\cd(\R^s)\\
     S_\eps:\cd(\R^s)\ni\vphi&\mapsto& S_\eps\vphi:=
         \frac{1}{\eps^s}\vphi\left(\frac{.}{\eps}\right)\in\cd(\R^s)\\
    S:(0,\infty)\times\cd(\R^s)\ni(\eps,\vphi)&\mapsto&S_\eps\vphi
    \in\cd(\R^s)\\
    T:\cd(\R^s)\times\R^s\ni(\vphi,x)
         &\mapsto&T(\vphi,x):=(T_x\vphi,x)
         \in\cd(\R^s)\times\R^s\\
     S^{(\eps)}:\cd(\R^s)\times\R^s\ni(\vphi,x)
         &\mapsto& S^{(\eps)}(\vphi,x):=(S_\eps\vphi,x)
         \in\cd(\R^s)\times\R^s.
\eea
\et

$T_x$ and $S_\eps$ are linear.
All the operators introduced in the preceding definition are
one-one and onto; moreover, they are continuous
and smooth with respect to the natural topologies (see 
section \rf{calc}).

In a next step, we take (\rf{iocf}) and (\rf{iojf}) as a starting point
for the determination of suitable domains for representatives
of generalized functions on an open subset $\Om$ of $\R^s$:
Assuming $x\in\Om$ in (\rf{iocf}) and (\rf{iojf}), it is immediate that
in (\rf{iojf}) $\vphi$ has to have its support in $\Om$, whereas for
(\rf{iocf}) to be well-defined for any smooth function $f$ on $\Om$,
the support of $\vphi$ must be contained in $\Om-x$. This
motivates the introduction of the following sets:

\bd
Let $\eps\in I$.
\beas
U(\Om)&:=&T^{-1}(\ca_0(\Om)\times\Om)\hphantom{(S^{(\eps)})}=
          \{(\vphi,x)\in\ca_0(\R^s)\times\Om\mid
          \spp\vphi\subseteq\Om-x\}\\
U_\eps(\Om)&:=&(S^{(\eps)})^{-1}(U(\Om))=\\
           &\hphantom{:}=&(TS^{(\eps)})^{-1}(\ca_0(\Om)\times\Om)
           =\{(\vphi,x)\in\ca_0(\R^s)\times\Om\mid
          \spp\vphi\subseteq\eps^{-1}(\Om-x)\}
\eeas
\et
The notation $U(\Om)$ is due to Colombeau (\cite{c2}, 1.2.1).
By definition, the maps $T:U(\Om)\to\ca_0(\Om)\times\Om$ and
$S^{(\eps)}:U_\eps(\Om)\to U(\Om)$ are algebraic isomorphisms
in the sense that they are bijective and linear in the first
argument. The question of topology, however, is somewhat subtle:
Let $\tau_\Om$ denote the product of the (LF)-topology of
$\cd(\Om)$ and the Euclidean topology on $\Om$; abbreviate
$\tau_{\R^s}$ as $\tau_0$. Then on $\ca_0(\Om)\times\Om$,
the topology $\tau_\Om$ without doubt is the appropriate one to consider,
rather than (the restriction of) $\tau_0$. For $U(\Om)$, on the
other hand, the topology $\tau_1$ induced by $\tau_0$ and
the topology $\tau_2:=T^{-1}\tau_\Om$ both seem to be natural
choices. (Note that $\tau_1$ can be obtained equally
as $T^{-1}\tau_0$, due to $T$ being a homeomorphism with respect to $\tau_0$.)
As the following example (which can easily be generalized to arbitrary
non-trivial open subsets of $\R^s$) shows,
$\tau_1$ is strictly coarser than $\tau_2$ in general.

\bex
Let $\Om:=\{x\in\R\mid x>-1\}$. Choose $\vphi\in\cd(\Om)$ with 
$\spp\vphi=[0,1]$ and $\int\vphi=0$. Pick any
$\rho\in\ca_0(\Om)$ such that $\spp\rho\subseteq[1,2]$.
Letting $\psi_n:=\rho+\frac{1}{n}\vphi(.+\frac{n-1}{n})\in\ca_0(\Om)$,
it is easy to check that $T^{-1}(\psi_n,0)=(\psi_n,0)\in U(\Om)$ 
tends to $(\rho,0)\in U(\Om)$ with respect to $\tau_1$, yet
is not convergent (in fact, not
even bounded) with respect to $\tau_2$.
\et

The situation is similar in the case of 
$U_\eps(\Om)$: Apart from the topology $\tau_{1,\eps}$
induced by the topology $\tau_0$
of $\ca_0(\R^s)\times\R^s$ via inclusion,
the natural topology $\tau_\Om$ of $\ca_0(\Om)\times\Om$
via $TS^{(\eps)}$ induces a topology $\tau_{2,\eps}$ 
which, in general, is strictly finer than $\tau_{1,\eps}$.

Now, in order to have the respective formalisms of Colombeau
and Jel\'\i nek equivalent,
we want $T:U(\Om)\to\ca_0(\Om)\times\Om$ and
$S^{(\eps)}:U_\eps(\Om)\to U(\Om)$ to be also topological
isomorphisms (hence diffeomorphisms). This amounts to endowing
$U(\Om)$ and $U_\eps(\Om)$ with the topologies $\tau_2$
resp.\ $\tau_{2,\eps}$ induced via $T$ resp.\ $TS^{(\eps)}$.
Thus we adopt the following convention:

{\it Whenever questions of topology (in particular,
boundedness) or smoothness on $U(\Om)$ or $U_\eps(\Om)$
are discussed, we consider their topologies to be $\tau_2$
resp.\ $\tau_{2,\eps}$, i.e., those induced
by the natural topology of $\ca_0(\Om)\times\Om$
via $T$ resp.\ $TS^{(\eps)}$.} 

To phrase it differently, $U(\Om)$ can be viewed as
(infinite-dimensional) smooth ma\-ni\-fold, modelled over
$\ca_0(\Om)\times\Om$, having an atlas consisting of a single chart
$T$. A similar statement is valid for $U_\eps(\Om)$ and
$TS^{(\eps)}$.
The importance as well as the subtlety of distinguishing
between $\tau_1$ and $\tau_2$ are highlighted
in example \ref{hilfe} below.

We are now able to introduce the basic spaces of smooth functions on which the
construction of diffeomorphism invariant Colombeau algebras is built.

\bd\lb{ceom}
\bea \ce^J(\Om)&:=&\cc^\infty(\ca_0(\Om)\times\Om)\\
     \ce^C(\Om)&:=&\cc^\infty(U(\Om))
\eea
\et

By our above choice of topologies,
$T^*$ indeed maps $\ce^J(\Om)$ bijectively onto $\ce^C(\Om)$.
The next definition shows how the space of distributions on $\Om$
is to be embedded into $\ce^J(\Om)$ resp.\ $\ce^C(\Om)$.

\bd\lb{io}
For $u\in\cd'(\Om)$, define
$$\ba{ll}
\io^J\kern+1pt:\cd'(\Om)\to\ce^J(\Om)\qquad\qquad&
             (\io^J u)(\vphi,x):=\lgl u,\vphi\rgl\\
\io^C:\cd'(\Om)\to\ce^C(\Om)\qquad\qquad&
             (\io^C u)(\vphi,x):=\lgl u,\vphi(.-x)\rgl
\ea$$
\et
By definition, $\io^C=T^*\circ\io^J$.

It remains to introduce the respective
extensions of partial differentiation 
from $\cd'(\Om)$ to $\ce^C(\Om)$ resp.\ $\ce^J(\Om)$
and the respective actions of a diffeomorphism.

\bd\lb{deriv}
For $i=1,\dots,s$, define
$$\ba{ll}
  D_i^C         :\ce^C(\Om)\to\ce^C(\Om)\qquad\qquad&
         D_i^C:=\pa_i,\\
  D_i^J\kern+1pt:\ce^J\kern+1pt(\Om)\to\ce^J\kern+1pt(\Om)
                                           \qquad\qquad&
  D_i^J\kern+1pt:=(T^*)^{-1}\circ\pa_i\circ T^*,
\ea$$
i.e., for $R\in\ce^J\kern+1pt(\Om),\ (\vphi,x)\in
                                 \ca_0(\Om)\times\Om$
we set
$$(D_i^JR)(\vphi,x):=-((\rmd_1 R)(\vphi,x))(\pa_i\vphi)
                   +(\pa_i R)(\vphi,x).$$
\et

Of course we have to demonstrate that for given $R\in\ce^C(\Om)$
and $(\vphi,x)\in U(\Om)$, $(D_i^C R)(\vphi,x)$
in fact exists and that $(\vphi,x)\mapsto (D_i^C R)(\vphi,x)$
is smooth on $U(\Om)$ with respect to $\tau_2$ (and similar
for $R\in\ce^J(\Om)$ and $D_i^J$).
This being a non-trivial task---in particular for
the case of the innocent-looking map $D_i^C=\pa_i$
[{\it sic}!]---requiring some technical prerequisites, we
have to defer it to the following section.

Commutativity of the following diagram is immediate:

\[
\begin{CD}
\D'(\Om)        @>\pa_i>>       \D'(\Om)        \\
                @VV\iota^CV     @VV\iota^CV                     \\
\E^C(\Om)       @>D_i^C>>       \E^C(\Om)       \\
                @AAT^*A         @AAT^*A                         \\
\E^J(\Om)       @>D_i^J>>       \E^J(\Om)       \\
\end{CD}
\]

\bd\lb{barmu}\lb{barmueps}
Let $\mu:\ti\Om\to\Om$ be a diffeomorphism. Define
$$\ba{rccc}
  \bar\mu^J\kern+1pt:&\ca_0(\ti\Om)\times\ti\Om&\to&
                   \ca_0(\Om)\times\Om\\
  \bar\mu^C:&         U(\ti\Om)&\to& U(\Om)\\
  \bar\mu_\eps\,:&U_\eps(\ti\Om)&\to& U_\eps(\Om)
\ea$$
by
\beas
\bar\mu^J\kern+1pt(\ti\vphi,\ti x)
     &:=&\left((\ti\vphi\circ\mu^{-1})\cdot
                                   |\det D\mu^{-1}|\,,\,\mu\ti
                                   x\right),\\
\bar\mu^C(\ti\vphi,\ti x)&:=&
    \left(T^{-1}\circ\bar\mu^J\circ T\right)(\ti\vphi,\ti x)\\
    &\,=&\left(\ti\vphi(\mu^{-1}(.+\mu\ti x)-\ti x)\cdot
        |\det D\mu^{-1}(.+\mu\ti x)|\,,\,\mu\ti x\right).\\
    \bar\mu_\eps(\ti\vphi,\ti x)&:=&
    \left((S^{(\eps)})^{-1}\circ T^{-1}\circ\bar\mu^J\circ T\circ
    S^{(\eps)}\right)(\ti\vphi,\ti x)\\
    &\,=&\left(\ti\vphi
    \left(\frac{\mu^{-1}(\eps.+\mu\ti x)-\ti x}{\eps}\right)\cdot
        |\det D\mu^{-1}(\eps.+\mu\ti x)|\,,\,\mu\ti x\right).
\eeas
\et

\bd\lb{hatmu}
Let $\mu:\ti\Om\to\Om$ be a diffeomorphism and $\eps\in I$. Define
$$\ba{rccc}
  \hat\mu^J\kern+1pt:&\ce^J(\Om)&\to&\ce^J(\ti\Om)\\
  \hat\mu^C:         &\ce^C(\Om)&\to&\ce^C(\ti\Om)\\
  \hat\mu_\eps\kern+2pt:      &\cc^\infty(U_\eps(\Om))&\to&
                               \cc^\infty(U_\eps(\ti\Om))
\ea$$
by $\hat\mu^J:=(\bar\mu^J)^*$, $\hat\mu^C:=(\bar\mu^C)^*$,
$\hat\mu_\eps:=(\bar\mu_\eps)^*$, i.e.,
\beas
(\hat\mu^J\kern+1pt R)(\ti\vphi,\ti x)&:=&R(\bar\mu^J
                \kern+1pt(\ti\vphi,\ti x))
    \qquad(R\in\ce^J\kern+1pt(\Om),\hphantom{mma}(\ti\vphi,\ti x)\in
                                 \ca_0(\ti\Om)\times\ti\Om),\\
(\hat\mu^C          R)(\ti\vphi,\ti x)&:=&R(\bar\mu^C(\ti\vphi,\ti x))
    \qquad(R\in\ce^C(\Om),\hphantom{mma}(\ti\vphi,\ti x)\in U(\ti\Om)),\\
(\hat\mu_\eps\kern+2pt R)(\ti\vphi,\ti x)&:=&R(\bar\mu_\eps
                \kern+2pt(\ti\vphi,\ti x))
    \qquad(R\in\cc^\infty(U_\eps(\Om)),\ (\ti\vphi,\ti x)\in U_\eps(\ti\Om)).
\eeas
\et

For   $X\in\{C,J\}$ we obtain

\[
\begin{CD}
\D'(\Om)        @>\mu^*>>       \D'(\tilde\Om)  \\
@VV\iota^XV     @VV\iota^XV                     \\
\E^X(\Om)       @>\hat\mu^X>>   \E^X(\tilde\Om) \\
\end{CD}
\]

\vskip2mm

\vskip2mm
where for $u\in\cd'(\Om)$, $\mu^*u$ is defined by
$\lgl\mu^*u,\vphi\rgl:=
                  \lgl u,(\vphi\circ\mu^{-1})\cdot|\det D\mu^{-1}|\rgl$
($\vphi\in\cd(\Om)$)
which extends $f\mapsto\mu^*f=f\circ\mu$ for $f\in\cc^\infty(\Om)$.

\begin{minipage}{.45\textwidth}
\[
\begin{CD}
U_\eps(\tilde\Om)@>\bar\mu_\eps>> U_\eps(\Om)
 \\ 
@VVS^{(\eps)}V
        @VVS^{(\eps)}V
\\
U(\tilde\Om)    @>\bar\mu^C>>   U(\Om)          
\\
        @VV T V                 @VV T V         
        \\
{\cal A}_0(\tilde\Om)\times\tilde\Om    
                @>\bar\mu^J>>   {\cal A}_0(\Om)\times\Om        
\end{CD}
\]
\end{minipage}
\hfill
\begin{minipage}{.45\textwidth}
\[
\begin{CD}
\CC^\infty(U_\eps(\tilde\Om))   @<\hat\mu_\eps<<        \CC^\infty(U_\eps(\Om)) \\
        @AA(S^{(\eps)})^*A      @AA(S^{(\eps)})^*A                              \\
\E^C(\tilde\Om)                 @<\hat\mu^C<<           \E^C(\Om)               \\
        @AA T^* A               @AA T^* A                                       \\
\E^J(\tilde\Om)                 @<\hat\mu^J<<           \E^J(\Om)               \\
\end{CD}
\]
\end{minipage}

\vskip5mm
Definitions \rf{barmu} and \rf{hatmu} reflect the fact
that in Definition {\bf (D5)} of section \rf{scheme1}, we
chose to regard distributions
(and, in the sequel, non-linear generalized functions) as
generalizations of {\it functions}, acting as functionals on test densities
(compare, e.g., \cite{Hoe}). This approach has to be distinguished
from constructing distributions as distributional
{\it densities}, acting on test functions (see, e.g., \cite{D3}).

In the following table, we compare the C-formalism and the
J-formalism regarding simplicity of the respective definitions
and, in the last item, the degree of familiarity.

\vskip3mm
\begin{center}
\begin{tabular}{|c|c|c|}
   \hline
      C-formalism&Feature&                                   J-formalism\\\hline\hline
        $-$   &domain of basic space $\ce(\Om)$&                $+$\\\hline
        $-$   &smoothness structure&                            $+$\\\hline
        $-$   &embedding of $\cd'(\Om)$&                        $+$\\\hline
        $+$   &formula for differentiation&                     $-$\\\hline
        $+$   &solving differential equations&                  $-$\\\hline
        $-$   &action induced by a diffeomorphism&              $+$\\\hline
        $+$   &formula for testing&                             $-$\\\hline
        $-$   &generalization to manifolds&                     $+$\\\hline
        $+$   &tradition&                                       $-$\\\hline
\end{tabular}
\end{center}

\vskip3mm
The distribution of the $+$'s and $-$'s in the table should be
rather obvious by inspecting the corresponding definitions.
Due to the absence of a linear structure on a general smooth manifold,
it is clear that the C-formalism does not lend itself to a
definition of non-linear generalized functions
on manifolds based only on intrinsic terms,
whereas the J-formalism in fact does
permit such a construction; see \cite{vi}.

We conclude this section by presenting an example 
that emphasizes the importance of
carefully distinguishing between the topologies $\tau_1$ and $\tau_2$
on $U(\Om)$.

\bex\lb{hilfe}
We will specify an open subset $\Om$ of $\R^2$, a line segment
of the form $\Phi(t):=
(\vphi,z)+t(\psi,v)$
($-1\le t\le1$) in $U(\Om)$
and a distribution $u$ on $\Om$ such that

$$\lim_{t\to0}\frac{(\io^Cu)(\Phi(t))-(\io^Cu)(\Phi(0))}{t}=\infty.$$

This seems
to suggest that in the point $(\vphi,z)$, the function
$\io^Cu$ on $U(\Om)$ (which ought to be smooth according to our
definitions) has no directional derivative with respect to $(\psi,v)$;
or, to phrase it differently, that the composition of the functions
$\io^Cu$ and $\Phi$ (both of which have the appearence of being smooth)
is not even differentiable in $(\vphi,z)$. We will leave the solution
to this puzzle for the end of the example. First we give the details of
the construction.

Let $\Om:=\{(x,y)\in\R^2\mid x>-y^2-1\}$, $z:=(0,0)$, $v:=(0,1)$.
Choose $\rho_1\in\cd(\R)$ such that $\spp\rho_1\subseteq[0,\frac{3}{2}]$
and $\rho_1(x)=\mbox{\rm exp}(-\frac{1}{x})$ on
$ I$. Let $c:=\int\rho_1$ and
choose $\rho_2\in\cd(\R)$ such that $\spp\rho_2\subseteq[\frac{3}{2},2]$
and $\int\rho_2=1$. Then $\rho:=\rho_1-c\rho_2$ has its support contained
in $[0,2]$, coincides with $\mbox{\rm exp}(-\frac{1}{x})$ on
$ I$ and satisfies
$\int\rho=0$. Pick $\om\in\cd(\R)$ with the properties $\spp\om=[-2,+2]$,
$0\le\om\le1$ and $\om\equiv1$ on $[-1,+1]$.
Now define $\psi\in\cd(\R^2)$ by $\psi(x,y):=\om(y)\cdot\rho(x+1-y^2).$

Finally, in order to obtain $\Phi$ as defined above, take any
$\vphi\in\ca_0(\R^2)$ whose support is located at the right hand side
of the line given by $x=6$. It is easy to check that 
$\Phi(t)=(\vphi+t\psi,z+tv)$ belongs to $U(\Om)$
for all $t\in[-1,+1]$.

Now there is still $u$ to be defined. To this end, let
$$f(x):=\left\{\ba{ll}
               \frac{1}{x^2}\mbox{\rm exp}
               (\frac{1}{x}+\frac{2}{x^2})\qquad&(0<x<1)\\
               0            &(x\ge1)\ea
        \right..$$

For $\si\in\cd(\Om)$ define the distribution $u\in\cd'(\Om)$ by
$$\lgl u,\si\rgl:=\int\limits_{-1}^0 f(x+1)\si(x,0)\,dx=
                  \int\limits_{0}^1 f(x)\si(x-1,0)\,dx.$$
For $0<|t|\le1$, it follows
\beas
\frac{1}{t}[(\io^Cu)(\vphi+t\psi,z+tv)-(\io^Cu)(\vphi,z)] =
  \frac{1}{t}t\lgl u,\psi(x,y-t)\rgl
  = \int\limits_0^1f(x)\rho(x-t^2)\,dx.
\eeas
We will show that for $0<|t|\le\frac{1}{\sqrt{2}}$,
the value of the last integral can be
estimated from below by
$\mbox{\rm exp}(\frac{1}{t^2}-1)-\mbox{\rm exp}(1)$, thus tending
to infinity for $t\to0$. Substituting $x=\frac{1}{u}$,
$t^2=\frac{1}{v}$, we obtain
\beas
\int\limits_0^1f(x)\rho(x-t^2)\,dx&=&
         \int\limits_{1}^{v}e^{u+2u^2}e^{-\frac{vu}{v-u}}\,du
     \ge \int\limits_{1}^{v-1}e^{2u^2}e^{-\frac{u^2}{v-u}}\,du
     \ge \int\limits_{1}^{v-1}e^{2u^2}e^{-u^2}\,du\\
    &\ge&\int\limits_{1}^{v-1}e^u\,du
     =   e^{v-1}-e
     =   e^{\frac{1}{t^2}-1}-e.
\eeas

The apparent inconsistencies mentioned at the beginning of the example
dissolve by taking into account that, in fact, both $\tau_1$ and $\tau_2$
are involved in the argument: The statement that $\Phi:[-1,+1]\to U(\Om)$
is smooth is true only if it refers to $\tau_1$ (the image of any
neighborhood of $0$ under $\Phi$ is even unbounded with respect to
$\tau_2$ since the supports of $T(\Phi(t))$ are not contained in
any compact subset of $\Om$ around $t=0$). The statement that
$\io^Cu$ is smooth is true only if $U(\Om$) is endowed with the
topology $\tau_2$ induced by the natural topology $\tau_\Om$ of
$\ca_0(\Om)\times\Om$ via $T$. $\tau_2$ being strictly finer than
$\tau_1$, we cannot infer the differentiability of $(\io^Cu)\circ\Phi$
from the actual smoothness properties of $\io^Cu$ resp.\ $\Phi$.
Another way of capturing the problem is by pointing out that
$(\psi,v)$ is not a member of the tangent space to $U(\Om)$
at $(\vphi,z)$ (in the sense of the following section)
since $\supp\psi$ is not contained in $\Om-z=\Om$.
\et
\section{Calculus on $U_\eps(\Om)$}\lb{calcuepsom}
The purpose of this section is to develop an appropriate framework for
defining and handling differentials of any order of a function
$f:U_\eps(\Om)\to\C$
which is smooth with respect to $\tau_{2,\eps}$ (by definition, $f$
is of the form $f_0\circ T\circ S^{(\eps)}$ where
$f_0\in\cc^\infty(\ca_0(\Om)\times \Om)$).
By choosing $\eps=1$, this includes the case of smooth functions on
$U(\Om)$, i.e., of elements of the basic space $\ce^C(\Om)$.
As a matter of fact, the author of \cite{JEL} has payed
only minor attention to these questions. However, it should be clear
even from a glimpse at sections
\rf{jelshort} and \rf{atoz}, in particular,
that a sound definition and a proper handling of the differentials
of $R_\eps:=R^J\circ T\circ S^{(\eps)}=R^C\circ S^{(\eps)}$
are crucial for the construction of a diffeomorphism
invariant Colombeau algebra.

To start with, we discuss
an important property of the sets $U_\eps(\Om)$ which
will be fundamental in the sequel at many places.
Loosely speaking, every subset of $\ca_0(\R^s)\times\Om$ which is
``not too large'' finally gets into $U_\eps(\Om)$ by scaling
and does not feel any difference between
$\tau_{1,\eps}$ and $\tau_{2,\eps}$.
To this end, we introduce the following notation:

\bd\lb{UK}
For every compact subset $K$ of $\Om$ define
\beas
     \ca_{0,K}(\Om)&:=&\{\vphi\in\ca_0(\Om)\mid\supp\vphi\subseteq
     K\},\\
     \ca_{00,K}(\Om)&:=&\{\vphi\in\ca_{00}(\Om)\mid\supp\vphi\subseteq
     K\},\\U_K(\Om)&:=&T^{-1}(\ca_{0,K}(\Om)\times\Om),\\
     U_{\eps,K}(\Om)&:=&(S^{(\eps)})^{-1}(U_K(\Om)).
\eeas
\et
By definition, we have
\beas
     U_K(\Om)&=&\{(\vphi,x)\in\ca_0(\R^s)\times\Om\mid
          \spp\vphi\subseteq K-x\},\\
     U_{\eps,K}(\Om)&=&(S^{(\eps)})^{-1}
                T^{-1}(\ca_{0,K}(\Om)\times\Om)\\
     &=&\{(\vphi,x)\in\ca_0(\R^s)\times\Om\mid
          \spp\vphi\subseteq\eps^{-1}(K-x)\}.
\eeas
Then it is immediate that
for $K\subset\subset\Om$, 
the topologies on $\ca_{0,K}(\Om)$ 
inherited from the natural topologies of $\ca_0(\R^s)$ and
$\ca_0(\Om)$, respectively, coincide. Consequently,
on $U_K(\Om)$ the topologies $\tau_1$ and $\tau_2$ are equal,
as are $\tau_{1,\eps}$ and $\tau_{2,\eps}$ on $U_{\eps,K}(\Om)$.

We now are in a position to
complement Definition \rf{deriv} by
establishing that the derivation operators $D_i^C$ and $D_i^J$
are in fact well-defined.
From 
the explicit formulas for $D_i^C$ resp.\
$D_i^J$ one is certainly tempted to view the former as being the
simpler one of them since it does not seem to
involve infinite-dimensional calculus. Yet appearances are
deceiving in this case: Since we 
have to view $U(\Om)$ as
a manifold modelled over $\ca_0(\Om)\times\Om$ the only legitimate
way of interpreting $(D_i^CR^C)(\varphi,x)=(\pa_iR^C)(\varphi,x)$
is to push forward the curve $t\mapsto(\varphi,x+te_i)$ via $T$ to
$\ca_0(\Om)\times\Om$ and to study the directional derivative
of $R^C\circ T^{-1}$ along $c:t\mapsto T(\varphi,x+te_i)=
(\varphi(\,.\,-(x+te_i)),x+te_i)$ at $t=0$. 

To this end, first note that for small absolute values of $t$, $c$ actually
takes values in $\ca_0(\Om)\times\Om$.
Moreover, $t\mapsto c(t)$ is a smooth curve in $\ca_0(\Om)\times\Om$
with respect to $\tau_0$, for the time being, according to
Proposition
\rf{smoothmaps}.
Since $c$ maps some interval
$[-\de,+\de]$ into $\ca_{0,K}(\Om)\times\Om$ for a suitable
$K\subset\subset\Om$, the restriction of $c$ to $(-\de,+\de)$
is smooth even with respect to $\tau_\Om$.
Therefore, the directional
derivative of $R^C\circ T^{-1}$ along $c:t\mapsto T(\varphi,x+te_i)=
(\varphi(\,.\,-(x+te_i)),x+te_i)$ at $t=0$ exists. Having
established existence, we can calculate its value as being
given by
$$\lim\limits_{t\to0}\frac{1}{t}
  [R\circ T^{-1}(c(t))-R\circ T^{-1}(c(0))]=
  \lim\limits_{t\to0}\frac{1}{t}
  [R(\vphi,x+te_i)-R(\vphi,x)].$$
Thus the usual formula works for $R^C\in\cc^\infty(U(\Om))$
and $D_i^C=\pa_i$, although $U(\Om)$ is not a linear space.

Finally, to see that $\pa_i R^C$ is again smooth, we have to convince
ourselves that
$(\pa_i R^C)\circ T^{-1}=(T^{-1})^*D_i^CR^C=D_i^J(T^{-1})^*R^C=
D_i^J(R^C\circ T^{-1})$ is smooth on $\ca_0(\Om)\times\Om$.
Since, by definition, $R^J:=R^C\circ T^{-1}$ is smooth on
$\ca_0(\Om)\times\Om$, so are its differential
$\rmd_1(R^C\circ T^{-1})$ and its partial derivative
$\pa_i(R^C\circ T^{-1})$ on their respective domains.
By Definition \rf{deriv}, $D_i^J(R^C\circ T^{-1})=(D_i^CR^C)\circ T^{-1}$
is smooth which is equivalent to the smoothness of $D_i^CR^C$ on
$U(\Om)$.
The smoothness of $D_i^JR^J$ for given $R^J\in\ce^J(\Om)$,
on the other hand, is immediate solely by the last part of the
argument given above.

Let us return to studying the sets $U_\eps(\Om)$.
For the purpose of reference, the following observation is
formulated as a lemma.
                    
\blem\lb{slL}
Let $K\subset\subset L\subset\subset\R^s$ 
and let $B$ be a subset of $\cd(\R^s)$ such that all $\varphi\in B$ have
their supports contained in some
bounded set. Then there exists $\eta>0$
such that $\supp S_\eps(\varphi)\subseteq L-x$
for all $\eps\le\eta$ and $\varphi\in B$, $x\in K$.
\et                      
\pr
Set $h:=\dist(K,\pa L)$. Then for each $x\in K$, $L$
contains the closed ball $\overline{B}_h(x)$ of radius $h$
around $x$.
If, on the other hand, the compact ball 
$D:=\overline{B}_r(0)$ contains the supports of all $\vphi\in B$
then putting $\eta:=\frac{h}{r}$ will do:
We have $\supp S_\eps(\varphi)+x\subseteq\eps D+x\subseteq L$
for $\eps\le\eta$, $\varphi\in B$, $x\in K$.
\ep

\bp \lb{slU}
Let 
$K\subset\subset L\subset\subset\Om$
and let $B$ be a subset of $\ca_0(\R^s)$ such
that all $\vphi\in B$
have their support contained in some fixed bounded subset of
$\R^s$.
Then there exists $\eta>0$
such that $B\times K\subseteq U_{\eps,L}(\Om)$
for all $\eps\le\eta$. In particular,
$B\times K\subseteq U_\eps(\Om)$
and the restrictions of
$\tau_{1,\eps}$ and $\tau_{2,\eps}$ to $B\times K$
are equal.
\et
\pr
$L$, $K$ and $B$ satisfying the assumptions of Lemma
\rf{slL}, we obtain $\eta>0$
such that $\supp S_\eps(\varphi)\subseteq L-x$,
i.e., $(\vphi,x)\in U_{\eps,L}(\Om)$
for all $\eps\le\eta$ and $\varphi\in B$, $x\in K$.
\ep

The fact that for small $\eps$ the topologies $\tau_{1,\eps}$
and $\tau_{2,\eps}$ agree on
sets of the form $B\times K$ as 
above is crucial to get the smoothness
properties right when it comes to testing for moderateness resp.\
negligibility, as we will see.

With these prerequisites at hand, we now are ready to introduce
the tangent space of $U_\eps(\Om)$ and to define differentials
of all orders of a smooth function defined on $U_\eps(\Om)$.
From an abstract point of view,
the tangent space of $U_\eps(\Om)$ with respect to $\tau_{2,\eps}$
at the point $(\vphi,x)\in U_\eps(\Om)$ is isomorphic to
$\ca_{00}(\Om)\times\Om$; to the tangent vector
\linebreak
$(\si,v)\in
\ca_{00}(\Om)\times\R^s$ at $(\rho,x)\in \ca_{0}(\Om)\times\Om$
there corresponds the ``tangent vector''
$(S_{\frac{1}{\eps}}T_{-x}(\si+\rmd\rho\cdot v),v)\in
\ca_{00}(\R^s)\times\R^s$ at
$(\vphi,x)=(T\circ S^{(\eps)})^{-1}(\rho,x)=
(S_\frac{1}{\eps}T_{-x}\rho,x)\in
U_\eps(\Om)$ where $\rmd\rho\cdot v$  denotes the directional
derivative of $\rho$ with respect to $v$.
The preceding formula is obtained by taking the derivative at
$t=0$ of the smooth curve $t\mapsto (T\circ S^{(\eps)})^{-1}(\rho+t\si,x+tv)$.
In this sense, the tangent space
to $U_\eps(\Om)$ at $(\vphi,x)\in U_\eps(\Om)$ can be identified
with the set of all $(\psi,v)\in\ca_{00}(\R^s)\times\R^s$
satisfying $\supp\psi\subseteq\frac{\Om-x}{\eps}$.
Note that in this case the kinematic tangent space
coincides with the operational one (the space of derivations defined on the 
smooth functions): In fact, by \cite{KM}, 28.7., and \cite{F}, this is true 
for the space $\cd (\Om)$ (more generally, for smooth sections with compact 
support of vector bundles over a manifold) and hence by \cite{F} for its 
complemented subspace $\ca _0(\Om)$.

Essentially, Proposition \rf{slU} also applies to tangent vectors:
\bp\lb{slUtang}
Let $K\subset\subset L\subset\subset\Om$ and let $B,C$ be subsets
of $\ca_0(\R^s)$ resp.\ $\ca_{00}(\R^s)$ such that all $\om\in B\cup
C$ have their supports contained in a fixed bounded subset of
$\R^s$.
Then there exists $\eta>0$
such that 
$B\times K\subseteq U_{\eps,L}(\Om)$
and $C\times\R^s$ is contained in the tangent space to
$U_\eps(\Om)$ at $(\vphi,x)$ for all
$(\vphi,x)\in B\times K$.
\et
The proof is virtually the same as for Proposition \rf{slU},
with $B$ now replaced by $B\cup C$; it even yields $\supp\psi
\subseteq\frac{L-x}{\eps}$ for all tangent vectors
$(\psi,v)$ with $\psi\in C$.

Now let $f:U_\eps(\Om)\to\C$ be a function
which is smooth with respect to $\tau_{2,\eps}$. 
Basically, $\rmd^n\!f$ 
ought to be defined on the
$n$-fold tangent space to $U_\eps(\Om)$, that is, on
$$\ct^n U_\eps(\Om):=\bigsqcup_{(\vphi,x)\in U_\eps(\Om)}
  \{(\vphi,x)\}\times\{(\psi,v)\in\ca_{00}(\R^s)\times\R^s
  \mid\supp\psi\subseteq\frac{\Om-x}{\eps}\}^n.
$$
$f$ being assumed as smooth with respect to $\tau_{2,\eps}$
by definition, we cannot use {\it a priori} the structure of the
surrounding space $\ca_0(\R^s)\times\Om$ to define $\rmd^n\!f$.
Instead,
we will decompose $U_\eps(\Om)$ (which has to be viewed as a manifold
modelled over $\ca_0(\Om)\times\Om$)
into a family of subsets which is 
characteristic for smoothness of a function with respect to
$\tau_{2,\eps}$ in the sense that $f$ is smooth on $U_\eps(\Om)$ if
and only if the restriction of $f$ to any of these subsets is
smooth, yet this time---due to equality of $\tau_{1,\eps}$
and $\tau_{2,\eps}$ on each of these subsets---either with
respect to $\tau_{1,\eps}$ or $\tau_{2,\eps}$. This allows the
calculus of $\ca_0(\R^s)\times \Om$ to be applied to $f$.
In particular, differentials of $f$ of any order can be computed
already from the restrictions of $f$ to these subsets; the chain
rule holds.

For the following, fix $\eps\in I$. We will simply write $S$ in
place of $S^{(\eps)}$.
\bp\lb{decompueps}
For given $x\in\Om$ and $L\subset\subset\Om$, set
$K:=K_{x,L}:=\frac{L-x}{\eps}(\subset\subset\frac{\Om-x}{\eps})$;
choose a compact set $M:=M_L$ and a positive number $h:=h_L$
such that $L\subset\subset M\subset\subset\Om$
and $0<h<\dist(L,\pa M)$.
Finally, set $U:=U_{x,L}:=B_h(x)$. Then
\bea\lb{ctnu}
\ct^nU_\eps(\Om)= \!
   \bigcup\limits_{x\in\Om \atop L\subset\subset\Om}
   \! \ca_{0,K}(\R^s)\times U\times\big(\ca_{00,K}(\R^s)
   \times\R^s\big)^n = \!
   \bigcup\limits_{x\in\Om \atop L\subset\subset\Om}
   \! \ct^n\big(\ca_{0,K}(\R^s)\times U\big).
\eea

In particular,
$
U_\eps(\Om)=
   \bigcup\limits_{x\in\Om\atop L\subset\subset\Om}
   \ca_{0,K}(\R^s)\times U
$.
Moreover, each set $\ca_{0,K}(\R^s)\times U$ is contained in
the corresponding set $U_{\eps,M}(\Om)$.
\et
\pr
For $y\in U$, we have
$K=\frac{L+(y-x)-y}{\eps}\subseteq\frac{M-y}{\eps}$.
Consequently, $\ca_{0,K}(\R^s)\times U$ is a subset of
$U_{\eps,M}(\Om)(\subseteq U_\eps(\Om))$
and any $(\psi,v)\in\ca_{00,K}(\R^s)\times\R^s$
belongs to the tangent space of $U_\eps(\Om)$ at $(\vphi,y)$
for arbitrary $(\vphi,y)\in\ca_{0,K}(\R^s)\times U$.
Conversely, if $(\vphi,x)\in U_\eps(\Om)$ and vectors $(\psi_i,v_i)$
($i=1\dots,n$) tangent to $U_\eps(\Om)$ at $(\vphi,x)$ are given, set
$L:=x+\eps\cdot(\supp\vphi\,\cup\bigcup\limits_{i=0}^n\supp\psi_i)$.
Noting that
$L\subset\subset\Om$, we obtain
$(\vphi,x)\in\ca_{0,K_{x,L}}(\R^s)\times U_{x,L}$
and $(\psi_i,v_i)\in\ca_{00,K_{x,L}}(\R^s)\times\R^s$.
This establishes (\rf{ctnu}). The last statement of the proposition
has already been shown above.
\ep

\bt\lb{smoothsub}
Let $f:U_\eps(\Om)\to\C$. $f$ is smooth
(with respect to $\tau_{2,\eps}$) if and only if the restriction
of $f$ to every set $\ca_{0,H}(\R^s)\times V$ is smooth with
respect to $\tau_{1,\eps}$ where $H\subset\subset\R^s$, $V$ is an
open subset of $\Om$ and $\ca_{0,H}(\R^s)\times V$
is contained in $U_{\eps,N}(\Om)$ for some $N\subset\subset\Om$.
\et
\pr
To begin with, note that it follows from
$\ca_{0,H}(\R^s)\times V
\subseteq U_{\eps,N}(\Om)$ that $\tau_{1,\eps}$ and
$\tau_{2,\eps}$ agree on $\ca_{0,H}(\R^s)\times V$.
Assuming $f$ to be smooth with respect to $\tau_{2,\eps}$
it is now clear that its restriction to any set
$\ca_{0,H}(\R^s)\times V$ is smooth with respect to
$\tau_{1,\eps}$.

Conversely, suppose that $f$ is smooth with respect to
$\tau_{1,\eps}$ on any set $\ca_{0,K}(\R^s)\times U$
as defined in Proposition \rf{decompueps}.
Let $L_1\subset\subset\Om$ and $x\in\Om$.
In a first step, we show that there exists an open neighborhood
$V_x$ of $x$ such that
$(T\circ S)^{-1}(\ca_{0,L_1}(\Om)\times V_x)
   \subseteq\ca_{0,K}(\R^s)\times U$
for some $K=K_{x,L}$, $U=U_{x,L}$:
Choose $L$ as to satisfy
$L_1\subset\subset L\subset\subset\Om$  and
define $K$, $M$, $h$, $U$ as it had been done in Proposition
\rf{decompueps}; finally, set $V_x:=B_r(x)$
where $r:=\min(h,\dist(L_1,\pa L))$.
For $\vphi\in\ca_{0,L_1}(\Om)$ and $y\in V_x$
it now follows that $y\in U$ and
$\supp S_{\frac{1}{\eps}}T_{-y}\varphi=
  \frac{1}{\eps}\cdot(\supp\vphi-y)\subseteq\frac{L_1-y}{\eps}
  =\frac{L_1+(x-y)-x}{\eps}\subseteq\frac{L-x}{\eps}=K$.
Altogether, $(T\circ S)^{-1}(\vphi,y)\in\ca_{0,K}(\R^s)\times U$.
By assumption and due to the inclusion relation just shown,
$f\circ(T\circ S)^{-1}$ is smooth on every set
$\ca_{0,L_1}(\Om)\times V_x$. Since $L_1\subset\subset\Om$
and $x\in\Om$ have been arbitrary, $f\circ(T\circ S)^{-1}$
is smooth on the whole of $\ca_0(\Om)\times \Om$
according to (an obvious modification of) Theorem
\rf{smoothonD}. By definition, $f$ is smooth.
\ep

An inspection of the preceding proof shows that
it is even sufficient for the smoothness of $f$ that the
condition stated in the theorem is satisfied for all
$H=K_{x,L}$, $V=U_{x,L}$ where
$H=K_{x,L}$, $V=U_{x,L}$ are defined as in Proposition
\rf{decompueps}.

Theorem \rf{smoothsub} allows us to define $\rmd^n\!f$
on every set $\ca_{0,H}(\R^s)\times V$ as above
which is contained in some set of the form
$U_{\eps,N}(\Om)$ with $N\subset\subset\Om$.
Since the sets $\ca_{0,H}(\R^s)\times V$ cover $U_\eps(\Om)$,
the differentials of $f$ are defined globally on $U_\eps(\Om)$.
Note that, in fact, they are well-defined:
Let  $(\vphi,x);(\psi_1,v_1),\dots,(\psi_n,v_n)$
be a set of data with 
\bea\lb{datadnr}
(\vphi,x)\in U_\eps(\Om),\ 
v_i\in\R^s,\ \psi_i\in\ca_{00}(\R^s),
\ \supp\psi_i\subseteq\frac{\Om-x}{\eps}\quad(i=1,\dots,n),
\eea
which is a member of
$\ct^n\big(\ca_{0,H_j}(\R^s)\times V_j\big)$ both for $j=0$
and $1$; then either way of restricting $f$ to
$\ca_{0,H_j}(\R^s)\times V_j$ 
gives the same value
for $\rmd^n\!f$ evaluated at these data
as the restriction of $f$ to
$\ca_{0,H}(\R^s)\times V$ would produce,
where $H:=H_1\cap H_2$ and $V:=V_1\cap V_2$.
In the particular case where $f$ is the
restriction of some $\ti f\in\cc^\infty(\ca_0(\R^s)\times\Om)$
to $U_\eps(\Om)$ then, of course, the differentials of $f$ will
agree with the restriction of the differentials of $\ti f$ to the
corresponding tangent spaces.

Now it is clear that the chain rule holds
for any composition of the form $f\circ\Phi$ where
$\Phi:W\to
U_\eps(\Om)$ is a
smooth function (no matter if with respect to $\tau_{1,\eps}$
or $\tau_{2,\eps}$) such that the domain $W$ of $\Phi$
can be covered by a family $W_\io$ of open sets with the property
that for each $\io$, $\Phi(W_\io)$ is a subset of a suitable
set $\ca_{0,H}(\R^s)\times V$ being, in turn, contained in
some $U_{\eps,N}(\Om)$. This is exactly the situation we will meet
in section \rf{jelshort} when constructing
the diffeomorphism invariant Colombeau algebra.

\brem \label{epsfree} In this final remark we drop the assumption of
$\eps\in I$ being fixed:
We demonstrate that
for $R_\eps:=R\circ S^{(\eps)}$
with $R\in\ce^C(\Om)=\cc^\infty(U(\Om))$ and
for a given family of sets of data
$(\vphi,x);(\psi_1,v_1),\dots,(\psi_n,v_n)$ for which
the supports of all $\vphi\in\ca_0(\R^s)$
and $\psi_1,\dots,\psi_n\in\ca_{00}(\R^s)$
occurring in this family are contained in a fixed
bounded set and $x$ ranges over some compact subset of
$\Om$ ($v_i$ being arbitrary from $\R^s$),
$\rmd^n\!R_\eps(\vphi,x)((\psi_1,v_1),\dots(\psi_n,v_n))$
is defined for all suffciently small $\eps$.
This will be the typical situation in the applications which are to come
along in the sequel.
To this end, let a subset $B$ of $\cd(\R^s)$ be given such that all
$\om\in B$ have their support contained in some fixed bounded set,
say, a closed ball $D:=\clb r(0)$. Choose $L$ satisfying
$K\subset\subset L\subset\subset\Om$ and let
$0<\eps_0<\frac{1}{r}\dist(L,\pa\Om)$.
For $\eps\le\eps_0$, set $M:=L+\clb{r\eps}(0)$; we have
$M\subset\subset\Om$. Then $\ca_{0,D}(\R^s)\times L^\circ
\subseteq U_{\eps,M}(\Om)$ since $\supp\vphi\subseteq D$
and $x\in L$ imply $\supp T_xS_\eps\vphi\in\clb{r\eps}(x)
\subseteq M$.
Consequently, $R_\eps$ is smooth with respect to $\tau_{1,\eps}$
on $\ca_{0,D}(\R^s)\times L^\circ$.
Moreover, $\ca_{00,D}(\R^s)$ is a subset of the
tangent space of $U_\eps(\Om)$ at
$(\vphi,x)\in\ca_{0,D}(\R^s)\times L$.
Hence we conclude that for all $\eps\le\eps_0$,
$\rmd^n\!R_\eps(\vphi,x)((\psi_1,v_1),\dots(\psi_n,v_n))$
is defined for all $x\in K$
(or even $x\in L^\circ$), $\varphi\in B\cap\ca_0(\R^s)$,
$\psi_1,\dots,\psi_n\in B\cap\ca_{00}(\R^s)$ and
$v_1,\dots,v_n\in\R^s$.
\et

\section{Construction of a diffeomorphism invariant \\ Colombeau algebra}\lb{jelshort}
The aim of this section is to complete Jel\'\i nek's approach to constructing
a diffeomorphism invariant Colombeau algebra, guided by the blueprint
sketched in section \rf{scheme1}. Contrary to \cite{JEL}, we base our
presentation on the C-formalism---for the convenience of those
readers who are acquainted best with the notation used in \cite{c2}
and \cite{CM}. Nevertheless, at each stage it should be possible
without difficulty to switch to the J-formalism of \cite{JEL}, due
to the translation apparatus described in section \rf{basics}.

Apart from closing a gap in Jel\'\i nek's construction we
supply those parts of the respective arguments which have
been included neither in \cite{CM} nor in \cite{JEL}, yet
which---according to our view---deserve due attention to be given.
This applies, in particular, to {\bf (T4)}, {\bf (T5)} and
{\bf (T6)}.
We also include a counterexample showing that, in
particular, the question of smoothness of the objects involved in
the construction in fact requires a careful treatment.
\subsection{The basis for the definition of the algebra}

The ``basic space'' $\ce^C(\Om)=\cc^\infty(U(\Om))$
and the embedding $\io^C:\cd'(\Om)\to\ce^C(\Om)$
have already been introduced in Definition~\rf{io}.
To complete {\bf(D1)}, it remains to introduce $\si$:
\bd\lb{defsi}
Let $\si:\cc^\infty(\Om)\to\ce^C(\Om)$ be the map defined by
$$(\si(f))(\vphi,x):=f(x)\qquad\qquad
  \mbox{\rm(}f\in\cc^\infty(\Om),\ (\vphi,x)\in U(\Om)\mbox{\rm)}.$$
\et
Also, ({\bf (D2)}) has already been taken care of in Definition~\rf{deriv}.
%
%
In the next step, we define the subspaces of moderate and
negligible members of $\ce^C(\Om)$, respectively.
From now on we will make free use of the convention that each (in)equality (E)
involving $R(S_\eps\varphi,x)$ with any arguments in place of $\varphi$ and $x$
is to be understood as 

``$R(S_\eps\varphi,x)$  is defined [i.e., $(\varphi,x)\in
U_\eps(\Om)$, i.e., $\supp\frac{1}{\eps^s}
   \varphi\!\left(\frac{.}{\eps}\right)
   \subseteq\Om -x$] and (E) holds''.

\bd\lb{defmod} \mbox{\rm{\bf (D3)} (\cite{JEL}, 8.)}
Let $R\in\ce^C(\Om)$. $R$ is called {\rm moderate} if the following
condition is satisfied:
\beas\forall K\subset\subset\Om\ \forall\al\in\N_0^s\ \exists N\in\N\ 
\forall \phi\in\cc^\infty_b( I\times\Om,\ca_0(\R^s)):
\eeas
$$\pa^\al(R(S_\eps\phi(\eps,x),x))=O(\eps^{-N})\qquad\qquad (\eps\to0)$$
uniformly for $x\in K$.
The set of all moderate elements of $R\in\ce^C(\Om)$
will be denoted by $\ce^C_M(\Om)$.
\et
There are several (mutually equivalent) variants of the above condition
defining moderateness;
six of them are listed below
in Theorem \rf{a--z}. The formulation in
Definition \rf{defmod} is condition (C) in
that theorem. Actually, Jel\'\i nek has chosen condition (A) of
Theorem \rf{a--z} for defining moderateness in \cite{JEL}, 8.

\bd\lb{defnegl} \mbox{\rm{\bf (D4)} (\cite{JEL}, 18.)}\\
Let $R\in\ce_M^C(\Om)$. $R$ is called {\rm negligible} if the following
condition (which, following \cite{JEL}, will be denoted by
$(3^\circ)$) is satisfied:
\beas\forall K\subset\subset\Om\ \forall\al\in\N_0^s\ \forall n\in\N\ 
   \exists q\in\N\ 
\forall \phi\in\cc^\infty_b( I\times\Om,\ca_q(\R^s)):
\eeas
$$\pa^\al(R(S_\eps\phi(\eps,x),x))=O(\eps^n)\qquad\qquad (\eps\to0)$$
uniformly for $x\in K$.
The set of all negligible elements of $R\in\ce^C(\Om)$
will be denoted by $\cn^C(\Om)$.
\et
Observe that in the preceding definition $R$ is presupposed to be
moderate in the sense of Definition \rf{defmod}.

Also for the condition given in Definition \rf{defnegl}
(without assuming $R$ to be moderate)
there are several 
equivalent
reformulations (see Theorem \rf{a--z/vm}). 

At this point it might be useful to observe that in the framework
of the J-formalism, the tests of {\bf (D3)} and {\bf (D4)}
have to be performed with
$R(T_xS_\eps\phi(\eps,x),x)$ in place of $R(S_\eps\phi(\eps,x),x)$,
due to $T^*$ being the appropriate bijection between
$\ce^J(\Om)$ and $\ce^C(\Om)$: If $R^C=T^*R^J$ then
$R^C(S_\eps\phi(\eps,x),x)=R^J(T_xS_\eps\phi(\eps,x),x)$;
thus $R^C$ is moderate resp.\ negligible in the C-frame if and only if
$R^J$ is moderate resp.\ negligible in the J-frame.

In order to make 
{\bf (D3)} and {\bf (D4)} meaningful,
$R(S_\eps\phi(\eps,x),x)$ has to depend in a smooth way on $x$.
This statement, though looking innocent at first glance, hides a
rather delicate question: The assumptions on $R$ resp.\ $\phi$
to be smooth refer to different topologies resp.\ bornologies.
In fact, as Example \rf{nonsmooth} shows, it can happen that
$R(S_\eps\phi(\eps,x),x)$ is not even locally bounded as a function
of $x$ for fixed, yet arbitrarily small $\eps$. This aspect has
been treated neither in \cite{CM} nor in \cite{JEL}. To put things
right, the following argument is needed:

Given $K\subset\subset\Om$,
choose $L$ such that $K\subset\subset L\subset\subset\Om$.
From the boundedness of $\phi$
we conclude that all $\phi(\de,x)$ with $\de\in I$,
$x\in L$ have their support contained in a suitable
fixed bounded subset of $\R^s$. According to
Proposition \rf{slU} there exists $\eta>0$ such that for all
$\eps\le\eta$,
$W:=\{\phi(\de,x)\mid\de\in I,x\in L\}\times L$
is a subset of $U_\eps(\Om)$
and the respective restrictions of
$\tau_{1,\eps}$ and $\tau_{2,\eps}$
to $W$ are equal.
Consequently, also the restrictions of
$\tau_1$ and $\tau_2$ to $S^{(\eps)}(W)=
\{S_\eps\phi(\de,x)\mid\de\in I,x\in L\}\times L$
agree. Let $L^\circ$ denote the interior of $L$.
On the set $(0,\eta)\times L^\circ$ (which is open in
$ I\times\Om$ and contains $(0,\frac{\eta}{2}]\times K$),
the map $(\eps,x)\mapsto (S_\eps\phi(\eps,x),x)$
is smooth with respect to $\tau_1$ by assumption,
hence also with respect to $\tau_2$. 
By definition, $R\in\cc^\infty(U(\Om))$ amounts to $R$ being smooth
with respect to $\tau_2$. 
Setting $\eps_0:=\frac{\eta}{2}$, we obtain that
$R(S_\eps\phi(\eps,x),x)$ is a smooth function
of $(\eps,x)$ on the open neighborhood
$(0,\eta)\times L^\circ$ of $(0,\eps_0]\times K$
which, finally, makes the test conditions in
{\bf (D3)} and {\bf (D4)} actually meaningful.

In this sense, we have to extend the convention
we made immediately preceding {\bf (D3)} by
requiring that whenever derivatives of a term like
$R(S_\eps\phi(\eps,x),x)$ on a set
$(0,\eps_0]\times K$ are under consideration,
it is to be understood that
$\eps_0$ is sufficiently small as to make sure
that $(\eps,x)\mapsto R(S_\eps\phi(\eps,x),x)$
is smooth on an open neighborhood of
$(0,\eps_0]\times K$.

In section \rf{scheme1}, definitions {\bf (D3)} and {\bf (D4)}
have been viewed as ``tests'' to be
performed on elements $R$ of $\ce^C(\Om)$, investigating their
behaviour on so-called ``test
objects'' $\phi(\eps,x)$. In the case at hand, the latter
take the form of smooth bounded (in the sense of section \rf{notterm})
maps from $ I\times\Om$ into $\ca_0(\R^s)$ for testing
moderateness of $R$ resp.\ into $\ca_q(\R^s)$ for testing
negligibility of $R$.

Returning to the exclusive use of the C-formalism,
we will drop the superscript $C$ in
$\io^C$, $D_i^C$, $\ce^C(\Om)$, $\cn^C(\Om)$ and $\ce^C_M(\Om)$
from now on. Moreover, 
note that in the sequel, by
$\pa_i=\frac{\pa}{\pa x_i}$ we will always denote the
corresponding derivative with respect to $x$, 
i.e., for example, 
$\pa_i\phi(\eps,x)=\frac{\pa}{\pa x_i}\phi(\eps,x)$
which must not be confused with $\frac{\pa}{\pa \xi_i}\phi(\eps,x)(\xi)$.

\bt \mbox{\rm{\bf (T1)}}\lb{th1} 
$$\begin{array}{rlrl}
        (i) &\io(\cd'(\Om))\subseteq\ce_M(\Om)\qquad
        &(ii)&\si(\cc^\infty(\Om))\subseteq\ce_M(\Om) \\
        (iii)&(\io-\si)(\cc^\infty(\Om))\subseteq\cn(\Om)\qquad
        &(iv)&\io(\cd'(\Om))\cap\cn(\Om)=\{0\}.
\end{array}$$
\et
\pr Since this theorem does not occur explicitly in \cite{JEL},
we include a proof; we will be more explicit on those aspects
which are new, compared to Colombeau algebras already
known, and more concise concerning the rest.  
To start with, let $u\in\cd'(\Om)$,
$\phi\in\cc^\infty_b( I\times\Om,\ca_0(\R^s))$ and let
$K\subset\subset L\subset\subset\Om$. By the boundedness
of $\phi$, there exists a bounded subset $C$ of $\R^s$ such that
$\spp\phi(\eps,x)\subseteq C$ for all $\eps\in I$, $x\in L$.
Consequently, for $x\in K$,
$\supp\pa^\al\phi(\eps,x)\subseteq C$
for all $\al\in\N_0^s$, $\eps\in I$.
Since for $\eps$ sufficiently small (say, for $\eps\le\eps_0$),
even $K+\eps C$ is contained in $L^\circ$, we obtain that
$\spp \phi(\eps,x)\left(\frac{.-x}{\eps}\right)
\subset\subset L^\circ$ for $\eps\le\eps_0$, $x\in K$.
Thus for the values taken by
$$\pa^\al
   ((\io u)(S_\eps\phi(\eps,x),x))
   =\left\lgl u,\pa^\al
   \frac{1}{\eps^s}\phi(\eps,x)\left(\frac{.-x}{\eps}\right)\right\rgl$$
on $(0,\eps_0]\times K$, only the restriction of $u$
to $L^\circ$ is relevant. Moreover, 
again by the boundedness of
$\phi$, 
each
$\pa^\al
   \frac{1}{\eps^s}\phi(\eps,x)\left(\frac{y-x}{\eps}\right)$
is of order at most $\eps^{-|\al|-s}$ as $\eps\to0$, uniformly for
$x\in K$, $y\in\R^s$. Finally, integrating the modulus of the latter
function over $\R^s$ with respect to $y$ gives values of order at
most $\eps^{-|\al|}$, uniformly for $x\in K$.

(i) 
Consider first the case 
$f \in \mathcal{C}(\Om)$.
Then for $\eps\le\eps_0$, 
$|\pa^\al((\io f)(S_\eps\phi(\eps,x),x))|$ is majorized by
$$\sup\limits_{L}|f|\cdot\int\limits_\Om\left|\pa^\al
   \frac{1}{\eps^s}\phi(\eps,x)\left(\frac{y-x}{\eps}\right)\right|\,dy
   =O(\eps^{-|\al|})
$$
uniformly for $x\in K$.
Since locally every distribution is a derivative of a suitable
continuous function, a similar estimate establishes the first
inclusion.

(ii) Let $f\in\cc^\infty(\Om)$. Then
$\pa^\al((\si f)(S_\eps\phi(\eps,x),x))=\pa^\al f(x)$
clearly is bounded on any compact set $K$.

(iii) Consider $(\io-\si)f$ for a given
$f\in\cc^\infty(\Om)$. Assume, in addition to the above, that
$\phi$ takes its values in $\ca_q(\R^s)$. Then 
for $\eps\le\eps_0$ and $x\in K$,
\beas&&\pa^\al(\io f-\si f)(S_\eps\phi(\eps,x),x))=
  \\&&\hphantom{\pa^\al(\io f-\si f)}
  \sum\limits_\beta {\al\choose\bet}\int\limits_{\frac{\Om-x}{\eps}}
  \left[(\pa^\bet f)(z\eps+x)-(\pa^\bet f)(x)\right]
  \pa^{\al-\bet}\phi(\eps,x)(z)\,dz.
\eeas
%
Taylor expansion of each $\pa^\bet f$ up to order 
$q$ yields that all terms containing a power of $\eps$ less or equal to $q$
vanish due to $\pa^{\al-\bet}\phi(\eps,x)\in\ca_q(\R^s)\cup
\ca_{q0}(\R^s)$. All the remainder terms are (smooth functions) of
order at most $\eps^{q+1}$, uniformly for $x\in K$, $z\in C$. Therefore,
$\pa^\al(\io f-\si f)(S_\eps\phi(\eps,x),x))=O(\eps^{q+1})$.

(iv) Suppose $\io u\in\cn(\Om)$ for some $u\in\cd'(\Om)$.
For $K\subset\subset\Om$
choose $q\in\N$ such that the condition in {\bf (D4)}
is satisfied for $\al=0$, $n=1$. Pick any $\vphi\in\ca_q(\R^s)$
and set $\phi(\eps,x):=\vphi$ for all $\eps\in I$, $x\in\Om$.
Then 
by the negligibility of $\io u$,
%
$(u*S_\eps\check\vphi)(x)=
  \lgl u,\eps^{-s}\vphi(\eps^{-1}(.-x))\rgl
  \to0$
as $\eps\to0$, uniformly on $K$. This shows that $u$, being the weak
limit of the smooth regular distributions $(u*S_\eps\check\vphi)$,
is equal to $0$.
\ep

Also, we immediately get
\bt \mbox{\rm{\bf (T2)} \lb{th2}(\cite{JEL}, 19.)}
$\ce_M(\Om)$ is a subalgebra of $\ce(\Om)$.
\et

\bt \mbox{\rm{\bf (T3)} \lb{th3}(\cite{JEL}, 19.)}
$\cn(\Om)$ is an ideal in $\ce_M(\Om)$.
\et

\subsection{The approach taken by J. Jel\'\i nek}
While the conditions given in Definitions \rf{defmod} and
\rf{defnegl} are adequate for proving {\bf (T1)}--{\bf (T3)},
we do need appropriate reformulations for establishing
the invariance of $\ce_M(\Om)$ and $\cn(\Om)$ under
differentiation ({\bf (T4)}, {\bf (T5)}) as well as under the
action induced by a diffeomorphism ({\bf (T6)}--{\bf (T8)}).
Suitable equivalent conditions allowing to prove {\bf (T4)} and {\bf
(T5)}, on the one hand, are given in Theorem 17 and in part
$(3^\circ)\!\Leftrightarrow\!(2^\circ)$ of Theorem 18 in \cite{JEL},
respectively. We will quote these theorems below as \rf{JT17} and
\rf{JT1823}.

To establish diffeomorphism invariance, on the other
hand, two problems
have to be coped with: First, transformed test objects in general are not
defined on the whole of $ I\times\Om$; 
secondly, the property
$\phi(\eps,x)\in\ca_q(\R^s)$ (as occurring in Definition
\rf{defnegl}) is not preserved under the action of a
diffeomorphism. The first of these aspects, though presenting
considerable
intricacies, is covered only by a few remarks in \cite{JEL} which,
in our view, do not provide a treatment as rigorous as these
questions require. The appropriate reformulations of Definitions
\rf{defmod} and \rf{defnegl} dealing with the poor domains of
transformed test objects are provided by
$(\mathrm{C})\!\Leftrightarrow\!(\mathrm{Z})$ of Theorem \rf{a--z} and
$(\mathrm{C}'')\!\Leftrightarrow\!(\mathrm{Z}'')$ of Corollary
\rf{a--z/avm}, respectively. In order to cope with the problem of
$\phi(\eps,x)\in\ca_q(\R^s)$ not being preserved by a
diffeomorphism, Jel\'\i nek claims in part
$(3^\circ)\!\Leftrightarrow\!(4^\circ)$ of \cite{JEL}, Theorem 18
that $R\in\ce_M(\Om)$ is negligible (condition $(3^\circ)$) if and
only if it passes the test on test objects $\phi$ having only
asymptotically vanishing moments of order $q$ on $K$
(condition $(4^\circ)$), as compared to
$\phi(\eps,x)\in\ca_q(\R^s)$ required by condition $(3^\circ)$. 
While $(4^\circ)\!\Rightarrow\!(3^\circ)$ is obvious, the
converse statement is not true (see Example \rf{excond4} below).
The error in the proof of $(3^\circ)\!\Rightarrow\!(4^\circ)$ consists
in passing from terms of the form
$\rmd_1\pa^\al[R_\eps(\dots)]$
to $[\rmd_1\pa^\al R_\eps](\dots)$ without applying the chain rule
with respect to the composition of $R_\eps$
with some ``inner'' function represented by the dots
(compare the proof of Theorem \rf{JT1834infty} given in section
\rf{atoz}). 

As a consequence, the construction of a diffeomorphism invariant Colombeau
algebra aimed at in \cite{JEL} is not complete in the following
sense: Eliminating condition $(4^\circ)$ from Theorem 18 deprives
one of the possibility of proving diffeomorphism invariance for the algebra
at hand. If, on
the other hand, $(4^\circ)$ is accepted as defining membership
in $\cn(\Om)$ (provided $R\in\ce_M(\Om)$) then the embedding of
$\cd'(\Om)$ into $\cg(\Om)$ does not preserve the product of smooth
functions (being considered as regular distributions) even in the
most simple cases, as can be
seen from part two of Example \rf{excond4} below. To overcome this difficulty,
we will present a substitute for condition $(4^\circ)$
(see Theorem \rf{JT1834infty} below) which in fact is
equivalent to $(3^\circ)$ under the assumption of moderateness
and, moreover, allows to deduce
diffeomorphism invariance of $\cn(\Om)$.
\bexs\lb{excond4}
(1) Let $\Om:=\R$ and denote by $u$ the regular distribution on $\R$
defined by $\lgl u,\vphi\rgl:=\int\xi\vphi(\xi)\,d\xi$
($\vphi\in\cd(\R)$).
According to part (iii) of Theorem \rf{th1}, $R:=\io u-\si u$
is a member of $\cn(\R)$, that is,
$R$ is moderate and satisfies the condition specified
in Definition \rf{defnegl}
which is condition $(3^\circ)$ of \cite{JEL}, Theorem 18.
We are going to show that $R$ in fact violates condition
$(4^\circ)$, thereby disproving
$(3^\circ)\!\Rightarrow\!(4^\circ)$.
It is immediate from the definitions that $R$
is given by $R(\vphi,x):=\int\xi\vphi(\xi)\,d\xi$
($\vphi\in\ca_0(\R)$, $x\in\R$).
Set K:=\{0\}, $\al:=1$ and $n:=2$.
For any given $q\in\N$, define a test object $\phi_q$ by
$\phi_q(\eps,x):=\vphi_q+x\cdot\psi_q$
($0<\eps\le1$, $x\in\R$) where $\vphi_q$ is an
arbitrary fixed member of $\ca_q(\R)$ and $\psi_q\in\ca_{00}(\R)$ is
chosen as to satisfy $\int\xi\psi_q(\xi)\,d\xi=1$. Then
$\phi_q$ belongs to $\cc^\infty_b( I\times\R,\ca_0(\R))$ and,
being equal to $\vphi_q$ on $K=\{0\}$, has asymptotically vanishing
moments of order $q$ on $K$, as required by condition $(4^\circ)$.
Yet 
$$\pa(R(S_\eps\phi_q(\eps,x),x))=\pa(\eps\cdot x)=\eps\neq
O(\eps^2),$$
no matter how large $q$ is chosen. This manifestly
contradicts condition $(4^\circ)$ for the choices of $K,\al,n$
made above. We also see that adopting
$(4^\circ)$ (together with moderateness, of course) as defining property for 
$\cn(\R)$ would invalidate part (iii) of Theorem {\bf (T1)}
which is the basis for $\io$ to preserve the product of smooth
functions. This is made explicit by the following item.

(2) Define $\phi_q$ as in Example (1), yet
this time requiring both $\int\xi\psi_q(\xi)\,d\xi=1$ and
$\int\xi^2\psi_q(\xi)\,d\xi=0$
for $\psi_q$ to be chosen from $\ca_{00}(\R)$ and, in addition,
$\vphi_q\in\ca_{\max(2,q)}(\R)$.
De\-no\-ting by $f$ the identity function
on $\R$, $f$ can be identified with the distribution $u$ introduced
previously. A straightforward calculation yields
$$\big(\io(f)\cdot\io(f)-\io(f\cdot f)\big)(\vphi,x)=
  \Big(\int\xi\vphi(\xi)\,d\xi\Big)^2-\int\xi^2\vphi(\xi)\,d\xi$$
where $\vphi\in\ca_0(\R)$, $x\in\R$. Substituting
$S_\eps\phi_q(\eps,x)$ for $\vphi$ and taking second derivatives
at $x=0$, we obtain
$$\pa^2\big((\io(f)\cdot\io(f)-\io(f\cdot f))
(S_\eps\phi_q(\eps,x),x)\big)\big|_{x=0}=
  \pa^2\big(\eps^2\cdot x^2\big)\big|_{x=0}=2\eps^2\neq O(\eps^3)$$
for any $q\in\N$. Since again $\phi_q$ is a test object having
asymptotically vanishing moments of order $q$ on $K=\{0\}$,
$\big(\io(f)\cdot\io(f)-\io(f\cdot f)\big)$
does not satisfy $(4^\circ)$, this time with respect to
$K=\{0\}$, $\al:=2$, $n:=3$.
Consequently, adopting $(4^\circ)$ in place of
$(3^\circ)$ as the defining property for negligibility would
prevent the restriction of $\io$ to $\cc^\infty(\R)$
from being an algebra homomorphism.
\et
To complete the prerequisites for establishing {\bf (T4)}--{\bf (T8)}
it remains to state the theorem replacing part
$(3^\circ)\!\Leftrightarrow\!(4^\circ)$ of Theorem 18 of \cite{JEL}.
To this end, we introduce the following terminology
(which, actually, is taken from the second paper of this series):
\bd\lb{avminfty}
Let $\phi\in\cc^\infty_b( I\times\R,\ca_0(\R))$,
$K\subset\subset\Om$, $q\in\N$.
$\phi$ is said to be of type $[\mathrm{A}_\mathrm{l}^\infty]_{K,q}$
if all derivatives $\pa_x^\bet\phi(\eps,x)$ ($\bet\in\N_0^s$)
have asymptotically vanishing moments of order $q$ on $K$.
\et
In the preceding definition,
``$\mathrm{A}$'', ``$\mathrm{l}$'' and ``$\infty$'' stand
for ``asymptotically
vanishing moments'', ``locally'' (i.e., only on the particular
compact set $K$ under consideration) and  ``derivatives of all
orders''.
\bt\lb{JT1834infty}
Let $R\in\ce_M(\Om)$. $R$ is negligible, i.e., $R$ satisfies the
condition specified in Definition \rf{defnegl} if and only if
it satisfies the following property (which will be referred to
as $(4^{\infty})$):
\beas\forall K\subset\subset\Om\ \forall\al\in\N_0^s\ \forall n\in\N\ 
   \exists q\in\N\ 
\forall \phi\in\cc^\infty_b( I\times\Om,\ca_0(\R^s)),\ 
   \phi\mbox{\rm\ of type\ }[\mathrm{A}_\mathrm{l}^\infty]_{K,q}:
\eeas
$$\pa^\al(R(S_\eps\phi(\eps,x),x))=O(\eps^n)\qquad\qquad (\eps\to0)$$
uniformly for $x\in K$.
\et
The proof of the preceding theorem is deferred to section
\rf{atoz}.
Restricting $\bet$ in the definition of
type $[\mathrm{A}_\mathrm{l}^\infty]_{K,q}$
to the value $0\in\N_0^s$ turns condition
$(4^\infty)$ into condition $(4^\circ)$ of Theorem 18
of \cite{JEL}.

\subsection{Stability under differentiation}
Having set up the prerequisites for the remaining
part of the construction in the previous subsection
we proceed to establish
Theorems {\bf (T4)} and {\bf (T5)}.
\bt \mbox{\rm{\bf (T4)}}\lb{th4}
For $R\in\ce_M(\Om)$, $\pa_i R\in\ce_M(\Om)$ $(i=1,\dots,s)$.
\et
\bt \mbox{\rm{\bf (T5)}}\lb{th5}
For $R\in\cn(\Om)$, $\pa_i R\in\cn(\Om)$ $(i=1,\dots,s)$.
\et
Curiously enough, the preceding Theorems
are not even mentioned in Jel\'\i nek's paper
\cite{JEL}. We regard them as highly non-trivial, however:
At first glance they might seem obvious since the respective
tests ask for a certain behaviour of {\it all} derivatives
$\pa^\al(R(\phi(\eps,x),x))$; thus, as one might be tempted to argue,
the moderateness or negligibility of $R$ implies the respective
property also for each $D_i^CR=\pa_iR$. This reasoning, however,
does not take into account the fact that the expression
$\pa^\al(R(\phi(\eps,x),x))$ used for testing $R$ involves
a certain sum of partial derivatives of $R$ multiplied by partial
derivatives of $\phi(\eps,x)$, according to the chain rule.
There is no easy relation between the respective expressions
for $\pa^\al$ and $\pa_i\pa^\al$ which could be used to
draw from the asymptoptic behaviour of the former to infer
the corresponding property for the latter.

The key tools for establishing Theorems {\bf (T4)} and {\bf (T5)}
are Jel\'\i nek's Theorem 17 and part
$(2^\circ)\!\Leftrightarrow\!(3^\circ)$ of Theorem 18 in \cite{JEL}.
For their ingenious proofs
we refer to the original \cite{JEL}.
We presume that the author was completely aware
of the r\^ole Theorems 17 and 18 had to play in that respect,
yet for some reasons he decided not to address this issue.

\bt \mbox{\rm (\cite{JEL}, Theorem 17)}\lb{JT17}
Let $R\in\ce(\Om)$. $R$ is moderate (i.e., a member of
$\ce_M(\Om)$) if and only if the following condition is satisfied:
\beas\forall K\subset\subset\Om\ \forall\al\in\N_0^s\ 
   \forall k\in\N_0\ \exists N\in\N
   \ \forall B\,(\mbox{bounded})\,\subseteq \cd(\R^s):\eeas
$$\pa^\al\rmd_1^k(R\circ S^{(\eps)})(\vphi,x)(\psi_1,\dots,\psi_k)=O(\eps^{-N})
   \qquad\qquad (\eps\to0)$$
uniformly for $x\in K$, $\vphi\in B\cap\ca_0(\R^s)$,
   $\psi_1,\dots,\psi_k\in B\cap\ca_{00}(\R^s)$.
\et
\bt \mbox{\rm (\cite{JEL}, Theorem 18,
       $(3^\circ)\!\Leftrightarrow\!(2^\circ)$)}\lb{JT1823}\\
Let $R\in\ce_M(\Om)$. $R$ is negligible (i.e., a member of
$\cn(\Om)$) if and only if the following condition is satisfied:
\beas\forall K\subset\subset\Om\ \forall\al\in\N_0^s\ 
   \forall k\in\N_0\ \forall n\in\N\ \exists q\in\N
   \ \forall B\,\mbox{(bounded)}\,\subseteq \cd(\R^s):\eeas
$$\pa^\al\rmd_1^k(R\circ
    S^{(\eps)})(\vphi,x)(\psi_1,\dots,\psi_k)=O(\eps^n)
   \qquad\qquad (\eps\to0)$$
uniformly for $x\in K$, $\vphi\in B\cap\ca_q(\R^s)$,
   $\psi_1,\dots,\psi_k\in B\cap\ca_{q0}(\R^s)$.
\et
Actually, the respective conditions in Definition \rf{defnegl} and
Theorem \rf{JT1823} (that is, conditions $(3^\circ)$ and $(2^\circ)$
of Theorem 18 of \cite{JEL}) are equivalent even without assuming
$R$ to be moderate: The proof is similar to the proof of Theorem 17
of \cite{JEL}, taking into account the equivalence of conditions
$(\mathrm{A}')$ and $(\mathrm{C}')$ of Theorem \rf{a--z/vm}.
 
It is a remarkable fact that in the condition of Theorem
\rf{JT1823}, $k$, $d_1^k$ and $\psi_1,\dots,\psi_k$ can be omitted
completely (i.e., only the case $k=0$ has to be taken into account)
without changing its content, provided $R$ is assumed to be moderate
(\cite{JEL}, Theorem 18, $(1^\circ)\!\Leftrightarrow\!(2^\circ)$).
We can interpret this  heuristically as the fact that the
moderateness condition takes care of the derivatives of $R$
to be small in the limit, provided only $R$ itself is
small in the appropriate sense.
A still stronger result will be shown in the second
paper of this series: Provided that $R$ is moderate, it
is even possible to omit the $x$-derivatives in the 
condition of Theorem \rf{JT1823}, yielding 
that $R\in\ce_M(\Om)$ is negligible
if and only if the following condition is satisfied:
\beas\forall K\subset\subset\Om\ \forall n\in\N\ \exists q\in\N
   \ \forall B\,\mbox{(bounded)}\,\subseteq \cd(\R^s):
   R(S_\eps\vphi,x)=O(\eps^n)
   \quad (\eps\to0) \eeas
uniformly for $x\in K$, $\vphi\in B\cap\ca_q(\R^s)$.

{\bf Proof of (T4) and (T5). }In order to
derive, for example, {\bf (T4)} from 
Theorem \rf{JT17}, assume $R\in\ce(\Om)$ to be moderate,
hence to satisfy the condition in \rf{JT17}.
Let $i\in\{1,\dots,s\}$;
due to $\pa_i(R\circ S^{(\eps)})=(\pa_iR)\circ S^{(\eps)}$
we obtain
$\pa^\al\rmd_1^k((\pa_iR)\circ S^{(\eps)})
  =\pa^{\al+e_i}\rmd_1^k(R\circ S^{(\eps)})$
where $e_i$ denotes the $i$-th standard unit vector in
$\R^s$. From the preceding equation it is immediate that
together with the differentials of $R$, also the differentials
of $\pa_iR$ are of order at most $\eps^{-N}$
for $\eps\to0$ in the appropriate sense.
Applying Theorem \rf{JT17} once more, we infer the moderateness of
$\pa_iR$. The proof of Theorem {\bf (T5)} proceeds
along the same lines, this time using Theorem \rf{JT1823}. \ep

\subsection{Diffeomorphism invariance}
For any diffeomorphism $\mu:\ti\Om\to\Om$ the requirements of {\bf (D5)} are
satisfied for 
$\bar\mu^C:U(\ti\Om)\to U(\Om)$ and
$\hat\mu^C:\ce^C(\Om)\to\ce^C(\ti\Om)$ as in Definitions \ref{barmu} and
\ref{hatmu}.

Again, we shall omit the superscript C from $\bar\mu^C$, $\hat\mu^C$.

Our next task is to establish the invariance of test objects
under the appropriate action induced by $\mu$.
This, of course, is at the very heart of the diffeomorphism
invariance of the Colombeau algebra to be constructed.
In the end, we
must be able to infer the moderateness of $\hat\mu R$
from the moderateness of $R$ (and, similarly, for
negligibility). 
Unfortunately, it need not be true in a strict sense that the class of
test objects $\phi(\eps,x)$ as in Definitions {\bf (D3)} resp.\
{\bf (D4)} is invariant under the action of a diffeomorphism:
The transformed test objects turn out to be defined not on the
whole of $ I\times\Om$ necessarily, i.e., they form a larger class than
the original test objects do. Due to $\hat\mu R=R\circ\bar\mu$, we
must start with the (formally stronger) assumption (let us denote
it by (Z))
that $R$ is moderate even with
respect to that larger class of test objects
to reach the conclusion that $\hat\mu R$ is moderate in the sense
of {\bf (D3)} (condition (C)). However, in Theorem \rf{a--z}
we will show that (C) and (Z) are, in fact,
equivalent so that, in the end, the property of moderateness
is shown to be preserved under the action of a diffeomorphism.

The following heuristic calculation clearly shows
which path is to be pursued: Let $\bar\mu_\eps$ be defined as in
section \rf{basics}.
For $\ti\phi\in\cc^\infty_b( I\times\ti\Om,\ca_0(\R^s))$ given,
we have to determine a function $\phi$
defined on a suitable subset of $ I\times\Om$
taking values in $\ca_0(\R^s)$
as to satisfy the following relation:
\beas
&& \hphantom{=} (\hat\mu R)(S_\eps\ti\phi(\eps,\ti x),\ti x)=
     R(\bar\mu(S_\eps\ti\phi(\eps,\ti x),\ti x))
     = R(\bar\mu S^{(\eps)}(\ti\phi(\eps,\ti x),\ti x))\\
&& = R(S^{(\eps)}(S^{(\eps)})^{-1}\bar\mu S^{(\eps)}(\ti\phi(\eps,\ti x),\ti
                x))
     = R(S^{(\eps)}\bar\mu_\eps(\ti\phi(\eps,\ti x),\ti x))
      = R(S_\eps\phi(\eps,\mu\ti x),\mu\ti x)
     \eeas

where $\phi$ is defined implicitly by the requirement of the last
equality to hold. Obviously, this is the case if and only if
$(\phi(\eps,x),x)=\bar\mu_\eps(\ti\phi(\eps,\mu^{-1}x),\mu^{-1}x))$
which, according to \rf{barmueps}, amounts to
(\ref{muepsphi}) in Theorem \ref{tht6} below.
To carry out the program outlined above,
three aspects of
$\phi$ have to be handled simultaneousely:
domain of definition, smoothness
and boundedness.

Starting with the first of these, observe
that the right hand side
of (\rf{muepsphi}) is only defined if
$\xi$ is an element of
$\frac{\Om-x}{\eps}$ whereas we would want
$\xi\mapsto\!\phi(\eps,x)(\xi)$ to be a test function on the whole of $\R^s$.
For the convenience of the reader, we include what essentially is
the argument in \cite{JEL}, Remark 25:
The right hand side of
(\rf{muepsphi})
(viewed as a smooth function on $\frac{\Om-x}{\eps}$)
can be extended to a smooth function
on the whole of $\R^s$ by setting it equal to 0 for
$\xi\notin\frac{\Om-x}{\eps}$, provided its support
is a compact subset of $\frac{\Om-x}{\eps}$.
This, in turn, is equivalent to
$\ti\phi(\eps,\mu^{-1}x)$ having compact support contained in
$\frac{\ti\Om-\mu^{-1}x}{\eps}$: Indeed,
$\xi\mapsto\frac{\mu^{-1}(\eps\xi+x)-\mu^{-1}x}{\eps}$
maps $\frac{\Om-x}{\eps}$ diffeomorphically onto
$\frac{\ti\Om-\mu^{-1}x}{\eps}$. Therefore, the largest
natural domain of definition for
$\phi$ is the set $D$ of all
$(\eps,x)\in I\times\Om$ for which 
$\spp\ti\phi(\eps,\mu^{-1}x)$ is contained in  
$\frac{\ti\Om-\mu^{-1}x}{\eps}$, i.e., for which
$(\ti\phi(\eps,\mu^{-1}x),\mu^{-1}x)\in U_\eps(\ti\Om)$.
We do not know the form of $D$ explicitly;
however, due to the boundedness of the map $\ti\phi$,
for each given $K\subset\subset\Om$
there exists $\eps_0>0$ (chosen appropriately
with respect to the compact set $\mu^{-1}K$ by Proposition
\rf{slU}) such that
$(\ti\phi(\eps,\mu^{-1}x),\mu^{-1}x)\in U_\eps(\ti\Om)$ for all
$(\eps,x)\in(0,\eps_0]\times K$ which amounts to
$(0,\eps_0]\times K$ being a subset of
$D$. Summarizing,
$\phi$ is defined at least on ``rectangles''
of the form $(0,\eps_0]\times K$ as a map
taking values in $\ca_0(\R^s)$.
This settles the problem
of the domain of $\phi$ in a satisfactory way,
as we will see shortly.

In the light of Example \rf{hilfe} as well as of Example
\rf{nonsmooth} at the end of this section,
it seems advisable to
give a careful treatment also to
the question of
smoothness of $\phi$. We defer this to the formal proof of
{\bf (T6)}.

\bt \label{tht6} \mbox{\rm{\bf (T6)}\lb{th6} (\cite{JEL},  25.)}
Let $\mu:\ti\Om\to\Om$ be a diffeomorphism.
Let $\ti\phi\in\cc^\infty_b( I\times\ti\Om,\ca_0(\R^s))$
and define 
$D(\subseteq I\times\Om$) by
$$D:=\{(\eps,x)\in I\times\Om\mid
   (\ti\phi(\eps,\mu^{-1}x),\mu^{-1}x)\in U_\eps(\ti\Om)\}.$$
For $(\eps,x)\in D$, set
$\phi(\eps,x):=\pro_1\bar\mu_\eps(\ti\phi(\eps,\mu^{-1}x),
\mu^{-1}x))$, i.e., 
\bea\lb{muepsphi}\phi(\eps,x)(\xi)&:=&
  \ti\phi(\eps,\mu^{-1}x)\left(\frac{\mu^{-1}(\eps\xi+x)-\mu^{-1}x}
   {\eps}\right)\cdot|\det D\mu^{-1}(\eps\xi+x)|.
\eea
Then $\phi$ satisfies the requirements 1) and 2) specified for test objects
in condition (Z) of Theorem \rf{a--z}.

If, in addition, all derivatives $\pa^\al_{\ti x}\ti\phi(\eps,\ti x)$
have asymptotically vanishing moments
of order $q$ on some compact subset $\ti L$ of $\ti\Om$
($\al\in\N_0^s$) then
all derivatives $\pa^\al_x\phi(\eps,x)$ of the
the function $\phi$ defined by (\rf{muepsphi})
have asymptotically vanishing moments of order
$\left[\frac{q+1}{2}\right]$ on the (compact) set $L=\mu(\ti L)$. 
\et

\pr
That $\phi$ is well-defined on $D$ has already been shown.
To establish the smoothness of $\phi$ on suitable open subsets of
$ I\times\Om$, expand $\bar\mu_\eps$ to obtain
$$(\phi(\eps,x),x)=(S^{(\eps)})^{-1}T^{-1}\bar\mu^J TS^{(\eps)}
       (\ti\phi(\eps,\mu^{-1}x),\mu^{-1}x)).$$
In a first step, we discuss the smoothness of
$\Phi(\eps,x):=TS^{(\eps)}
       (\ti\phi(\eps,\mu^{-1}x),\mu^{-1}x))$.
$\Phi$ involves the maps
$\mu^{-1}$, $\ti\phi$, $S$ and $T$, all of which are smooth by the results of
Proposition \rf{smoothmaps}, provided
$\ca_0(\R^s)$ is endowed with the natural locally convex topology
inherited from $\cd(\R^s)$.
Let $K\subset\subset\Om$;
we are going to show that for suitable $\eps_0>0$
and $\ti M\subset\subset\Om$,
$\Phi$ actually maps some open neighborhood of
$(0,\eps_0]\times K$ into ${\cal A}_{0,\tilde M}
(\tilde\Om)\times\tilde\Om$.
To this end, choose $L,M$ such that
$K\subset\subset L\subset\subset M\subset\subset\Om$.
Set $\ti L:=\mu^{-1}L$, $\ti M:=\mu^{-1}M$ and
$h:=\dist(L,\pa M)$, $\ti h:=\dist(\ti L,\pa\ti M)$.

Due to the
boundedness of $\ti\phi$,
there is $r\ge \tilde{h} \ (>0)$ such that the supports of all 
$\ti\phi(\eps,\ti x)$ with $\eps\in I$, $\ti x\in\ti L$
are contained in the closed ball $\clb r(0)$. 
Setting $\eta:=\frac{\ti h}{r}$, Proposition
\rf{slU} and a glance at the proof of Lemma \rf{slL}
show that
$\ti\phi( I\times\ti L)\times\ti L\subseteq U_{\eps,\ti M}
  (\ti\Om)$
for all $\eps\le\eta$.
In particular, for all $x\in L$ and $\eps\le\eta$,
\be\lb{tiphiUM}
(\ti\phi(\eps,\mu^{-1}x),\mu^{-1}x)\in U_{\eps,\ti M}(\Om).
\end{equation}
Therefore, $\Phi$ maps the open set $U:=(0,\eta)\times L^\circ$
into $\ca_{0,\ti M}(\ti\Om)\times\ti\Om$.
On the latter, however, the topologies
$\tau_0$ and $\tau_{\ti\Om}$
introduced in section \rf{basics} coincide.
From the smoothness of the restriction of $\Phi$ to $U$ with
respect to $\tau_0$ (which was established above) we conclude the
smoothness with respect to $\tau_{\ti\Om}$. Now we are ready to go on
with the proof of the smoothness of
$\phi=\pro_1\circ(S^{(\eps)})^{-1}\circ T^{-1}\circ\bar\mu^J\circ\Phi$,
observing that $\bar\mu^J$ is smooth
if the domain $\ca_0(\ti\Om)\times\ti\Om$ and the range space
$\ca_0(\Om)\times\Om$ carry the topologies $\tau_{\ti\Om}$ and
$\tau_\Om$, respectively.
(Note that in general, $\bar\mu^J$ is {\it not}
$\tau_0$-$\tau_0$-smooth as can be seen
from Example \rf{munonsmooth} below.)
Weakening this conclusion by replacing $\tau_\Om$
by $\tau_0$ on $\ca_0(\Om)\times \Om$ and using once more the smoothness
of $T$ and $S$ with respect to the usual topology of $\ca_0(\R^s)$,
we finally obtain that for $\eps_0:=\frac{1}{2}\eta$
(if $\eta=1$ we may choose $U:= I\times L^\circ$, being
open in $ I\times \Om$, and $\eps_0:=1$),
$\phi$ is smooth on the open neighborhood $U$ of
$(0,\eps_0]\times K$, as claimed by condition 1)
of Theorem \ref{a--z} (Z).

For the proof of boundedness of $\phi$, we extend the argument of
\cite{JEL}, 25., Proposition.
Note that, by (\rf{tiphiUM}) above, $\phi$ is defined at least on
$(0,\eta]\times L$. Let
$l:=\max(1,\sup\{\|D\mu_{\ti x}\|\mid\ti x\in\ti M\})$.
Then for $\ti x\in\ti L$, we have
\be\lb{lip}
  \|\ti x-\ti y\|\le\ti h\Rightarrow 
   \|\mu\ti x-\mu\ti y\|\le l\|\ti x-\ti y\|,
\end{equation}
due to $\clb{\ti h}(\ti x)\subseteq \ti M$.
$\clb r(0)$ containing the support of every $\tilde\phi(\eps,\ti x)$ for
$\eps\in I$, $\ti x\in\ti L$, we have
$$\supp T_{\mu^{-1}x}S_\eps\ti\phi(\eps,\mu^{-1}x)\subseteq
  \clb{r\eps}(\mu^{-1}x)\subseteq\ti M$$
for $\eps\le\eta$, $x\in L$. Applying $\bar\mu^J$
we obtain, by (\rf{lip}), 
$$\supp {\pro}_1\,\bar\mu^JTS^{(\eps)}(\ti\phi
  (\eps,\mu^{-1}x),\mu^{-1}x)\subseteq
  \clb{lr\eps}(x)\cap M$$
and, finally, $\supp\phi(\eps,x)\subseteq\clb{lr}(0)$
for $\eps\le\eta$, $x\in L$. It follows that for each
$\al\in\N_0^s$, $\supp\pa^\al\phi(\eps,x)\subseteq\clb{lr}(0)$
for $\eps<\eta$, $x\in L^\circ$.
For the boundedness of
$\{\pa^\al\phi(\eps,x)\mid(\eps,x)\in(0,\eta_1)\times L^\circ\}$
where $\eta_1:=\min(\eta,\frac{h}{rl})=
\min(\frac{\ti h}{r},\frac{h}{rl})$´ it now suffices to show that
for each fixed $\bet\in\N_0^s$,
$$\sup\{|\pa^\bet_\xi\pa^\al_x\phi(\eps,x)(\xi)|\,\big|\,
  \eps<\eta_1,\ x\in L,\ \|\xi\|\le lr\}$$
is finite.
For $\eps<\eta_1,\ x\in L^\circ,\ \|\xi\|\le lr$ 
(observe that $\supp\pa^\al_x\phi(\eps,x)\subseteq\clb{lr}(0)
\subseteq\frac{\Om-x}{\eps}$)
$\pa^\bet_\xi\pa^\al_x\phi(\eps,x)(\xi)$  is a sum of terms of the
form
\bea\lb{giga}
   \pa^{\bet_0}_{\ti\xi}((\pa^{\al_0}_{\ti x}\ti\phi)(\eps,\mu^{-1}x))
   \left(\frac{\mu^{-1}(\eps\xi+x)-\mu^{-1}x}
   {\eps}\right)
   \cdot g_0(x)\cdot g_1(\eps,x,\xi)\cdot\,\dots\,\cdot g_p(\eps,x,\xi)\cdot
   \nn\\ 
   \cdot\pa^{\bet'}_\xi\pa^{\al'}_x|\det
   D\mu^{-1}(\eps\xi+x)|\hphantom{www}
\eea
where $g_0$ is a certain product of derivatives of
components of $\mu^{-1}$
(hence bounded on $L$) and
each $g_j$ ($j=1,\dots,p$)
is some derivative of some component of
$\frac{\mu^{-1}(\eps\xi+x)-\mu^{-1}x}{\eps}$, i.e.,
$$g_j(\eps,x,\xi)=\pa^{\bet_j}_\xi\pa^{\al_j}_x
  \left(\frac{\mu^{-1}_{i_j}(\eps\xi+x)-\mu^{-1}_{i_j}x}{\eps}\right).$$
For $\eps\le\eta_1$, $x\in L$, $\|\xi\|\le lr$, we have
$x+\eps\xi\in\clb h(x)\subseteq M$. Thus the last factor
in (\rf{giga}) is uniformly bounded, as is
$\pa^{\bet_0}_{\ti\xi}((\pa^{\al_0}_{\ti x}\ti\phi)(\eps,\mu^{-1}x))
(\ti\xi)$ for $\eps\le1$, $x\in L$, $\ti\xi\in\R^s$. It remains to
discuss the boundedness of the factors $g_j$ ($j=1,\dots,p$).
For the sake of simplicity, we replace
$\al_j,\bet_j,i_j$ by $\al,\bet,i$. Considering first the
case $|\bet|>0$ (say,
$\bet_k\ge1$), the uniform boundedness of
$$\pa^\bet_\xi\pa^\al_x
  \frac{\mu^{-1}_i(\eps\xi+x)-\mu^{-1}_ix}{\eps}=
  \pa^{\bet-e_k}_\xi\pa^\al_x
  \left((\pa_k\mu^{-1}_i)(\eps\xi+x)\right)$$
on $\eps\le\eta_1$, $x\in L$, $\|\xi\|\le lr$ is evident.
If, on the other hand, $|\bet|=0$, choose a Lipschitz constant
$l_{\al,i}$ for $\pa^\al_x\mu^{-1}_i$ with respect to
$L,h$ (in the same way as $l$ was chosen for $\mu$ with respect to
$\ti L,\ti h$ as to satisfy (\rf{lip})). It follows that
$$\left|\frac{(\pa^\al_x\mu^{-1}_i)(\eps\xi+x)
   -(\pa^\al_x\mu^{-1}_i)(x)}{\eps}\right|\le
   \frac{l_{\al,i}\|\eps\xi\|}{\eps}\le l_{\al,i}lr$$
which establishes the uniform boundedness on
$\eps\le\eta_1$, $x\in L$, $\|\xi\|\le lr$
also of this term.
Replacing $\eps_0=\frac{1}{2}\eta$ by $\eps_0:=\frac{1}{2}\eta_1$,
we have shown altogether that $\phi$ is smooth and
each derivative $\pa^\al\phi$ is bounded on the open
neighborhood $(0,\eta_1)\times L^\circ$ ($\subseteq D$)
of $(0,\eps_0]\times K$, as required for satisfying
conditions 1) and 2) of Theorem \ref{a--z}, (Z).

Finally, assume that 
all derivatives $\pa^\al_{\ti x}\ti\phi(\eps,\ti x)$
have asymptotically vanishing moments
of order $q$ on some compact subset $\ti L$ of $\ti\Om$
($\al\in\N_0^s$). We have to show that
for all $\bet\in\N_0^s$ satisfying
$|\bet|\leq \left[ \frac{q+1}{2}\right]$ and for arbitrary
$\al\in\N_0^s$, 
\be\lb{th6avm}
  \langle\xi^\bet,\pa^\al_x\phi(\eps,x)(\xi)\rangle
  =\pa^\al_x\int\left(\frac{\mu(\eps\ti\xi+\ti x)-\mu\ti x}{\eps}\right)^\bet
                \ti\phi(\eps,\ti x)(\ti\xi)\,d\ti\xi
        =O(\eps^{\left[\frac{q+1}{2}\right]})
\end{equation}
uniformly for $x=\mu\ti x \in L$
(or $\ti x=\mu^{-1}x\in\ti L$, respectively).
Note that the preceding equation is meaningful since there exists
$\eps_0>0$ such that all the terms occurring therein are defined
for $\eps\le\eps_0$, $x\in L$, $\xi\in\R^s$, $\ti\xi\in\R^s$
resp.\ $\ti\xi$ ranging over a fixed compact set containing the
supports of all $\ti\phi(\eps,\ti x)$ with $\eps\le\eps_0$,
$\ti x\in\ti L$ in its interior.

Let us consider the case $\al=0$ first.
Expanding $\mu$ into a Taylor series up to order $n$ at $ \ti x$
we may write the first factor in the integral
as a finite sum of terms of the form
\be\lb{moment}
        \frac{1}{\eps^{|\bet|}}\,
        \frac{\partial^{\al_1}\mu_{i_1}
           (\ti x+\eta_{|\al_1|}
           \theta_{i_1}\eps\ti\xi)}{\al_1!}\,
           \ti\xi^{\al_1}\eps^{|\al_1|}
        \!\dots\!
        \frac{\partial^{\al_{|\bet|}}\mu_{i_{|\bet|}}
           (\ti x+\eta_{|\al_{|\bet|}|}
           \theta_{i_{|\bet|}}\eps\ti\xi)}{\al_{|\bet|}!}\,
           \ti\xi^{\al_{|\bet|}}\eps^{|\al_{|\bet|}|}
\end{equation}
where $1\leq|\al_j|\leq n+1$, $0<\theta_i<1$, $\eta_j=1$ if $j=n+1$
(i.e., if the respective factor is a remainder term of the Taylor series)
and $\eta_j=0$ otherwise.
Observe that in the present context,
by $\al_1,\dots,\al_{|\bet|}$ we are denoting
$|\bet|$ variables taking values in $\N_0^s$, yet {\it not}
components of a single variable $\al\in\N_0^s$.
Letting $\gamma:=\sum\limits_{j=1}^{|\bet|}
\al_j$,
the above expression (\rf{moment})
contains a factor $\ti\xi^\gamma\,\eps^{|\gamma|-|\bet|}$. 
(Note that since all $|\al_j|\geq 1$ we have $|\gamma|\geq |\bet|$.) If $|\gamma|-|\bet|
\geq\left[\frac{q+1}{2}\right]$ we are done with that particular
term, taking into account that
the integral has to be taken 
over a fixed compact set only. If,
on the other hand, $|\gamma|-|\bet|<\left[\frac{q+1}{2}\right]$ 
and all $\eta_j$ vanish we may use the assumption on $\ti\phi(\eps,\ti x)$ since 
in this case $|\gamma|<\left[\frac{q+1}{2}\right]+|\bet|\leq q+1$.
Finally, if there is at least one $\eta_j$ nonvanishing then at
least one $|\al_j|=n+1$,
implying $|\gamma|\geq n+|\bet|$. Hence, choosing $n\geq\left[\frac{q+1}{2}\right]$
completes the proof for the case $\al=0$.

To deal with the general case,
express the operator $\pa^\al_x$ occurring in (\rf{th6avm})
in terms of operators $\pa_{\ti x}^{\al'}$, according to the chain
rule. 
Now apart from certain partial derivatives of components of
$\mu^{-1}$ (which are bounded on $L$),
Leibniz' rule yields a sum of terms which are similar to
those considered above, with certain derivatives
$\pa_{\ti x}^{\al''}\pa^{\al_j}\mu_{i_j}$ and
$\pa_{\ti x}^{\al'''}\ti\phi$ replacing
$\pa^{\al_j}\mu_{i_j}$ and
$\ti\phi$, respectively. The powers of $\xi$ resp.\ $\eps$
remaining unchanged, the same reasoning as above establishes
(\rf{th6avm}) for arbitrary $\al\in\N_0^s$.
\ep

Note that the conclusion of the preceding theorem is also obtained if, instead
of $\tilde\phi \in \cc^\infty_b(I\times \tilde\Om, {\mathcal A}_0(\R^s))$, 
$\tilde\phi$ is only assumed to satisfy the analogs (for $\tilde\Om$) of
conditions 1) and 2) of \ref{th6}. Moreover, ${\mathcal A}_0(\R^s)$ can
be replaced by $\D(\R^s)$ throughout. 

Now {\bf (T7)} and {\bf (T8)} follow from {\bf (T6)} and {\bf (D5)},
due to the particular form of {\bf (D3)} and {\bf (D4)}:
Assuming, for
example, $R$ to be moderate, $R$ also satisfies condition (Z) of
Theorem \rf{a--z}. Given $\ti K\subset\subset\ti\Om$,
$\al\in\N_0^s$ and $\ti\phi\in 
\cc^\infty_b( I\times\ti\Om,\ca_0(\R^s))$,
define $\phi$ as in {\bf (T6)}.
According to the chain rule,
\beas\pa_{\ti x}^\al
     \left((\hat\mu R)(S_\eps\ti\phi(\eps,\ti x),\ti x)\right)&=&
      \pa_{\ti x}^\al
      \left(R(S_\eps\phi(\eps,\mu\ti x),\mu\ti x)\right)\\
      &=&\sum\limits_{\bet:\,|\,\bet|\le|\al|}
      \pa_{x}^\bet
      (R(S_\eps\phi(\eps,x),x))\Big|_{x=\mu\ti x}
      \cdot g_\bet(\ti x)
     \eeas
where each function $g_\bet$ is a certain sum of products of
partial derivatives of components of $\mu$, hence bounded on
$\ti K$. $R$ satisfying (Z) of Theorem \rf{a--z}, it follows
that for some $N\in\N$, $\pa_{\ti x}^\al
\left((\hat\mu R)(S_\eps\ti\phi(\eps,\ti x),\ti x)\right)=O(\eps^{-N})$
uniformly on $K$. This shows that also $\hat\mu R$ is
moderate.
If, on the other hand, $R$ is assumed to be negligible, $R$ even
passes the negligibility test on test objects $\phi$
being of type $[\mathrm{A}_\mathrm{l}^\infty]_{K,q}$,
according to
Theorem \rf{JT1834infty}.
Now a similar reasoning as in the case of moderateness, this time
using Corollary \rf{a--z/avm} in place of Theorem \rf{a--z}, establishes
the invariance of negligibility under the action
of a diffeomorphism. Thus we have shown
\bt \mbox{\rm{\bf (T7)} \lb{th7}(\cite{JEL}, 25.)}\\
$\ce_M$ is invariant under $\hat\mu$, i.e., $\hat\mu^C$ maps
                     $\ce_M(\Om)$ into $\ce_M(\ti\Om)$.
\et
\bt \mbox{\rm{\bf (T8)} \lb{th8}(\cite{JEL}, 25.)}\\
$\cn$ is invariant under $\hat\mu$, i.e., $\hat\mu^C$ maps
                     $\cn(\Om)$ into $\cn(\ti\Om)$.
\et
Having  completed  all the steps of the general construction scheme in
section  3  we finally reach the goal of this section, the definition of
the algebra itself:
\bd \mbox{\rm{\bf (D6)} (\cite{JEL}, 19.)}
$$\cg(\Om):=\ce_M(\Om)\big/\cn(\Om)$$
\et
Since the respective ideals of negligible functions are
invariant under  
$D_i:\ce_M(\Om)$  $\to$ $\ce_M(\Om)$ as well as under
$\hat\mu:\ce_M(\Om)\to\ce_M(\ti\Om)$, both these maps factorize 
via the respective quotients to yield maps (which we denote
by the same symbols)  $D_i:\cg(\Om)\to\cg(\Om)$ and
$\hat\mu:\cg(\Om)\to\cg(\ti\Om)$. This completes the (functorial)
construction of a differential algebra containing $\cd'(\Om)$
in such a way as to extend the usual product on $\cc^\infty(\Om)$.

If we had decided to perform this construction in the J-frame
we would have obtained objects $\cg^J(\Om)$ isomorphic to the
$\cg^C(\Om)$ above: Indeed, also $T^*$ factorizes via
quotients with respect to $\cn^J(\Om)$ resp.\ $\cn^C(\Om)$,
thereby inducing a bijection between the J- and C-variant
of the diffeomorphism invariant Colombeau algebra
at hand.

Next, we give an example of a test object
$\phi(\eps,x)$ in the
sense of {\bf (D3)} and a distribution $u\in\cd'(\Om)$ such that
on every strip
$((0,\eps_0]\times \Om)\cap D$, the map
$x\mapsto(\io u)(S_\eps\phi(\eps,x),x)$ is not even locally bounded
(hence, {\it a fortiori, }neither smooth) where $\eps_0\in I$
is arbitrary
and
$$D:=\{(\eps,x)\mid(\phi(\eps,x),x)\in U_\eps(\Om)\}
     =\{(\eps,x)\mid \spp S_\eps\phi(\eps,x)\subseteq\Om-x\}$$
is the natural maximal domain of definition of
$(\io u)(S_\eps\phi(\eps,x),x)$. This phenomenon is due to the
mismatch of the respective smoothness notions for $\io u$ and
$\phi$. It cannot occur on sets of the form $(0,\eps(K)]\times K$
where $K\subset\subset\Om$ and $\eps(K)$ is chosen
suitably with respect to $K$, according to the discussion
following {\bf (D4)} (Definition \rf{defnegl}).

\bex \lb{nonsmooth}
We employ the notation introduced in
Example \rf{hilfe}. In particular, by
$x,y$; $\xi,\eta$ we now denote coordinates of points
$(x,y),(\xi,\eta)\in\R^2$, $(\xi,\eta)$ replacing
the former $(x,y)$ due to the need for additional
variables $x,y$ in $\phi(x,y)(\xi,\eta)$.
Let $\Om,\psi,\vphi,u$ be defined as in Example \rf{hilfe}.
Choose a smooth non-decreasing function
$\nu:\R\to\R$ taking the constant value $k$ on each of the
intervals $I_k:=[k-\frac{1}{4},k+\frac{1}{4}]$
($k\in\Z$), respectively. Define
$\phi(\eps,x,y)(\xi,\eta):=
  \sin\pi y\cdot S_{|\nu(y)|}\psi(\xi,\eta+\nu(y))
  \ +\ \vphi(\xi,\eta).$
Then, obviously,
$\phi\in\cc^\infty_b( I\times\Om,\ca_0(\R^s))$ (note that also
$|\nu|$ is smooth). Letting $x:=0$, $y\in I_k\setminus\{k\}$,
$\eps:=\frac{1}{|k|}$ ($0\neq k\in\Z$), $(\xi,\eta)\in\Om$, we obtain
\beas
(\io u)(S_\eps\phi(\eps,x,y),x,y)&=&
   \sin\pi y\cdot\lgl u,S_\eps S_{|\nu(y)|}\psi(\xi,\eta+\nu(y)-y)\rgl
   \\&=&\sin\pi y\cdot\lgl u,\psi(\xi,\eta+k-y)\rgl.\eeas
Substituting $t:=y-k$
in the last expression
(note that $0<|t|\le\frac{1}{4}$) yields \linebreak
$(-1)^k\sin\pi t\cdot\lgl u,\psi(\xi,\eta-t)\rgl$,
the modulus of which tends to infinity as $t\to0$ (i.e.,
as $y\to k$)
according to Example \rf{hilfe}.
\et
\bex\lb{munonsmooth}
We demonstrate that $\bar\mu^J$ in general is not
$\tau_0$-$\tau_0$-continuous. To this end, it is sufficient to show
that $\Phi_\mu:\vphi\mapsto(\vphi\circ\mu^{-1})\cdot|(\mu^{-1})'|=
(\vphi\circ\mu^{-1})\cdot\frac{1}{|\mu'\circ\mu^{-1}|}$ is not
a continuous map from $\ca_0(\ti\Om)$ into $\ca_0(\Om)$ with
respect to the topology $\tau$ induced by the (LF)-topology
of $\cd(\R)$, for some
open subsets $\ti\Om,\Om$ of $\R$ and $\mu:\ti\Om\to\Om$
a suitable diffeomorphism. Consider $\ti\Om:=\Om:=(0,\infty)$;
choose $\rho\in\cd(\R)$ as in Example \rf{hilfe}, that is,
$\supp\rho\subseteq[0,2]$, $\int\rho=0$ and
$\rho(x)=\exp(-\frac{1}{x})$ for $0<x\le1$. Further, fix
any $\psi\in\ca_0(\Om)$. Then $\vphi_n(\xi):=
\frac{1}{n}\rho(\xi-\frac{1}{n})+\psi(\xi)$
defines a sequence converging to $\psi$ in $\ca_0(\Om)$ with
respect to $\tau$. Now consider $\mu:\ti\Om\to\Om$ defined by
$\mu(\xi):=\frac{\xi}{3}\exp(-\frac{3}{\xi})$. Then
the sequence formed by
$\Phi_\mu(\vphi_n)$ is not even bounded with respect to $\tau$:
Evaluating $\Phi_\mu(\frac{1}{n}\rho(\xi-\frac{1}{n}))$
at $\mu(\frac{2}{n})$ yields
\beas
&& \frac{1}{n}\cdot\rho\Big(\frac{2}{n}-\frac{1}{n}\Big)\cdot
\frac{1}{\mu'(\frac{2}{n})} =
\frac{1}{n}\cdot e^{\frac{n}{2}}\cdot\frac{6}{2+3n} \to \infty \quad (n\to
\infty) .
\eeas
\et
To  conclude  this  section,  we  briefly  introduce  the  notion  of
association into the diffeomorphism-invariant setting:
\bd $R_1,\ R_2 \in \mathcal{G}(\Om)$ are called associated ($R_1\approx R_2$) if the 
following condition holds:  $\forall \psi\in \D(\Om) \exists q\in\N\ 
\forall \phi\in\cc^\infty_b( I\times\Om,\ca_q(\R^s))$:
\beas
&& \lim_{\eps\to 0}\int (R_1-R_2)(S_\eps\phi(\eps,x),x)\psi(x) \,dx = 0
\eeas
\et
(where we have used the C-formalism). The  concept  of associated distribution 
as well as the basic properties of association are analogous to the 
non-diffeomorphism invariant case.

\section{Sheaf properties}\lb{sheaves}
Localization  properties  of  elements  of  ${\mathcal  G}$  are  most
conveniently  formulated  in  terms  of  sheaves. This section presents the
relevant notations and results.

Let $\Om \subseteq \Om'$. Then $U(\Om') \subseteq U(\Om)$ and for $R\in \G(\Om)$
we denote by $R|_{\Om'}$ the restriction of $R$ to $U(\Om')$. 
\bt
$\Om \to \mathcal{G}(\Om)$ is a fine sheaf of differential algebras.
\et
\pr
Let $\Om = \bigcup\limits_{\al\in A} \Om_\al$. We have to show that
\begin{itemize}
\item[(S1)] If $R_1$, $R_2 \in \G(\Om)$ and $R_1|_{\Om_\al} = R_2|_{\Om_\al}$ 
for all $\al\in A$ then $R_1 = R_2$.
\item[(S2)] If for each $\al\in A$ we are given $R_\al \in \G(\Om_\al)$ such
that $R_\al|_{\Om_ \al\cap \Om_\beta} = R_\beta|_{\Om_ \al\cap \Om_\beta}$ for
all $\al, \beta$ with $\Om_\al \cap \Om_\beta \not= \emptyset$ then there exists
some $R\in \G(\Om)$ with $R|_{\Om_\al} = R_\al$ for all $\al \in A$.
\item[(F)] If $(\Om_\al)_\al$ is locally finite there exists a family of sheaf
morphisms $\eta_\al: \G \to \G$ such that 
\begin{itemize}
\item[(i)]  $\sum_{\al\in A} \eta_\al = \mathrm{id}$.
\item[(ii)] $\eta_\al(\G_x) = 0$ for all $x$ in a neighborhood of $\Om\setminus
\Om_\al$ (where $\G_x$ denotes the stalk of $\G$ at $x$).
\end{itemize}
\end{itemize}
Noting that any $K\comp \Om$ can be written as $K=\bigcup_{\al\in A}K_\al$,
$K_\al \comp \Om_\al$, $K_\al = \emptyset$ $\forall \al\in A\setminus H$,
$|H| < \infty$, (S1) follows directly from
Definition \rf{defnegl}. For proving (S2) we adapt 
a construction from
\cite{JEL}, 21, Theorem 1. Choose a locally finite covering $(W_j)_{j\in
\N}$ of $\Om$ such that for each $j\in \N$ there exists $\al(j)$ with
$\overline{W}_j\comp \Om_{\al(j)}$. Let $(\chi_j)_{j\in\N}$ be a partition of
unity subordinate to $(W_j)_{j\in \N}$. Moreover, for each $j\in \N$ let
$\theta_j \in \D(\Om_{\al(j)})$, $\theta_j \equiv 1$ in a neighborhood of 
$\overline{W}_j$ and let $\psi_j\in \A_0(\Om_{\al(j)})$.
The map
\[
\pi_j(\vphi,x) := (\theta_j(\,.\,+x)\vphi - (\int\theta_j(\xi)\vphi(\xi-x)d\xi - 1) 
\psi_j(\,.\,+x),x)
\]
is smooth from $U(\Om)$ to $T^{-1}(\A_0(\Om_{\al(j)}) \times \Om))$ and
$\pi_j|_{U(W_j)} = \mathrm{id}$. 
Then for each $j\in \N$ $R_j:U(\Om)\to \C$,
\[
R_j(\vphi,x) = \left\{ \begin{array}{cl} 
                        \chi_j(x)R_{\al(j)}(\pi_j(\vphi,x)) & x \in  \Om_{\al(j)}\\
            0 & x \not \in \Om_{\al(j)}
            \end{array} \right.
\]
is smooth and $R_j|_{W_j} = R_{\al(j)}|_{W_j}$. Since $(W_j)_{j\in \N}$ is locally finite
\[
R(\vphi,x) := \sum_{j\in \N} R_j(\vphi,x)
\]
is an element of $\E(\Om)$. To show that $R$ is moderate we first note that in a neighborhood
of any $K\comp \Om$ only finitely many $R_j$ do not vanish identically, so it is enough to 
estimate one single $R_j$. Let $\phi \in \cc_b^\infty( I\times\Om,\A_0(\R^s))$ and choose
$L \comp W_j$ with $\supp(\chi_j)\comp L$. There exists $\eps_0>0$ such that for all $\eps\le\eps_0$ 
and all $x$ in a compact neighborhood of $L$ in $W_j$ $\supp(S_\eps\phi(\eps,x)) \subseteq W_j -x$,
so $(S_\eps\phi(\eps,x),x) \in U(W_j)$. On this set, $\pi_j = \mathrm{id}$ from which the claim 
follows by our assumption on $R_j$.

To establish (S2), by (S1) it suffices to show that
$R|_{\Om_\al\cap W_k} = R_{\al(k)}|_{\Om_\al\cap W_k}$  for all 
$k \in \N$ and all $\al \in A$ (Note that  $R_{\al(k)}|_{\Om_\al\cap W_k} = R_\al|_{\Om_\al\cap W_k}$
for any $\al$ by the assumption in (S2)). Now
\begin{equation} \label{se}
R|_{\Om_\al\cap W_k} - R_{\al(k)}|_{\Om_\al\cap W_k} = \sum_{j\not=k} \chi_j(R_{\al(j)}\circ \pi_j
- R_{\al(k)})|_{\Om_\al\cap W_k}
\end{equation}
For $K\comp \Om_\al \cap W_k$ and $j\not=k$ set 
$L=K_1 \cap \supp(\chi_j)$. Let $\phi \in \cc_b^\infty( I\times(\Om_\al\cap W_k),
\A_0(\R^s))$. Then for $x$ in a neighborhood $M\comp \Om_\al\cap W_j\cap W_k$ of $L$
 and sufficiently small $\eps$, 
$(S_\eps\phi(\eps,x),x) \in U(\Om_\al\cap W_j \cap W_k)$, so $\pi_j(S_\eps\phi(\eps,x),x) =
(S_\eps\phi(\eps,x),x)$. Hence the $\mathcal{N}$-estimates for the $j$-th term in (\ref{se})
follow from those of $R_{\al(j)}|_{\Om_{\al(j)}\cap \Om_{\al(k)}} - R_{\al(k)}|_{\Om_{\al(j)}\cap 
\Om_{\al(k)}}$ on $M$ and the fact that $\chi_j$ vanishes identically in a neighborhood of
$K\setminus M^\circ$. Finally, for proving (F) set (for $\Om' \subseteq \Om$)
\[
\eta_\beta|_{\Om'} := \G(\Om') \ni R \to 
\sum_{\{j | \beta = \al(j)\}} \!\!\!\chi_j\, (R|_{\Om'\cap W_j}\circ\pi_j|_{\Om'})
\]
\ep

\section{Separating the basic definition from testing}\lb{deftest}
Having introduced a  diffeomorphism  invariant  Colombeau  algebra 
in section \ref{jelshort}, we briefly return to the
general discussion of full Colombeau algebras.
Regarding the definitions of moderateness resp.\ negligibility
({\bf (D3)},{\bf (D4)}), we have adopted the terminology of
``testing'' in section \rf{scheme1}: By Definition {\bf (D1)},
certain ``objects'' (i.e., functions) $R$ are specified; those
which are singled out by Definition {\bf (D3)}
as being moderate serve as representatives
of elements of the algebra $\cg=\ce_M\left/\cn\right.$ of generalized functions.
The process of deciding whether an object $R$ belongs to $\ce_M$ (or to $\cn$,
respectively) has been called ``testing for moderateness resp.\
negligibility''.
It is performed by scaling ``test objects'' of the appropriate type
by the operator $S_\eps$ (as well as translating them
appropriately, whenever the J-formalism is used), plugging them into $R$
and analyzing the resulting behaviour of $R$ on these ``paths''
as $\eps\to0$.
Depending on the type of Colombeau algebra that is to be
constructed, test objects
take different forms, for example, $\vphi$ (\cite{c2}), $\phi(\eps)$
(\cite{CM}), $\phi(x)$ (\cite{JEL}) or $\phi(\eps,x)$ (\cite{JEL}, \cite{vw2}).
As opposed to that, the objects $R$ themselves do {\it not} depend in any way on
$\eps$; neither do they depend on $x$ via the first argument (the ``$\vphi$-slot'').
In other words, $R$ accepts only certain pairs $(\vphi,x)$ as arguments where
$\vphi\in\ca_0(\R^s)$. $x\in\R^s$. Summarizing, we adopt the following policy:

{\bf Defining the objects $R\in\ce(\Om)$
is to be separated strictly from testing them.}

This decision is based on the following reasons:
                 
First, it makes the objects simpler and the theory easier to comprehend,
yet without restricting its potential. Second, it provides a unifying
framework and a common terminology by means of which the different versions of
Colombeau algebras and their relations to each other can be analyzed. 
Finally, it is crucial for the development of algebras of nonlinear generalized
functions on smooth manifolds, if this is to be achieved in terms of
intrinsic objects;  this task is deferred to a subsequent paper
(\cite{vi}, jointly with J. Vickers).

Supposing $R$ to be the image of a non-smooth distribution $u$ unter the
corresponding embedding $\io$ into $\ce_M$, $R(S_\eps\phi,x)$ can be thought
of as a regularization of $u$: Indeed, $S_\eps\phi$ tends
to the delta distribution weakly, due to $\int\phi\equiv1$. In this sense, $S_\eps\phi$
(e.g., $S_\eps\phi(\eps,x)$) represents a ``smoothing process'' in its totality,
for all $\eps\in I$ and on the whole $x$-domain $\Om$. Separating the
definition of the objects $R$ from testing them thus amounts to assuming
that $R$ does not respond to the smoothing process as a whole but only
to its particular stages (represented by single elements
$\vphi$ of $\ca_0(\R^s)$). 

In the literature, three variants of increasing complexity can be
distinguished:

\begin{enumerate}
\item The objects $R$ take (certain) pairs $(\vphi,x)$ as arguments; testing
   is performed by inserting $(S_\eps\phi(\eps,x),x)$ into $R$;
   the behaviour of $R(S_\eps\phi(\eps,x),x)$ has to be studied.
\item Each object $R$ is given by a family $(R_\eps)_{\eps\in I}$
   of functions $R_\eps$ as in 1.; testing is performed by investigating
   $(R_\eps(S_\eps\phi(\eps,x),x))_\eps$.
\item The objects $R$ are defined on some set of pairs $(\cs,x)$ where
   $\cs=((\eps,x)\mapsto S_\eps\phi(\eps,x))$ resp.\
   $\cs=((\eps,x)\mapsto \phi(\eps,x))$ represents some ``smoothing process''.
   For testing $R$, $R(\cs,x)$ (which, in turn, has to be dependent on $\eps$!)
   has to be studied as $\eps\to0$.
\end{enumerate}
The first of the above variants is the one corresponding to ``separation
of definitions from testing''. Any object from level $i$ gives rise to
an object of level $(i+1)$ ($i=1,2$) by the following assignments:
\vspace{2.5mm}
$$\begin{array}{lrcll} 
\mbox{\rm level 1 $\to$ level 2}\qquad&R_\eps           &:=&R\qquad\qquad\qquad&
                                              \mbox{\rm(for all $\eps$)}\\
\mbox{\rm level 2 $\to$ level 3} &(R(\cs,x))(\eps,x)&:=&
                           R_\eps(S_\eps\phi(\eps,x),x).&
\end{array}$$

\vskip1mm
Jel\'\i nek in \cite{JEL}, Definition 5,
definitely chose level 1 for performing his
construction: This is made explicit in the last paragraph of item 2 of
\cite{JEL} (see also the discussion in item 3 of \cite{JEL}).
As opposed to that, Definition 5 of \cite{vw2}, e.g., clearly aims
at level 3 (the following definition of moderateness
(Definition 6), however, is ambiguous since it is not clear in which way
$R(\cs,x)$ (using our notation) depends on $\eps$).
The authors of \cite{CM}, on the other hand,
introduced their basic objects ${\cal R}\in\ce(\Om)$
in Definition 2
as smooth maps ${\cal R}:\ca_0\times\Om\to\C^{ I}$
where $\ca_0$ denotes a certain set of bounded paths
$\phi(\eps)$ (i.e., smoothing processes that are independent of
$x\in\Om$). At first glance,
this seems to be a clear indication that it was level 3 they
had in mind. In the following line, however,
${\cal R}(\phi,x)_\eps$ is specified to be of the form
$R(S_\eps\phi(\eps),x)$ (using our notation $S_\eps$ for the
scaling operator)
which has the appearance of level 2, generated by an object $R$ of
level 1 in the way described above.
As the case may be, using $\C^{ I}$
as range space for ${\cal R}$ instead of $\C$ 
definitely incorporates a certain part of the testing
procedure into the definition of the basic objects by introducing
$\eps$ as parameter from the very beginning.

\section{Characterization results}\lb{atoz}
The aim of this section is to derive several characterizations of moderateness
and negligibility, respectively, which
turned out to be indispensable tools in establishing the
diffeomorphism invariance of the algebra constructed in section
\rf{jelshort}.
Moreover, these characterizations
will serve as a basis for an intrinsic formulation of the theory on
manifolds (\cite{vi}).

We begin by proving Theorem \rf{JT1834infty}. To
this end, we introduce ``descending'' sequences of
linear projections $P_0,\dots,P_m$ with the property that
$P_0$ acts as the identity operator on $\ca_0(\R^s)$, $P_m$ projects
$\ca_0(\R^s)$ onto $\ca_q(\R^s)$ and the range of $P_j$ is of
codimension 1 in the range of $P_{j-1}$ ($j=1,\dots,m$).

Fix $q\in\N$ and $r>0$. Enumerate
$\{\bet\mid1\le|\bet|\le q\}$ in an arbitrary manner as
$\{\bet_1,\dots,\bet_m\}$. Since the family
$\{\xi^\bet\mid0\le|\bet|\le q\}$ is linearly independent in
$\cd'(B_r(0))$, there exist $\vphi_1,\dots,\vphi_m\in
\ca_{00}(\R^s)\cap\cd(B_r(0))$ satisfying 
$\int\xi^{\bet_i}\vphi_j(\xi)\,d\xi=\de_{ij}$
($1\le i,j\le m$). Now set
$$P_j:=\id{}_{\ca_0(\R^s)}-\sum_{i=1}^{j}\vphi_i\otimes\xi^{\bet_i},$$
that is,
$$P_j(\vphi):=\vphi-\sum_{i=1}^{j}\Big(\int\xi^{\bet_i}\vphi(\xi)\,d\xi\
\Big)
\cdot\vphi_i$$
for $\vphi\in\ca_0(\R^s)$, $j=0,\dots,m$. Obviously, the
operators $P_j$ satisfy the properties mentioned above.
As to using projections as defined above and the mean value theorem
in the subsequent proof we follow \cite{JEL}.

{\bf Proof of Theorem \rf{JT1834infty}. }
It is clear that condition $(4^\infty)$ implies the condition
specified in Definition \rf{defnegl}. So let us prove the converse,
assuming, in addition, that $R$ is moderate.
Let $K\subset\subset\Om$, $\al\in\N_0^s$ and $n\in\N$ be given.
In order to derive the required estimates for
$\pa^\al(R(S_\eps\phi(\eps,x),x))$ ($\phi$ being a test object of
type $[\mathrm{A}_\mathrm{l}^\infty]_{K,q}$) we have to provide
sets of the form $\ca_{0,M_1}(\R^s)\times U\subseteq
U_{\eps,M_2}(\Om)$ ($M_1\subset\subset\R^s$,
$M_2\subset\subset\Om$, $U$ an open subset of $\Om$), the former
of which will serve as domain for the operations of calculus to be
performed, according to section \rf{calcuepsom}.
Choose $L$ with $K\subset\subset\ L\subset\subset\Om$ and $r>0$
such that $\supp \phi(\eps,x) \subseteq B_r(0)$ 
for $\eps\in I$, $x\in L$. Then 
$\supp \pa^\bet\phi(\eps,x)\subseteq B_r(0)$ 
for $\eps\in I$, $x\in K$,
$\bet\in\N_0^s$. Let $M_1:=\clb r(0)$ and pick $M_2$
satisfying $L\subset\subset M_2\subset\subset\Om$.
According to Proposition \rf{slU} there exists
$\eps_0>0$ such that for $\eps\le\eps_0$, $U_{\eps,M_2}(\Om)$ contains
$\ca_{0,M_1}(\R^s)\times L$. Hence $\ca_{0,M_1}(\R^s)\times L^\circ$
is an appropriate domain for our purpose, 
supposing $\eps\le \eps_0$ in the sequel.

To obtain a suitable value of $q$, corresponding to $K,\al,n$
fixed above, we proceed as follows:

\ben
\item  By Theorem \rf{JT17} (note that $R$ was assumed to be
       moderate), there exists $N\in\N$ such that
       for every $k=0,\dots,|\al|$, for every $\bet\in\N_0^s$
       with $0\le|\bet|\le|\al|+1$ and for every bounded subset
       $B$ of $\cd(\R^s)$,
       $$\pa^\bet\rmd_1^k(R\circ S^{(\eps)})(\vphi,x)(\psi_1,\dots,\psi_k)
       =O(\eps^{-N})\qquad\qquad (\eps\to0)$$
       uniformly for $x\in K$, $\vphi\in B\cap\ca_0(\R^s)$,
       $\psi_1,\dots,\psi_k\in B\cap\ca_{00}(\R^s)$.
\item  By Definition \rf{defnegl}
       there exists $q\in\N$ such that for every
       $\phi\in\cc^\infty_b( I\times\Om,\ca_q(\R^s))$,
       $$\pa^\al(R(S_\eps\phi(\eps,x),x))=O(\eps^n)\qquad\qquad (\eps\to0)$$
       uniformly for $x\in K$.
\item  Without loss of generality, we may suppose $q\ge n+N$.
\een
Now let $\phi\in\cc^\infty_b( I\times\Om,\ca_0(\R^s))$ be of
type $[\mathrm{A}_\mathrm{l}^\infty]_{K,q}$.
For the values of $q$ and $r$ determined so far, consider $m$,
$\vphi_1\dots,\vphi_m$ and the projections $P_0,\dots,P_m$
introduced above. Defining $c_i\in\cc^\infty_b( I\times\Om,\C)$
by $c_i(\eps,x):=\int\xi^{\bet_i}\phi(\eps,x)(\xi)\,d\xi$
($i=1,\dots,m$), there exists a constant $C\ge1$ such that
$|\pa^\ga c_i(\eps,x)|\le C\eps^q(\le C)$ for all $x\in K$,
$0\le|\ga|\le|\al|$, $i=1,\dots,m$, due to $\phi$ being of type
$[\mathrm{A}_\mathrm{l}^\infty]_{K,q}$. In order to benefit from
the moderateness of $R$ in form of the differential condition
above, there remains an appropriate bounded subset of $\cd(\R^s)$
to be specified. To this end,
let
$$B:=\Gamma\, \{\pa^\ga\phi(\eps,x)\mid
   x\in K,\ 0\le|\ga|\le|\al|,\ \eps\in I\}+
   mC\cdot\Gamma\, \{\vphi_i\mid i=1\dots,m\}$$
where $\Gamma A$ denotes the absolutely convex hull of the subset
$A$ of an arbitrary linear space.
Since $B\cap\ca_0(\R^s)\subseteq\ca_{0,M_1}(\R^s)$ and
$B\cap\ca_{00}(\R^s)\subseteq\ca_{00,M_1}(\R^s)$, we may safely
use differentials of $R\circ S^{(\eps)}$ for $\eps\le\eps_0$,
evaluated at the respective vectors and $x\in K$.

For $j=0,\dots,m$, set                                 
$\psi_j(\eps,x):=P_j\phi(\eps,x):=\phi(\eps,x)-\sum_{i=1}^{j}
   c_i(\eps,x)\vphi_i.$
Then, in particular, $\psi_0=\phi$ and $\psi_m\in
\cc^\infty_b( I\times\Om,\ca_q(\R^s))$.
Restricting our attention to $x\in K$, $\eps\le\eps_0$,
the following statements are easily verified:
$$\pa^\ga\psi_m(\eps,x)\in\ca_{q,M_1}(\R^s)
     \qquad\qquad(0\le|\ga|\le|\al|),$$
$$\pa^\ga\psi_j(\eps,x)\in B,\ \vphi_j\in B
\qquad\qquad (j=0,\dots,m,\ 0\le|\ga|\le|\al|),$$ and, for every
$t\in[0,1]$, $j=1,\dots,m$,
$$\psi_j(\eps,x)+t\cdot c_j(\eps,x)\vphi_j=
  \phi(\eps,x)-\sum_{i=1}^{j-1}c_i(\eps,x)\vphi_i
              -(1-t)\cdot c_j(\eps,x)\vphi_j\in B.$$

By our choice of $q$,
$\sup_{x\in K}|\pa^\al(R(S_\eps\psi_m(\eps,x),x))|=O(\eps^n)$ and
we are going to show in the sequel that also
$$\sup_{x\in K}|\pa^\al(R(S_\eps\phi(\eps,x),x))-
   \pa^\al(R(S_\eps\psi_m(\eps,x),x))|=O(\eps^n)$$
which will complete the proof of $R$ satisfying $(4^\infty)$.
For the sake of simplicity, we will omit the argument
$(\eps,x)$ for the functions $\phi,\psi_j,c_j$ in the following.
According to the chain rule (the proof for the case $\al=0$
is contained trivially in the argument to follow), we obtain
\beas
   &&\kern-10pt
   \pa^\al(R(S_\eps\phi,x))-
   \pa^\al(R(S_\eps\psi_m,x))=
   \sum_{j=1}^{m}
   \big[\pa^\al(R(S_\eps\psi_{j-1},x))-
   \pa^\al(R(S_\eps\psi_j,x))\big]=\\
   &&\kern-10pt
   \sum_{j=1}^{m}
   \sum_{\bet,p}\big[(\pa^\bet\rmd_1^p(R\circ S^{(\eps)}))
       (\psi_{j-1},x)
       (\pa^{\ga_1}\psi_{j-1},\dots,\pa^{\ga_p}\psi_{j-1})-\\
   &&\kern-10pt\hphantom{mmmmmmmmmmmmmmmmm}
       (\pa^\bet\rmd_1^p(R\circ S^{(\eps)}))
       (\psi_j,x)
       (\pa^{\ga_1}\psi_j,\dots,\pa^{\ga_p}\psi_j)\big]
\eeas
where the second sum extends over certain
$\bet,p;\ga_1,\dots,\ga_p$. Obviously it is sufficient to derive the
desired estimate for each of the terms in square brackets
separately; thus fix $\bet,p;\ga_1,\dots,\ga_p$. Substituting
$\pa^{\ga_i}\psi_{j-1}=\pa^{\ga_i}\psi_j+\pa^{\ga_i}c_j\vphi_j$
($i=1,\dots,p$) and using multilinearity and symmetry of the iterated
differential transforms the square-bracket term into the sum of
\bea\lb{typ1}
     &&\kern-10pt
       (\pa^\bet\rmd_1^p(R\circ S^{(\eps)}))
       (\psi_{j-1},x)
       (\pa^{\ga_1}\psi_j,\dots,\pa^{\ga_p}\psi_j)-\nn\\
     &&\kern-10pt\hphantom{mmmmmmmmmm}
       (\pa^\bet\rmd_1^p(R\circ S^{(\eps)}))
       (\psi_j,x)
       (\pa^{\ga_1}\psi_j,\dots,\pa^{\ga_p}\psi_j)
\eea     
and of $2^p-1$ terms of the form
\bea\lb{typ2}
     &&\kern-10pt
       (\pa^\bet\rmd_1^p(R\circ S^{(\eps)}))
       (\psi_{j-1},x)
       (\pa^{\ga_{i_1}}c_j\vphi_j,\dots,\pa^{\ga_{i_l}}c_j\vphi_j,
        \pa^{\ga_{i_{l+1}}}\psi_j,\dots,\pa^{\ga_{i_p}}\psi_j)
\eea
where $\{i_1\dots,i_p\}=\{1,\dots,p\}$ and $1\le l\le p$.
Let us consider (\rf{typ2}) first. Observing that
$\sup\limits_{x\in K}|\pa^{\ga_{i_1}}c_j|\le C\eps^q$ and that
$\psi_{j-1},\vphi_j,\pa^{\ga_{i_2}}c_j\vphi_j,\dots,
\pa^{\ga_{i_l}}c_j\vphi_j,
\pa^{\ga_{i_{l+1}}}\psi_j,\dots,\pa^{\ga_{i_p}}\psi_j$ all are
members of $B$ (provided $\eps\le\eps_0$ and $x\in K$)
we can make use of the moderateness of $R$ in form of
the property derived previously by means of Theorem \rf{JT17}
to conclude that
\beas
    &&\kern-13pt\sup_{x\in K}
      |(\pa^\bet\rmd_1^p(R\circ S^{(\eps)}))
      (\psi_{j-1},x)
      (\pa^{\ga_{i_1}}c_j\vphi_j,\dots,\pa^{\ga_{i_l}}c_j\vphi_j,
       \pa^{\ga_{i_{l+1}}}\psi_j,\dots,\pa^{\ga_{i_p}}\psi_j)|=\\
      &&\kern-10pt\hphantom{mmmmmmmmmmmmmmmmmmmmm}
      O(\eps^q)\cdot O(\eps^{-N})=O(\eps^{q-N})\le O(\eps^n),
\eeas
due to our choice of $q$. To handle (\rf{typ1}), on the other hand,
we apply the mean value theorem \rf{mvth} to obtain
that the value of the term (\rf{typ1}) is contained in the closed
convex hull of the set of all (complex) numbers
\bea\lb{typ1mvth}
       (\pa^\bet\rmd_1^{p+1}(R\circ S^{(\eps)}))
       (\psi_j+tc_j\vphi_j,x)
       (\pa^{\ga_1}\psi_j,\dots,\pa^{\ga_p}\psi_j,c_j\vphi_j)
\eea
where $0<t<1$, $x\in K$. Taking into account that
$\sup\limits_{x\in K}|c_j(\eps,x)|\le C\eps^q$ and that
each of $\psi_j+tc_j\vphi_j,
\pa^{\ga_1}\psi_j,\dots,\pa^{\ga_p}\psi_j,\vphi_j$ is a member of
$B$ (provided $\eps\le\eps_0$ and $x\in K$)
the modulus of each of the complex numbers given by (\rf{typ1mvth})
can be estimated by $C'\eps^q\eps^{-N}$ for some positive constant
$C'$ being independent of $t,x$, again by the differential
condition derived from the moderateness of $R$. Consequently,
\beas
     &&\kern-10pt\sup_{x\in K}
       |(\pa^\bet\rmd_1^p(R\circ S^{(\eps)}))
       (\psi_{j-1},x)
       (\pa^{\ga_1}\psi_j,\dots,\pa^{\ga_p}\psi_j)-\\
     &&\kern-10pt\hphantom{mmmmmmmmmm}
       (\pa^\bet\rmd_1^p(R\circ S^{(\eps)}))
       (\psi_j,x)
       (\pa^{\ga_1}\psi_j,\dots,\pa^{\ga_p}\psi_j)|=O(\eps^n),
\eeas
thereby concluding the proof.
\ep

The remaining part of this section makes available the means
allowing to test moderateness resp.\ negligibility on test
objects $\phi$ which are defined only on certain subsets of
$ I\times\Om$. To this end we first have to provide
the technical toolkit for manipulating smooth bounded
paths.

\blem {\rm (Partition of unity on $ I$)}\lb{pu01}
Let $1>\eps_1>\eps_2>\eps_3>\ldots\to 0$, $\eps_0=2$.
Then there exist $\lambda_j\in\cd(\R)$ ($j=1,2,\ldots$) having the following
properties:

\beas
&& 1) \ {\rm supp}\, \lambda_j=[\eps_{j+1},\eps_{j-1}] 
\quad 2) \ \lambda_j(x)>0 \  \mbox{ for } \ x\in(\eps_{j+1},\eps_{j-1})  \\
&& 3) \ \sum_{j=1}^{\infty}\lambda_j(x)\equiv 1 \  \mbox{ for } \ x\in I \quad
4) \ \lambda_j(\eps_j)=1 \quad 5) \ \lambda_1(x)=1 \  \mbox{ for } \ x\in  [\eps_1,1]
\eeas
\et
\pr
For $j\in\N$, choose $\lambda_j^\circ\in\cd(\R)$ such that
$\mbox{supp}\,\lambda_j^\circ=[\eps_{j+1},\eps_{j-1}]$, $\lambda_j^\circ>0$ on
$(\eps_{j+1},\eps_{j-1})$ and $\lambda_0^\circ\in\cd(\R)$ such that
$\mbox{supp}\,\lambda_0^\circ=[1,3]$, $\lambda_0^\circ>0$ on
$(1,3)$. Define $\lambda^\circ:=
\sum\limits_{j=0}^{\infty}\lambda_j^\circ$ and
${\displaystyle\lambda_j(x):=\frac{\lambda_j^\circ(x)}{\lambda^\circ(x)}}$ for $x\in(0,3)$.
Then
$ \sum_{j=1}^{\infty}\lambda_j(x)=\sum_{j=0}^{\infty}\lambda_j(x)\equiv1$
for   $x\in I$
and it is easy to see that also the remaining four conditions are satisfied.
\ep

\blem\lb{inteps}
Assume that for each
$K\subset\subset\Om$ there exists $\eps\in I$ such that
the pair ($\eps,K$) satisfies a certain property (P) which is stable with respect
to decreasing $\eps$ and $K$ in the following sense:
\beas\quad\mbox{\rm If ($\eps_1,K_1$) satisfies (P) and $\eps_2\le\eps_1$, $K_2\subseteq
K_1$, then also ($\eps_2,K_2$) satisfies (P)}\hphantom{mm}&&\\ \mbox{\rm ($0<\eps_1,\eps_2\le1$,
$K_1,K_2\subset\subset\Om$).}&&
\eeas
Then there exists $\psi\in\cc^\infty(\Om)$ satisfying $0<\psi(x)\le1$ for each $x\in\Om$
such that an $\eps$ as above which is appropriate for $K$ with respect to (P) can be chosen as
$\eta_K:=\min\limits_{x\in K}\psi(x)$,
i.e., ($\eta_K,K$) satisfies (P) for each $K\subset\subset\Om$.
\et
\pr
Let $K_n$ be an increasing sequence of compact subsets of $\Om$ satisfying
$K_n\subset K_{n+1}^\circ$ which exhausts $\Om$ ($n\in\N$),
e.g.\ $K_n:=\{\lambda\in\Om\mid |\lambda|\le n\ {\rm and
\ dist}(\lambda,\pa\Om)\ge\frac{1}{n}\}$. Set $K_0:=\emptyset$ and define
open sets $G_n$ ($n\in\N$) by $G_1:=\emptyset$,
$G_n:=K_n^\circ\setminus K_{n-2}$ ($n\ge2$).
Now choose a partition of unity $(\varphi_n)_{n\ge1}$ subordinate to the
open covering $(G_n)_{n\ge1}$ of $\Om$.
For $n\ge2$, let $\eps_n$ be such that ($\eps_n,K_n$) satisfies
(P), $\eps_n\le\eps_{n-1}$ (we put $\eps_1:=1$) and define
$\psi(x):=\sum\limits_{n=1}^{\infty}\eps_n\varphi_n(x)$. Then it is clear that
$\psi$ is smooth and takes its values in $ I$. To complete the proof,
let $\emptyset\neq K\subset\subset\Om$. Let $N\in\N$  be minimal
such that $K$ is contained in $K_N^\circ$ and
consider $x\in K\setminus K_{N-1}^\circ$: Since $x\in K_N$, it cannot be an
element of $G_n$ for $n\ge N+2$. On the other hand, since $x\notin
K_{N-1}^\circ$, it neither can belong to $G_n$ for $n\le N-1$ (if $N\ge2$).
Altogether, $x\in G_N\cup G_{N+1}$ which results in
$\psi(x)=\eps_N\varphi_N(x)+\eps_{N+1}\varphi_{N+1}(x)\in[\eps_{N+1},\eps_N].$
Therefore,
$\eta_K:=\min_{x\in K}\psi(x)\le\min_{x\in K\setminus
K_{N-1}^\circ}\psi(x)\le\eps_N.$
So finally from $\eta_K\le\eps_N$ and $K\subseteq K_N$ we conclude 
that ($\eta_K,K)$ satisfies (P).
\ep
   
The proof of the preceding lemma also gives the following result
(just omit $K$ from the argument and consider
$x\in K_N\setminus K_{N-1}^\circ$ instead of $x\in K\setminus K_{N-1}^\circ$):

\blem\lb{sbe}
Let 
$K_n$ be an increasing sequence of compact subsets of $\Om$ satisfying
$K_n\subset K_{n+1}^\circ$ which exhausts $\Om$ ($n\in\N$)
and let $\eps_1\ge\eps_2\ge\ldots>0$ be given. Then there exists
$\psi\in\cc^\infty(\Om)$ satisfying $0<\psi(x)\le\eps_n$ for
$x\in K_n\setminus K_{n-1}^\circ$ ($K_0:=\emptyset$).
\et

\parskip=0mm
\bp {\rm (Extension of bounded paths)}\lb{ebp}

Let $\phi:D\to\ca_0(\R^s)$
where \linebreak $D\subseteq I\times\Om$ 
such that for each $K\subset\subset\Om$ there
exists $\eps_0>0$ and a subset $U$ of $D$ which is
open in $ I\times\Om$ having the following properties:
\begin{itemize}
\item[1)] $(0,\eps_0]\times K\subseteq U(\subseteq D)$
          and $\phi$ is smooth on $U$;
\item[2)] for all $\al\in\N_0^s$,
          $\{\pa^\al\phi(\eps,x)\mid0<\eps\le\eps_0,\ x\in K\}$ is bounded in $\cd(\R^s)$.
\end{itemize}
Then there exist a smooth map
$\ti\phi: I\times\Om\to\ca_0(\R^s)$
and $\si\in\cc^\infty(\Om)$ $(0<\si(x)\le1$ for all $x\in\Om)$
satisfying
\begin{itemize}
\item[$1'$)] $\ti\phi=\phi$ on
           $\{(\eps,x)\in  I\times\Om \mid \eps\le\frac{2}{3}\si(x)\}$;
\item[$2'$)] 
$\tilde \phi \in \cc^\infty_b( I\times\Om,\ca_0(\R^s))$.
\item[$3'$)] $(\ti\phi(\eps,x),x)\in U_\de(\Om)$ for all
           $(\eps,x)\in I\times\Om$ and $\de\le\si(x)$.
\end{itemize}
In particular, conditions $1'$) and $3'$) imply that
for each $K\subset\subset\Om$ there exists $\eps_1:=\min
\limits_{x\in K}\si(x)$ such that $\ti\phi=\phi$ on
an open neighborhood of $(0,\frac{\eps_1}{2}]\times K$
and $(\ti\phi(\eps,x),x)\in U_\de(\Om)$ for all
$(\eps,x)\in I\times K$ and $\de\le\eps_1$.
\et
\parskip=2mm
\pr
First we show that without loss of generality it can be assumed that
$\eps_0$ occurring in conditions 1) and 2) above also satisfies the following property
3), in addition to 1) and 2):
\begin{itemize}
\item[{\it 3)}]  \mbox{\it $(\phi(\eps,x),x)\in U_\de(\Om)$ for all $0<\eps,\de\leq\eps_0$ and
          $x\in K$.}
\end{itemize}       
In fact, according to Proposition \rf{slU} there exists $\eta>0$ such that
$(\phi(\eps,x),x)\in U_\de(\Om)$ for all $0<\eps\le\eps_0$, $x\in K$ and
$0<\de\leq\eta$ (observe that 
$\{\phi(\eps,x)\mid0<\eps\le\eps_0,\ x\in K\}$ is bounded by 2)).
Replacing both $\eps_0$ and $\eta$ by $\min(\eps_0,\eta)$, we see that 1)--3)
can be assumed to hold simultaneously.

Now let us say that ($\eps_0,K$) satisfies (P) if 1)--3) are valid for this
particular pair $(\eps_0,K)$. Then Lemma \rf{inteps} can be applied and provides
a function $\si\in\cc^\infty(\Om)$ satisfying $0<\si(x)\le1$ for each $x\in\Om$
such that 1)--3) hold for $\min\limits_{x\in K}\si(x)$ in place of $\eps_0$.
Now let $\lambda_1$ be smooth on $\R$, $0\le\lambda_1\le1$ and $\lambda_1\equiv1$
on $(-\infty,\frac{2}{3}]$, $\lambda_1\equiv0$ on $[\frac{5}{6},+\infty)$. Set
$\lambda_2:=1-\lambda_1$ and define
$$\ti\phi(\eps,x):=\lambda_1\left(\frac{\eps}{\si(x)}\right)\cdot\phi(\eps,x)+
                      \lambda_2\left(\frac{\eps}{\si(x)}\right)\cdot\phi(\si(x),x).$$
$\ti\phi$ is defined on $ I\times\Om$ since the formula above actually
involves only values of $\phi$ on pairs $(\eps,x)$ satisfying $\eps\le\si(x)$
and $(0,\si(x)]\times\{x\}\subseteq D$ by setting $K:=\{x\}$ in 1).
In order to show that $\ti\phi$ is smooth
we start by observing that
$\phi$ is smooth on some open neighborhood $U$($\subseteq D$)
of $D_\si:=\{(\eps,x)\mid x\in\Om,\ 0<\eps\leq\si(x)\}$:
Setting $K:=\{x\}$ in 1) once more yields an open neighborhood
$U_x$ of $(0,\si(x)]\times\{x\}$ on which $\phi$ is smooth. It
suffices to take
$U:=\bigcup\limits_{x\in\Om}U_x$. Now, since $x\mapsto(\si(x),x)$ is a
smooth map from $\Om$ into $U$, $\phi(\si(x),x)$ is smooth as a
function of $x$. Taking into account that
$\spp\lambda_1(\frac{\eps}{\si(x)})$
is a subset of
$\{(\eps,x)\mid\eps\le\frac{5}{6}\si(x)\}$
and noting that $\si(x)>0$ for all $x$
we see that also the first term in the definition of $\ti\phi$
and, hence, also $\ti\phi$ itself are smooth.

Obviously, $\ti\phi(\eps,x)=\phi(\eps,x)$ for $\eps\le\frac{2}{3}\si(x)$.
Thus $1'$) is proved.

To show $2'$), we have to consider derivatives of $\ti\phi$ with respect to $x$ on
sets of the form $ I\times K$ where $K\subset\subset\Om$ is given.
Again we set $\eps_1:=\min\limits_{x\in K}\si(x)$.
First, observe that on $(0,\eps_1]\times K$ all derivatives
$\pa^\bet\phi(\eps,x)$ are bounded by 2).
Since they are clearly also bounded on the compact set
$\{(\eps,x)\mid x\in K,\ \eps_1\le\eps\le\si(x)\}$,
they are bounded on the whole of
$K_\si:=\{(\eps,x)\mid x\in K,\ 0\le\eps\le\si(x)\}$. Now fix
$\al\in\N_0^s$; to discuss
$\pa^\al\ti\phi$ in detail, we set
$$\ti\phi_1(\eps,x):=\lambda_1\left(\frac{\eps}{\si(x)}\right)\cdot\phi(\eps,x)
\quad \mbox{resp.}\quad
\ti\phi_2(\eps,x):=\lambda_2\left(\frac{\eps}{\si(x)}\right)\cdot\phi(\si(x),x).$$
Now
$$\pa^\al\ti\phi_1(\eps,x)= \sum_{\bet+\ga=\al}{\al\choose\bet}\,
  \pa^\bet\lambda_1\left(\frac{\eps}{\si(x)}\right)\cdot\pa^\ga\phi(\eps,x).$$
As we have seen above, all derivatives $\pa^\ga\phi(\eps,x)$ are bounded on
$K_\si$. Expanding $\pa^\bet\lambda_1(\frac{\eps}{\si(x)})$ according to the
chain rule gives a finite number of terms of the form
$$\lambda_1^{(l)}\left(\frac{\eps}{\si(x)}\right)\cdot\eps^l \cdot
\pa^{\ga_1}\left(\frac{1}{\si(x)}\right)\cdot\ldots\cdot\pa^{\ga_l}\left(\frac{1}{\si(x)}\right)$$
where $1\le l\le|\bet|$ and
$\ga_1,\ldots,\ga_l\in\N_0^s$ satisfy $\sum\limits_{i=1}^{l}|\ga_i|=|\bet|$. Each of these
terms is bounded on $ I\times K$. Taking into account that all
$\pa^\bet\lambda_1(\frac{\eps}{\si(x)})$ vanish for
$\eps\ge\frac{5}{6}\si(x)$ and that $K_\si$ is characterized by
$\eps\le\si(x)$, $\pa^\al\ti\phi_1$ is bounded on $ I\times K$.

The derivatives of $\ti\phi_2$ take the form
$$ \pa^\al\ti\phi_2(\eps,x)=\sum_{\bet+\ga=\al}{\al\choose\bet}\,
  \pa^\bet\lambda_2\left(\frac{\eps}{\si(x)}\right)\cdot\pa^\ga\phi(\si(x),x).
$$
The above reasoning showing
$\pa^\bet\lambda_1(\frac{\eps}{\si(x)})$
to be bounded also applies to
$\pa^\bet\lambda_2(\frac{\eps}{\si(x)})$.
For $\pa^\ga\phi(\si(x),x)$ the chain
rule gives a finite sum of terms of the form
\be\lb{dedxphi}
\pa_\eps^k\pa_x^{\ga_0}\phi(\si(x),x)\cdot\pa_x^{\ga_1}\si(x)\cdot\ldots\cdot
  \pa_x^{\ga_k}\si(x)
\end{equation}
where $0\le k\le|\ga|$ and $\ga_0,\ga_1,\ldots,\ga_k\in\N_0^s$ satisfy
$\sum\limits_{i=0}^{k}|\ga_i|=|\ga|$.
Since
$\pa_\eps^k\pa_x^{\ga_0}\phi$ is bounded on the compact subset
$\{(\si(x),x)\mid x\in K\}$ of $U$,
all the factors in
(\rf{dedxphi}) are bounded on $K$.
Combining this with the boundedness of
$\pa^\bet\lambda_2(\frac{\eps}{\si(x)})$
shows that also $\pa^\al\ti\phi_2$ (and hence $\pa^\al\ti\phi$) is
bounded on $ I\times K$ which completes the proof of $2'$).

Finally, to show $3'$) let $x\in \Om$, $\de\le\si(x)$ and conclude from 3) that in
particular $(\phi(\si(x),x),x)\in U_\de(\Om)$. Now for $\eps\le\si(x)$, both
$(\phi(\eps,x),x)$ and $(\phi(\si(x),x),x)$ belong to $U_\de(\Om)$ and so
also $(\ti\phi(\eps,x),x)$ does. On the other hand, for  $\si(x)$ $<$ $\eps$
$\le$ $1$
we have $(\ti\phi(\eps,x),x)=(\phi(\si(x),x),x)\in U_\de(\Om)$.
Therefore, for all $\eps\in I$, $(\ti\phi(\eps,x),x)\in U_\de(\Om)$.

The last statement of the proof follows from the fact that
$\ti\phi$ and $\phi$ coincide on the open set 
$\{(\eps,x)\mid \eps<\frac{2}{3}\si(x)\}$ which clearly contains
$(0,\frac{\eps_1}{2}]\times K$. 
\ep

In the following, we will identify each function
$\hat\phi: I\to\cc^\infty(\Om,\ca_0(\R^s))$
in the natural way with the corresponding function
$\phi\in\cc^{[\infty,\Om]}( I\times\Om,\ca_0(\R^s))$
where $\hat\phi(\eps)(x)=\phi(\eps,x)$. This identification
respects the properties of $\hat\phi$ resp.\ $\phi$
being smooth (see Theorem \rf{exp}) and/or bounded (in the sense
specified in section~\rf{notterm}).

\bt\lb{a--z}$(\mathrm{A}$ to $\mathrm{Z})$
The moderateness of an element $R$ of
$\cc^\infty(U(\Om))$ can be tested equivalently
on boun\-ded subsets of $\cc^\infty(\Om,\ca_0(\R^s))$, on
arbitrary bounded paths
$\phi: I\to\cc^\infty(\Om,\ca_0(\R^s))$
or on paths of the same form depending smoothly on
$\eps$ (conditions {\rm(A)}, {\rm(B)}, {\rm(C)}, respectively).
\smallskip\\
Moreover, equivalent moderateness tests can be performed on larger resp.\
smaller classes of smooth paths which are distinguished by better resp.\ poorer
properties with respect to the domain of definition of $R(S_\eps\phi(\de,x),x)$
(conditions {\rm(D)}, {\rm(E)}; {\rm(Z)}).
Formally:\smallskip\\
Let $R\in\cc^\infty(U(\Om))$. Then the following conditions are equivalent:
\ms
{\rm(A)}
\quad$
\forall K\subset\subset\Om\ \forall\al\in\N_0^s\ \exists N\in\N\ 
\forall \cb\,(bounded)\subseteq\cc^\infty(\Om,\ca_0(\R^s))\\
{}\hskip15mm\exists C>0\ \exists \eta>0\ \forall\phi\in\cb\ \forall\eps\,(0<\eps<\eta)
    \ \forall x\in K$:
\beas
|\pa^\al(R(S_\eps\phi(x),x))|\leq C\eps^{-N}
\eeas
\vskip1mm
{\rm(B)}
\quad$
\forall K\subset\subset\Om\ \forall\al\in\N_0^s\ \exists N\in\N\ 
\forall \phi
\in\cc^{[\infty,\Om]}_b( I\times\Om,\ca_0(\R^s))\\
{}\hskip15mm\exists C>0\ \exists \eta>0\ \forall\eps\,(0<\eps<\eta)
    \ \forall x\in K$:
\beas
|\pa^\al(R(S_\eps\phi(\eps,x),x))|\leq C\eps^{-N}
\eeas
\vskip1mm
{\rm(C)}
\quad$
\forall K\subset\subset\Om\ \forall\al\in\N_0^s\ \exists N\in\N\ 
\forall \phi\,\in\cc^\infty_b( I\times\Om,\ca_0(\R^s))\\
{}\hskip15mm\exists C>0\ \exists \eta>0\ \forall\eps\,(0<\eps<\eta)
    \ \forall x\in K$:
\beas
|\pa^\al(R(S_\eps\phi(\eps,x),x))|\leq C\eps^{-N}
\eeas
\vskip1mm
{\rm(D)}
\quad\kern-3pt
as condition (C), yet only paths 
$\phi
\in\cc^\infty_b( I\times\Om,\ca_0(\R^s))$ are considered \\
{}\hphantom{{\rm(F),}}\quad such that $(\phi(\eps,x),x)\in U_\de(\Om)$
          for all $\eps,\de\in I$ and $x\in K$.
\vskip1mm
{\rm(E)}
\quad\kern+1pt
\hskip-3pt as condition (C), yet only paths 
$\phi
\in\cc^\infty_b( I\times\Om,\ca_0(\R^s))$ are considered \\
{}\hphantom{{\rm(F),}}\quad such that
                      for each $L\subset\subset\Om$ there exists $\de_0$ having the
                      the property\\ 
{}\hphantom{{\rm(F),}}\quad $(\phi(\eps,x),x)\in U_\de(\Om)$
                      for all $(\eps,x)\in I\times L$ and $\de\le\de_0$.
\vskip1mm
{\rm(Z)}
\quad$
\forall K\subset\subset\Om\ \forall\al\in\N_0^s\ \exists N\in\N
\ \forall\phi
:D\to\ca_0(\R^s))\ \mbox{($D,\phi$ as described
          below})\\
{}\hskip15mm\exists C>0\ \exists \eta>0\ \forall\eps\,(0<\eps<\eta)
    \ \forall x\in K$: $(\eps,x)\in D$ and
\beas
|\pa^\al(R(S_\eps\phi(\eps,x),x))|\leq C\eps^{-N}
\eeas
{}\hphantom{{\rm(F),}}\quad
   where $D\subseteq I\times\Om$ and for
   $D,\vphi$ the following holds: For each
   $L\subset\subset\Om$\\
{}\hphantom{{\rm(F),}}\quad there exists $\eps_0$ and a
          subset $U$
          of $D$ which is open in $ I\times\Om$ such that\\
{}\hphantom{{\rm(F),}}\quad
\vbox{\begin{itemize}
\item[1)] 
          $(0,\eps_0]\times L\subseteq U(\subseteq D)$
          and $\phi$ is smooth on $U$;
\item[2)] for all $\bet\in\N_0^s$,
          $\{\pa^\bet\phi(\eps,x)\mid0<\eps\le\eps_0,\ x\in L\}$ is bounded in $\cd(\R^s)$.
\end{itemize}}
\et

\pr
(A) $\Longrightarrow$ (B) is clear since the image of $\phi$,
the latter being viewed as a function from $ I$ into
$\cc^\infty(\Om,\ca_0(\R^s))$, forms a bounded
subset of $\cc^\infty(\Om,\ca_0(\R^s))$.\ms
(B) $\Longrightarrow$ (C) is trivial.
\ms
(C) $\Longrightarrow$ (A) Assume to the contrary  $\neg$(A), i.e.,
\ms
\hphantom{}\hskip10mm$\exists K\subset\subset\Om\ \exists\al\in\N_0^s\ \forall N\in\N\ 
\exists \cb\,\mbox{(bounded})\subseteq\cc^\infty(\Om,\ca_0(\R^s))\\
\hphantom{}\hskip50mm\forall C>0\ \forall \eta>0\ \exists\phi\in\cb\ \exists\eps\,(0<\eps<\eta)
    \ \exists x\in K$:
\beas
|\pa^\al(R(S_\eps\phi(x),x))|> C\eps^{-N}
\eeas
Fix $K,\al$ as given by $\neg$(A). (C) yields $N\in\N$ with the property described there
($\forall \phi\,(\cc^\infty\mbox{, bounded})$ \ldots). Now for this $N$,
our assumption $\neg$(A) produces a bounded subset $\cb$ of
$\cc^\infty(\Om,\ca_0(\R^s))$ on which $R$ behaves ``badly'' in the sense that we
can inductively define sequences $\eps_j\in I$, $\phi_j\in \cb$ and
$x_j\in K$ from which a path $\phi_0(\eps,x)$ can be constructed giving a
contradiction to (C) if $R$ is tested on it. Explicitly:

Set $C:=1$, $\eta:=1$ and conclude from $\neg$(A):
$$\exists \eps_1<1,\hphantom{\min(,\eps_2)}\ \phi_1\in\cb,\ x_1\in K:
\big|\pa^\al(R(S_{\eps_1}\phi_1(x),x))\left|_{x=x_1}\right.\big|>
 1\cdot\eps_1^{-N}.$$
Set $C:=2$, $\eta:=\min(\frac{1}{2},\eps_1)$ and conclude from $\neg$(A):
$$\exists \eps_2<\min(\frac{1}{2},\eps_1),\ \phi_2\in\cb,\ x_2\in K:
\big|\pa^\al(R(S_{\eps_2}\phi_2(x),x))\left|_{x=x_2}\right.\big|>
 2\cdot\eps_2^{-N}.$$
Continuing this way, we get $\eps_0:=1>\eps_1>\eps_2>\ldots\to0$,
$\phi_j\in \cb$, $x_j\in K$ satisfying
$$\big|\pa^\al(R(S_{\eps_j}\phi_j(x),x))\left|_{x=x_j}\right.\big|>
 j\cdot\eps_j^{-N}.$$
Take a partition of unity $(\lambda_j)_j$ on $ I$ as provided by Lemma
\rf{pu01} and define
$\phi_0(\eps,x)$ $:=$ ${\displaystyle\sum_{j=1}^{\infty}\lambda_j(\eps)\phi_j(x)}$
($\eps\in I$, $x\in\Om$). Clearly, $\phi_0$ is smooth from $ I$ into
$\cc^\infty(\Om,\ca_0(\R^s))$ and its image is bounded since it is contained in the convex
hull of $\cb$. By construction, we get
$$\big|\pa^\al(R(S_{\eps_j}\phi_0(\eps_j,x),x))\left|_{x=x_j}\right.\big|=
  \big|\pa^\al(R(S_{\eps_j}\phi_j(x),x))\left|_{x=x_j}\right.\big|> j\cdot\eps_j^{-N},$$
contradicting $\pa^\al(R(S_\eps\phi_0(\eps,x),x))=O(\eps^{-N})$ as required by (C).
\ms
(C) $\Longrightarrow$ (D) is trivial.
\ms
(D) $\Longrightarrow$ (E) Let $K,\al$ be given and choose $N$ by (D). Consider a
path $\phi(\eps,x)$ and a constant $\de_0>0$ appropriate for $K$, both as
specified in (E). Define
$$\ti\phi(\eps,x):=S_{\de_0}\phi(\de_0\eps,x)\qquad\qquad(\eps\in I,\ x\in\Om).$$
Since $S_{\de}\ti\phi(\eps,x)=S_{\de\de_0}\phi(\de_0\eps,x)$,
$(\ti\phi(\eps,x),x)\in U_\de(\Om)$ for all $\eps,\de\in I$ and $x\in
K$. (D) now gives
$|\pa^\al (R(S_\eps\ti\phi(\eps,x),x))|$ $=$ $O(\eps^{-N})$
uniformly on $K$. However, $\pa^\al (R(S_\eps\ti\phi(\eps,x),x))=
\pa^\al (R(S_{\de_0\eps}\phi(\de_0\eps,x),x))$ which implies
$$|\pa^\al
(R(S_\eps\phi(\eps,x),x))|=O\left(\left(\frac{\eps}{\de_0}\right)^{-N}\right)=
   O(\eps^{-N}).$$
(E) $\Longrightarrow$ (Z) Again let $K,\al$ be given; this time, of course,
choose $N$ according to (E). By
Proposition \rf{ebp}, to a given path $\phi(\eps,x)$ as in (Z) there exists a
bounded path $\ti\phi(\eps,x)$ satisfying the condition in (E) such that
$\ti\phi$ and $\phi$ coincide on an open neighborhood of
$(0,\frac{\eps_1}{2}]\times K$, where $\eps_1$ is some positive constant.
Therefore
$$|\pa^\al (R(S_\eps\phi(\eps,x),x))|=
  |\pa^\al (R(S_\eps\ti\phi(\eps,x),x))|=O(\eps^{-N}).$$
%
(Z) $\Longrightarrow$ (C) Setting $D:= I\times\Om$ turns a path in the sense
of (C) in a path as required for the application of (Z).
\ep

There is an analog to the preceding theorem giving rise to equivalent
definitions of negligibility. Observe, however, that each of the
conditions in the following theorem is equivalent to the condition
in Definition \rf{defnegl} even {\it without} assuming $R$ to be
moderate. This latter property has to be assumed in addition to
obtain the definition of negligibility.

\bt\lb{a--z/vm}$(\mathrm{A}'$ to $\mathrm{Z}')$
In properties $(\mathrm{A})$--$(\mathrm{Z})$ of Theorem \rf{a--z},
insert
``$\forall n\in\N$'' after ``$\forall\al\in\N_0^s$''
and replace ``$\exists N\in\N$'' by ``$\exists q\in\N$'',
``$\ca_0(\R^s)$'' by ``$\ca_q(\R^s)$'' and
``$C\eps^{-N}$'' by ``$C\eps^n$'', throughout. 
Then the resulting six conditions $(\mathrm{A}')$--$(\mathrm{Z}')$ are
mutually equivalent.
If $R$, in addition, is
supposed to be moderate, each of them is equivalent to the
negligibility of $R$.
\et

The proof of the preceding theorem is analogous to that of Theorem
\rf{a--z}. Finally,
in the proof of {\bf (T8)} in section \rf{jelshort}, we have made
use of the following variant of (a part of) \rf{a--z/vm}:
\bc\lb{a--z/avm}
In each of conditions $(\mathrm{C}')$ and $(\mathrm{Z}')$
of Theorem \rf{a--z/vm},
replace ``$\ca_q(\R^s)$'' by ``$\ca_0(\R^s)$'' and
consider only test objects $\phi$ all whose derivatives
$\pa^\al_x\phi$ ($\al\in\N_0^s$) have asymptotically vanishing
moments of order $q$ on the compact set $K$ at hand. Then the resulting
conditions $(\mathrm{C}'')$ and $(\mathrm{Z}'')$ are equivalent.
\et
\pr
$(\mathrm{Z}'')\!\Rightarrow\!(\mathrm{C}'')$ being trivial, we 
have to show the reverse implication. To this end, note that
for the property of a test object
to have asymptotically vanishing moments
on some $K\subset\subset\Om$, only sets of the form
$(0,\eps_0]\times K$
(for some $0<\eps_0\le1$) are relevant.
Yet it is part of
the statement of Proposition \rf{ebp} that the extended path
$\ti\phi$ agrees with the given path $\phi$ on sets of this form, for
every given $K\subset\subset\Om$. Thus
the property of having asymptotically vanishing moments
is preserved by the extension process
$\phi\mapsto\ti\phi$. 
Now an argument analogous to that used to prove
$(\mathrm{C})\!\Rightarrow\!(\mathrm{Z})$ in Theorem \rf{a--z}
establishes $(\mathrm{C}'')\!\Rightarrow\!(\mathrm{Z}'')$.
\ep

Conditions (A)--(C) in Theorem \rf{a--z} are due to J. Jel\'\i nek (\cite{JEL},
the Remark following Definition 8; the proof of equivalence is only indicated
there).
The   equivalence   of  condition  (Z)  with  (A)--(C) has to be
considered as the technical cornerstone of the diffeomorphism invariance of
the  Colombeau algebra constructed in section \rf{jelshort}.
Apart from that, it can be of advantage in certain situations (for
example, when dealing with applications)
not to have to bother too much about the domains of definition of
$R(S_\eps\phi(\de,x),x)$ being too small. On the other hand,
it  can  be useful to have guaranteed a certain minimum size of
these  domains; for this reason (D) and (E)
have been included in the theorem. Last, but not least we felt the need to
give precise meaning to
statements like ``only $\eps>0$ small enough is relevant''
(\cite{CM},
p.~361) or ``$\ldots$ in the case when the maps $\Phi_\eps$ [$\ldots$]
are not defined on the same set. We only need that [$\ldots$] for all
$\eps>0$ sufficiently small.'' (\cite{JEL}, item 7).
In our view, the considerable technical expense required for
establishing the results of this
section clearly shows the necessity of a rigorous treatment 
to be given to these matters.

\section{Differential Equations}\lb{diffeq}
The main application of Colombeau algebras so far has been in the field
of  differential  equations.  It  is  therefore of considerable interest to
explore  how  the  changes  in the construction of the algebra necessary to
ensure diffeomorphism invariance affect the process of solving differential
equations  in ${\mathcal G}$. To illustrate these changes, in the following
we  are  going  to  discuss  two  prototypical examples. 
\bex  \label{odesec}
Consider the initial value problem 
\bea\label{ode}
&& \ddot x (t) = f(x(t))\delta(t) \nonumber\\
&& x(-1) = x_0 \\
&& \dot x(-1) = \dot x_0   \nonumber
\eea
($f: \R \to \R$ smooth) in $\mathcal{G}(\R)$. Equations of type~(\ref{ode}) arrise, e.g.,
in geodesic equations in impulsive gravitational waves, cf. \cite{KS}. For simplicity, 
we assume that the initial values $x_0$ and $\dot x_0$ are real numbers. 

As in the case of non-diffeomorphism invariant Colombeau theory, solving (\ref{ode}) 
requires existence and uniqueness results
for the classical equation with $\delta(t)$ replaced by $\varphi(-t)$ (recall that 
a representative of $\iota(\delta)$ is given by $(\vphi,t) \to \vphi(-t)$). 
By a standard fixed point argument (cf. \cite{KS}, Lemma 1) this initial value problem has
unique global solutions provided $\supp(\vphi)$ is contained in a sufficiently small
neighborhood $U$ (depending on $f$ and the initial conditions) of $0$. For $\vphi \in
\D(U)$ let $R_1(\vphi,\,.\,)$ be this unique solution. Choose some $\chi \in \D(U)$, $\chi
\equiv 1$ in a neighborhood $V$ of $0$. We claim that $R: (\vphi,t) \to R_1(\chi\vphi,t)$ 
is a representative of a locally bounded solution to (\ref{ode}) and that this solution 
is unique. (Note that independence of $\cl[R]$ from $\chi$ follows already from the
fact that $S_\eps \vphi \chi = S_\eps \vphi$ for any $\vphi$ and $\eps$ sufficiently small).
First of all, in order to show that $R$ is smooth on $(U(\R),\tau_2)$ it obviously suffices 
to establish smoothness of $R_1$ on that space. To this end, let $s\to (\vphi_s,t_s)$ be 
a  $\tau_1$-smooth curve into some $\mathcal{A}_{0,H}(\R)\times V \subseteq
U_{N}(\R)$ as in Theorem \ref{smoothsub}. 
Then smoothness of solutions of ODEs with respect to a real parameter
at once shows that $s \to R_1(\vphi_s,t_s)$ is $C^\infty$, from which the claim follows
by definition of smoothness and Theorem \ref{smoothsub}. (We feel that the 
ease of this kind of 
argument is a decisive advantage of calculus in convenient vector spaces as used 
here compared to earlier approaches to differential calculus in locally convex 
spaces.)

In order to show
that $R\in \mathcal{E}_M(\R)$, note that for $\eps$ sufficiently small
\[
R(S_\eps\vphi,t) = x_0 + \dot x_0(t+1) +\int\limits_{-1}^t\int\limits_{-1}^s f(R(S_\eps \vphi,r))
S_\eps \vphi(-r) dr ds\,.
\] 
The key to proving the desired estimates is the characterization of moderateness given in 
Theorem  \ref{JT17}  on  the one hand and Remark \ref{epsfree} on
the    other:    For   $\vphi$,   $\psi$   varying   bounded   subsets   of
$\mathcal{A}_0(\R)$  (resp. $\mathcal{A}_{00}(\R)$) and $\eps$ sufficiently
small, iterated differentials of $R\circ S^{(\eps)}$ are well defined and
can  be  calculated  according  to  the  usual  differentiable structure of
$\mathcal{A}_0(\R)\times \R$. In particular,
differentiation with respect to $\vphi$ can be interchanged with integration
(see the proof of \ref{part}),
and the chain rule gives, e.g.:
\beas
&& d_1(R(S_\eps\vphi,t))[\psi] = 
\int\limits_{-1}^t\int\limits_{-1}^s f'(R(S_\eps\vphi,r))
d_1(R(S_\eps\vphi,r))[\psi] S_\eps\vphi(-r) drds \\
&& + \int\limits_{-1}^t\int\limits_{-1}^s f(R(S_\eps\vphi,r)) S_\eps\psi(-r) drds
\eeas
so the result follows (using Gronwall's inequality) by induction. Even in this rather simple
example it is quite obvious that to check the moderateness condition \ref{defmod} directly would
be extremely tedious (resp.\ unmanageable for more complicated equations). Moreover, uniqueness
can be established similarly without even having to perform any differentiations owing to
the remark following Theorem \ref{JT1823}.
\bex The semilinear wave equation 
\bea \label{wave}
(\pa_t^2 - \Delta) u &=& F(u) + H\\
u|_{\{t<0\}} &=& 0 \nonumber 
\eea
\et
with $F\in \mathcal{O}_M(\R)$ globally Lipschitz, $F(0)=0$ and $H\in \mathcal{G}$, $\supp H
\subseteq \{t\ge 0\}$ has been treated in the Colombeau 
framework in \cite{OR}. Therefore, we will only indicate those modifications that allow to carry 
over the existence and uniqueness results achieved there into the current setting. Also,
we only treat space dimension $3$. Let $R_H^\eta$ be a representative of $H$ supported in 
$\{t>-\eta\}$ ($\eta > 0$). For each $\vphi \in \mathcal{A}_0(\R^4)$ let $(x,t) \to 
R(\vphi,x,t)$ be the smooth solution to
\beas
  (\pa_t^2 - \Delta) u(x,t) &=& F(u(x,t)) + R_H^\eta(\vphi,x,t) \\
 u|_{\{t\le -\eta\}} &=& 0
\eeas
Then from Kirchhoff's formula we obtain for $t\ge -\eta$
\begin{equation} \label{waveint}
R(\vphi,x,t) = \frac{1}{4\pi} \int\limits_{-\eta}^t \frac{1}{t-s} \int\limits_{|x-y|=t-s}
(F(R(\vphi,y,s)) + R_H^\eta(\vphi,y,s)) d\sigma(y)ds
\end{equation}
Composing $R$ in this formula with a smooth curve as in \ref{odesec}, smooth dependence of solutions
of Volterra integral equations on real parameters implies smoothness of $R$ in $\vphi$ and as in the
classical  case,  smooth extension to $t<-\eta$ is possible since $F(0)=0$.
Again using Remark \ref{epsfree},
$\vphi$-differentiation  (for  $R\circ  S^{\eps}$,  $\eps$  small, $\vphi$,
$\psi_i$ (as in Theorem \ref{JT17}) varying in bounded sets) interchanges  
with  integration in (\ref{waveint}).
Hence derivation of the necessary $\mathcal{E}_M$- resp.\ $\mathcal{N}$-estimates for existence
resp. uniqueness of solutions is carried out analogously to the classical case, again due
to  Theorem \ref{JT17} and the remark following Theorem \ref{JT1823}.   
\et
{\bf Acknowledgements.} 
The work on this series of papers was initiated during a visit of the authors
at the department of mathematics of the university of Novi Sad in July 1998. 
We would like to thank the faculty and staff
there, in particular Stevan Pilipovi\'c and
his group for many helpful 
discussions and for their warm hospitality. Also, we are indebted to Andreas 
Kriegl for sharing his expertise on infinite dimensional calculus.

{\small
{\it Electronic Mail:} 

\begin{tabular}{ll}
E.F.: &  {\tt eva@geom2.mat.univie.ac.at}\\
M.G.: & {\tt michael@mat.univie.ac.at}\\
M.K.: & {\tt Michael.Kunzinger@univie.ac.at}\\
R.S.: & {\tt Roland.Steinbauer@univie.ac.at}
\end{tabular}
}


\begin{thebibliography}{99}

\bibitem{we} H.\ Balasin, 
        {\it Colombeau's generalized functions on arbitrary manifolds,} 
        gr-qc Preprint\footnote[1]{available
        electronically under {\tt http://xxx.lanl.gov/archive/gr-qc}}
        {\bf 9610017} (1996).
\bibitem{hb} H.\ Balasin, 
        {\it Distributional aspects of general relativity: 
        The example of the energy-momentum tensor of the
        extended Kerr-geometry,} 
        in M.\ Grosser, G.\ H\"ormann, M.\ Kunzinger, M.\ Oberguggenberger (Eds.), 
        Nonlinear Theory of Generalized Functions, Chapman \& Hall/CRC
        Res.\ Notes Math. {\bf 401}, Chapman \& Hall/CRC, Boca Raton, FL, 1999, 
        pp.\ 275--290.
\bibitem{b} H. A.\ Biagioni,
        {\it A Nonlinear Theory of Generalized Functions,} (2nd ed.) 
        Lecture Notes in Math. {\bf 1421}, Springer, New York, 1990. 
\bibitem{bmo} H. A.\ Biagioni, M.\ Oberguggenberger,
        {\it Generalized solutions to the Korteweg-de Vries and the regularized 
        long-wave equations,}  SIAM J.\ Math.\ Anal.\ {\bf 23} (1992), 923--940. 
\bibitem{bmo2} H. A.\ Biagioni, M.\ Oberguggenberger,
        {\it Generalized solutions to Burgers' equation,} 
        J.\ Differential Equations {\bf 97} (1992), 263--287. 
\bibitem{bc} F.\ Berger, J. F.\ Colombeau,
        {\it Numerical solutions of one-pressure models in multifluid flows,} 
        SIAM J.\ Numer.\ Anal.\ {\bf 32} (1995), 1139--1154.
\bibitem{bcm} F.\ Berger, J. F.\ Colombeau, M.\ Moussaoui,
        {\it Solutions mesures de Dirac de systemes de lois de conservation 
        et applications numeriques,} C.\ R.\ Acad.\ Sci.\ Paris S\'er. I Math.
        {\bf 316} (1993), 989--994. 
\bibitem{c0} J. F.\ Colombeau, 
        {\it Differential Calculus and Holomorphy, Real and
        Complex Analysis in Locally Convex Spaces,} North Holland, Amsterdam, 1982.
\bibitem{c1} J. F.\ Colombeau, 
        {\it New Generalized Functions and Multiplication
        of Di\-stri\-bu\-tions,} North Holland, Amsterdam, 1984.
\bibitem{c2} J. F.\ Colombeau, 
        {\it Elementary Introduction to New Generalized Functions,} 
        North Holland, Amsterdam, 1985.
\bibitem{cbull} J. F.\ Colombeau, 
        {\it Multiplication of distributions,} 
        Bull.\ Amer.\ Math.\ Soc.\ (N.S.)
        {\bf 23} (1990), 251--268.
\bibitem{clec} J. F.\ Colombeau, 
        {\it Multiplication of Distributions. A Tool in Mathematics, Numerical 
        Engineering and Theoretical Physics,} 
        Lecture Notes in Math. {\bf 1532}, Springer, New York, 1992.
\bibitem{CM} J. F.\ Colombeau, A.\ Meril, 
        {\it Generalized functions and multiplication
        of distributions on ${\mathcal C}^\infty$ manifolds,} 
        J.\ Math.\ Anal.\ Appl.\ {\bf 186} (1994), 357--364.
\bibitem{chmo} J. F.\ Colombeau, A.\ Heibig, M.\ Oberguggenberger,
        {\it Le probleme de Cauchy dans un espace de fonctions generalisees I,} 
        C.\ R.\ Acad.\ Sci.\ Paris S\'er. I Math.
        {\bf 317} (1993), 851--855.
\bibitem{chmo2} J. F.\ Colombeau, A.\ Heibig, M.\ Oberguggenberger,
        {\it Le probleme de Cauchy dans un espace de fonctions generalisees II,} 
        C.\ R.\ Acad.\ Sci.\ Paris S\'er. I Math.
        {\bf 319} (1994), 1179--1183.
\bibitem{cmo} J. F.\ Colombeau, M.\ Oberguggenberger,
        {\it On a hyperbolic system with a compatible quadratic term: 
        Generalized solutions, delta waves, and multiplication of distributions,} 
        Comm.\ Partial Differential Equations
        {\bf 15} (1990), 905--938.
\bibitem{dd} J. W.\ de Roever, M.\ Damsma, 
        {\it Colombeau algebras on a ${\mathcal C}^\infty$-manifold,} 
        Indag.\ Math.\ (N.S.) {\bf 2} (1991), 341--358.
\bibitem{tschapitsch} N.\ Dapi\'c, S.\ Pilipovi\'c, 
        {\it Microlocal analysis of Colombeau's generalized functions on a manifold,} 
        Indag.\ Math.\ (N.S.) {\bf 7} (1996), 293--309.
\bibitem{trio} N.\ Dapi\'c, S.\ Pilipovi\'c, D.\ Scarpal\'ezos,
        {\it Microlocal analysis of Colombeau's generalized functions---Propagation of
        singularities,} J.\ Anal.\ Math.\ {\bf 75} (1998), 51--66.
\bibitem{D3} J.\ Dieudonn\'e, 
        {\it \'El\'ements d'Analyse,} Vol {\bf 3}, 
        Gauthier-Villars, Paris,  1974.
\bibitem{fo} [Void; in the second part of the series this item refers to the present
        article.]
\bibitem{F} E.\ Farkas, 
        {\it Approximation properties of convenient vector spaces,} 
        Preprint, Wien, 1996.\footnote[1]{available
        electronically under {\tt
        http://diana.mat.univie.ac.at/\symbol{126}diana/dianapub.html}}
\bibitem{FK} A.\ Fr\"olicher, A.\ Kriegl, 
        {\it Linear Spaces and Differentiation Theory,} Wiley, Chichester, 1988.
\bibitem{gmot} H.\ Grosse, M.\ Oberguggenberger, I. T.\ Todorov, 
        {\it Generalized functions for quantum fields obeying 
        quadratic exchange relations,} ESI Preprint\footnote[2]{available 
        electronically under {\tt http://www.esi.ac.at/ESI-Preprints.html}}   
        {\bf 653}, (1999).
\bibitem{ft} M.\ Grosser, 
        {\it On the foundations of nonlinear generalized functions
        II,} Preprint,
        Wien, 1999.
\bibitem{vi} M.\ Grosser, M.\ Kunzinger, R.\ Steinbauer, J.\ Vickers,
        {\it  A global theory of nonlinear generalized functions,} Preprint,
        Wien, 1999.
\bibitem{Hoe} L.\ H\"ormander, 
        {\it The Analysis of Linear Partial Differential Operators I,} 
        Grund\-lehren Math.\ Wiss. {\bf 256}, Berlin
        1990.
\bibitem{JEL} J.\ Jel\'\i nek, 
        {\it An intrinsic definition of the Colombeau generalized functions,} 
        Comment.\ Math.\ Univ.\ Carolin. {\bf 40} (1999), 71--95.
\bibitem{kk} B. L.\ Keyfitz, H. C.\ Kranzer,
        {\it Spaces of weighted measures for conservation laws with 
        singular shock solutions,} 
        J.\ Differential Equations {\bf 118} (1995), 420--451.
\bibitem{KM} A.\ Kriegl, P. W.\ Michor,  
        {\it The Convenient Setting of Global 
        Analysis,}  Math.\ Surveys Monogr.\ {\bf 53}, 
        Amer.\ Math.\ Soc., Providence, RI, 1997.
\bibitem{kunzdiss} M.\ Kunzinger,  
        {\it Lie Transformation Groups in Colombeau Algebras,} 
        doctoral thesis, University of Vienna, 1996.
\bibitem{KS} M.\ Kunzinger, R.\ Steinbauer, 
        {\it A rigorous solution concept for  geodesic and geodesic deviation 
        equations in impulsive gravitational waves,} 
        J.\ Math.\ Phys.\ {\bf 40} (1999), 1479--1489.
\bibitem{landau} E.~Landau,
        {\it Einige Ungleichungen f\"ur zweimal differentiierbare
        Funktionen,}
        Proc.\ London Math.\ Soc. Ser.\ 2, {\bf 13} (1913--1914),
        43--49.
\bibitem{serbenbuch} M.\ Nedelkov, S.\ Pilipovi\'c, D.\ Scarpal\'ezos,
        {\it The Linear Theory of Colombeau Generalized Functions,} 
        Pitman Res.\ Notes Math.\ Ser.\ {\bf 385}, Longman, Harlow, 1998. 
\bibitem{MObook} M.\ Oberguggenberger, 
        {\it Multiplication of Distributions and
        Applications to Partial Differential Equations,} 
        Pitman Res.\ Notes Math.\ Ser.\ {\bf 259}, Longman, Harlow, 1992.
\bibitem{OR} M.\ Oberguggenberger, F.\ Russo, 
        {\it Nonlinear SPDEs, Colombeau solutions and
        pathwise limits,} in L.\ Decreusefond, J.\ Gjerde, B.\ \O ksendal, A. S.\ \"Ust\"unel 
        (Eds.), 
        Stochastic Analysis and Related Topics VI., Birkh\"auser, Boston, 1998, pp.\ 319--332. 
\bibitem{r1} E. E.\ Rosinger, 
        {\it Distributions and Nonlinear Partial Differential
        Equations,} Lecture Notes Math. {\bf 684}, Springer, New York, 1978.
\bibitem{r2} E. E.\ Rosinger, 
        {\it Non-Linear Partial Differential Equations.\ An
        Algebraic View of Generalized Solutions,} 
        North Holland, Amsterdam, 1990.
\bibitem{schae} H. H.\ Schaefer, 
        {\it Topological Vector Spaces }(5th ed.), 
        Grad.\ Texts in Math., Springer, 1986.
\bibitem{Schw} L.\ Schwartz, 
        {\it Sur l'impossibilite de la multiplication des distributions,} 
        C.\ R.\ Acad.\ Sci.\ Paris S\'er. I Math. 
        {\bf 239} (1954), 847--848. 
\bibitem{roldiss} R.\ Steinbauer 
        {\it Distributional Methods in General Relativity,} 
        doctoral thesis, University of Vienna, 1999.
%
\bibitem{vw1} J. A.\ Vickers, J. P.\ Wilson, 
        {\it Invariance of the distributional curvature of the 
        cone under smooth diffeomorphisms,}  
        Classical Quantum Gravity {\bf 16} (1999), 579--588. 
\bibitem{vw2} J. A.\ Vickers, J. P.\ Wilson, 
        {\it A nonlinear theory of tensor distributions,} ESI Preprint\footnote[1]
        {available electronically under {\tt http://www.esi.ac.at/ESI-Preprints.html}}  
        {\bf 566} (1998).
\bibitem{v} J. A.\ Vickers, 
        {\it Nonlinear generalised functions in general relativity,}  in 
        M.\ Grosser, G.\ H\"ormann, M.\ Kunzinger, M.\ Oberguggenberger (Eds.),
        Nonlinear Theory of Generalized Functions, Chapman \& Hall/CRC Res.
        Notes Math. {\bf 401}, Chapman \& Hall/CRC, Boca Raton, FL, 1999,
        pp.\ 275--290.
\bibitem{w}  J. P.\ Wilson, 
        {\it Distributional curvature of time dependent
        cosmic strings,} Classical Quantum Gravity {\bf 14} (1997), 2485--2498.
\bibitem{Yama} S.Yamamuro, 
        {\it Differential Calculus in Topological Linear
        Spaces,} Lecture Notes in Math. {\bf 374}, Springer, New York, 1974.


\end{thebibliography}
\end{document}